\documentclass[twoside,12pt]{article}
\usepackage[sumlimits,intlimits,namelimits]{amsmath}
\usepackage{amsthm,amssymb,amstext,amsfonts}
\usepackage[mathscr]{eucal}
\usepackage{makeidx}
\usepackage{exscale}
\newcommand{\Em}[1]{\emph{#1}} 
\newcommand{\De}[1]{\emph{#1}} 
\newcommand{\unity}{{1\!\!\!\:\mathrm{l}}}

\newcommand{\Var}{\Delta}
\DeclareMathOperator{\vol}{vol}

\DeclareMathOperator*{\ord}{ord}
\newcommand{\comp}{} 
\newcommand{\Op}[1]{\mathsf{#1}} 
\DeclareMathAlphabet{\Mat}{U}{eur}{m}{n} 
\DeclareMathAlphabet{\Set}{U}{eur}{m}{n} 
\newcommand{\Spa}[1]{\mathrm{#1}} 
\newcommand{\Sh}[1]{\mathcal{#1}} 
\newcommand{\Func}[1]{\mathcal{#1}} 
\newcommand{\lattice}{\Lambda}

\newcommand{\willmore}{\mathscr{W}}

\newcommand{\moduli}{\mathcal{M}}
\newcommand{\sobolev}[1]{W^{#1}}
\newcommand{\banach}[1]{L^{\!#1}}

\renewcommand{\qed}{{\bf\hspace{\fill}q.e.d.}}
\newtheorem{Lemma}{Lemma}[section]
\newtheorem{Theorem}[Lemma]{Theorem}
\newtheorem{Corollary}[Lemma]{Corollary}
\newtheorem{Proposition}[Lemma]{Proposition}

\newtheorem{Remark}[Lemma]{Remark}
\begin{document}
\date{}

\title{Existence of Minimizing Willmore Surfaces of Prescribed
Conformal Class}

\author{Martin U. Schmidt\\
Max--Planck--Institut f\"ur Gravitationsphysik\\
Albert--Einstein--Institut\\
Am M\"uhlenberg 1\\
D--14476 Golm\\
email: {mschmidt@aei-potsdam.mpg.de}}
\begin{titlepage}
\maketitle
\tableofcontents 
\end{titlepage}
\newpage

\section{Introduction}\label{section introduction}
The Weierstra{\ss} formula describes conformal minimal immersions of a
Riemann surface into the 3--dimensional Euclidean space in terms of a
spinor in the kernel of the Dirac operator on the Riemann
surface. Before the Dirac operators was invented, a local
generalization of this formula, which is now called Weierstra{\ss}
representation, was already known by Eisenhardt \cite{Ei}.
The global version describes a conformal immersion of a surface into the
3--dimensional Euclidean space again in terms of a spinor in the
kernel of the Dirac operator with potential on the surface
\cite{Kon2,Ta1,Ta2,Fr2}.
Pinkall and Pedit generalized this Weierstra{\ss} representation to
immersion into 4--dimensional Euclidean space and invented the
`quaternionic function theory' \cite{PP,BFLPP,FLPP}. From their point
of view conformal immersion of Riemann surfaces into the
4--dimensional Euclidean space (identified with the quaternions) are
essentially sections of holomorphic quaternionic line bundles.
These holomorphic quaternionic line bundles are build form an usual
holomorphic complex line bundle on the Riemann surface together with a
Hopf field. Due to an observation of Taimanov,
the Willmore functional is equal to four times the integral
over the square of the potential \cite{Ta1}.
Our main subject is the investigation of these
holomorphic quaternionic line bundles, whose Hopf fields are
square integrable. The holomorphic sections of such quaternionic
line bundles form the maximal domain of definition
of the Willmore functional on the space of conformal mappings of a
Riemann surface into $\mathbb{H}\simeq\mathbb{R}^4$.
In the second section we extend Cauchy's integral formula to
these holomorphic section of quaternionic holomorphic line bundles.
In the fourth section we show that the corresponding
sections define sheaves, and that the \u{C}ech cohomology groups of these
sheaves obey the Riemann--Roch Theorem and S\'{e}rre duality.
In the fifth section extend those B\"acklund transformations to
square integrable Hopf fields, which relate the infinitesimal
quaternionic Weierstra{\ss} representation to the Kodaira embedding of
`quaternionic function theory' \cite{PP}. This yields in the sixth
section a general proof of the Pl\"ucker formula \cite{FLPP}
for these holomorphic quaternionic line bundles with square integrable
Hopf fields. In the seventh section we show that any bounded sequence
of square integrable Hopf fields has a convergent subsequence, and that
the limit is again the Hopf field of a
holomorphic quaternionic line bundle, but the holomorphic structure
might have singularities. This is used in the last section to proof
that the Willmore functional has on the space of all conformal
immersions of a compact Riemann surface into the 3--dimensional and
4--dimensional Euclidean space a minimum. Moreover, even the
restrictions of the Willmore functional to all conformal immersions
into the 4--dimensional Euclidean space has a minimum, whose underlying
holomorphic complex line bundle (compare \cite{PP,FLPP}) is fixed.
The existence of a minimizer on the space of all immersions from a
Riemann surface of prescribed genus into the
$n$----dimensional Euclidean spaces was proven by Simon for genus one
\cite{Si1,Si2}, and recently by Bauer and Kuwert for all finite genera
\cite{BK}.

We identify the quaternions with all complex $2\times 2$--matrices of
the form $\left(\begin{smallmatrix}
a & b\\
-\Bar{b} & \Bar{a}
\end{smallmatrix}\right)$. If we consider a $\mathbb{C}^2$--valued
function $\psi=\left(\begin{smallmatrix}
\psi_1\\
\psi_2
\end{smallmatrix}\right)$ on an open set
$\Omega\subset\mathbb{C}$ as a quaternionic valued function
$\left(\begin{smallmatrix}
\psi_1 & -\Bar{\psi}_2\\
\psi_2 & \Bar{\psi}_1
\end{smallmatrix}\right)$, then the operator $\left(\begin{smallmatrix}
\Bar{\partial} & -\Bar{U}\\
U & \partial
\end{smallmatrix}\right)$ defines a quaternionic holomorphic structure
in the sense of \cite[Definition~2.1.]{FLPP} on the trivial
quaternionic line bundle on $\Omega$ endowed with the complex
structure of multiplication on the left with complex numbers
$\mathbb{C}\subset\mathbb{Q}$. In particular, the action of
$\sqrt{-1}$ is given by left--multiplication with
$\left(\begin{smallmatrix}
\sqrt{-1} & 0\\
0 & -\sqrt{-1}
\end{smallmatrix}\right)$.
The corresponding holomorphic sections are defined as the elements
of the kernel of this operator, which agrees with the elements
of the kernel of the Dirac operator
\begin{equation*}
\begin{pmatrix}
U & \partial\\
-\Bar{\partial} & \Bar{U}
\end{pmatrix} =\begin{pmatrix}
0 & \unity\\
-\unity & 0
\end{pmatrix}\comp\begin{pmatrix}
\Bar{\partial} & -\Bar{U}\\
U & \partial
\end{pmatrix}.
\end{equation*}
The corresponding Hopf field is equal to $Q=-\Bar{U}d\Bar{z}$.
The space of holomorphic sections is invariant under
right--multiplication with quaternions and therefore
a quaternionic vector space.

The holomorphic structure is an operator from the sections of a
quaternionic line bundle into the space of sections of this
quaternionic line bundle tensored with the line bundle of
anti--holomorphic forms ($\simeq\Sh{O}_{-K}$) \cite[\S2.2]{FLPP}.
We represent the underlying holomorphic complex line bundle
on a Riemann surface $\Spa{X}$ as the trivial complex line bundles
on all members of an open covering together with a cocycle in the
corresponding multiplicative first \u{C}ech Cocomplex,
which represents an element in $H^1(\Spa{X},\Sh{O}^{\ast})$.
We shall state how the holomorphic structure transforms
under these cocycles and coordinate transformations
$z\mapsto z'=z'(z)$.
The multiplication with a non-vanishing function $f$ acts on the
spinors as $\psi\mapsto\left(\begin{smallmatrix}
f & 0\\
0 & \Bar{f}
\end{smallmatrix}\right)\psi$.
Therefore this multiplicative cocycle acts on the holomorphic structure as
$$\begin{pmatrix}
\Bar{\partial} & -\Bar{U}\\
U & \partial\\
\end{pmatrix}\mapsto
\begin{pmatrix}
f & 0\\
0 & \Bar{f}
\end{pmatrix}\comp\begin{pmatrix}
\Bar{\partial} & -\Bar{U}\\
U & \partial
\end{pmatrix}\comp\begin{pmatrix}
f & 0\\
0 & \Bar{f}
\end{pmatrix}^{-1}=
\begin{pmatrix}
-\Bar{\partial} & -\frac{f}{\Bar{f}}\Bar{U}\\
\frac{\Bar{f}}{f}U & \partial
\end{pmatrix}.$$
The corresponding potential $U$ and Hopf field $Q$ transforms as
$U\mapsto\frac{\Bar{f}}{f}U$ and $Q\mapsto\frac{f}{\Bar{f}}Q$.
The coordinate transformation $z\mapsto z'=z'(z)$
acts on the holomorphic structure as
$$\begin{pmatrix}
\Bar{\partial} & -\Bar{U}\\
U & \partial
\end{pmatrix}\mapsto
\begin{pmatrix}
\Bar{\partial}' & -\Bar{U}'\\
U' & \partial'
\end{pmatrix}=
\begin{pmatrix}
\overline{\frac{dz'}{dz}} & 0\\
0 & \frac{dz'}{dz}
\end{pmatrix}^{-1}\begin{pmatrix}
\Bar{\partial} & -\Bar{U}\\
U & \partial
\end{pmatrix}=
\begin{pmatrix}
\Bar{\partial}' & -\overline{\frac{dz}{dz'}}\Bar{U}\\
\frac{dz}{dz'}U & \partial'
\end{pmatrix}.$$
Therefore the potentials transforms as
$U\mapsto U'=\frac{dz}{dz'}U$ and the corresponding
Hopf field $Q=-\Bar{U}d\Bar{z}=-\Bar{U}'d\Bar{z}'$ does
not change.
Summing up, for any holomorphic complex line bundle,
which is represented by the trivial line bundles on all members of
an open covering together with a cocycle in $H^1(\Spa{X},\Sh{O}^{\ast})$,
this cocycle defines also cocycles for the corresponding spinors
$\psi$ and potentials $U$.

\newtheorem{Weierstrass}[Lemma]{Quaternionic Weierstra{\ss} Representation}
\begin{Weierstrass}\label{weierstrass}
\cite[Theorem~4.3.]{PP}
For any conformal immersion $f:\Spa{X}\rightarrow\mathbb{H}$
of a Riemann surface $\Spa{X}$ there
exist two quaternionic holomorphic line bundles
with two holomorphic sections $\psi$ and $\phi$,
such that the derivative of $f$ is given by
$$d\begin{pmatrix}
f_1 & -\Bar{f}_2\\
f_2 & \Bar{f}_1
\end{pmatrix}=\begin{pmatrix}
\phi_1 & \phi_2\\
-\Bar{\phi}_2 & \Bar{\phi}_2
\end{pmatrix}\begin{pmatrix}
dz & 0\\
0 & d\Bar{z}
\end{pmatrix}\begin{pmatrix}
\psi_1 & -\Bar{\psi}_2\\
\psi_2 & \Bar{\psi}_1
\end{pmatrix}\text{ with}$$
\begin{align*}
\begin{pmatrix}
\Bar{\partial} & -\Bar{U}\\
U & \partial
\end{pmatrix}\begin{pmatrix}
\psi_1 & -\Bar{\psi}_2\\
\psi_2 & \Bar{\psi}_1
\end{pmatrix}&=0&
\begin{pmatrix}
\Bar{\partial} & U\\
-\Bar{U} & \partial
\end{pmatrix}\begin{pmatrix}
\phi_1 & -\Bar{\phi}_2\\
\phi_2 & \Bar{\phi}_1
\end{pmatrix}&=0.
\end{align*}\qed
\end{Weierstrass}

We remark that the product of the underlying complex line bundles has
to be equal to the anti--canonical line bundle
(i.\ e.\ the line bundle of anti--holomorphic forms)
and that the potentials of both holomorphic
structures are determined by each other.
Immersion into $\mathbb{R}^3$ are obtained as immersion into the pure
imaginary quaternions $\simeq\mathbb{R}^3$. This is realized by the
additional reality conditions $\Bar{U}=U$, $df^{\ast}=-df$
\begin{align*}
\begin{pmatrix}
\phi_1 & \phi_2\\
-\Bar{\phi}_2 & \Bar{\phi}_1
\end{pmatrix}^{\ast}&=\begin{pmatrix}
0 & \sqrt{-1}\\
\sqrt{-1} & 0
\end{pmatrix}\begin{pmatrix}
\psi_1 & -\Bar{\psi}_2\\
\psi_2 & \Bar{\psi}_1
\end{pmatrix}&
\begin{pmatrix}
\psi_1 & -\Bar{\psi}_2\\
\psi_2 & \Bar{\psi}_1
\end{pmatrix}^{\ast}&=
\begin{pmatrix}
\phi_1 & \phi_2\\
-\Bar{\phi}_2 & \Bar{\phi}_1
\end{pmatrix}\begin{pmatrix}
0 & \sqrt{-1}\\
\sqrt{-1} & 0
\end{pmatrix}.
\end{align*}

\section{Local behaviour of holomorphic spinors}\label{section local}

Dolbeault's Lemma \cite[Chapter~I Section~D 2.~Lemma]{GuRo} implies
that the operator  $\Op{I}_{\mathbb{C}}(0)$ with the integral kernel
$$\begin{pmatrix}
(z-z')^{-1} & 0\\
0 & (\Bar{z}-\Bar{z}')^{-1}
\end{pmatrix}\frac{d\Bar{z}'\wedge dz'}{2\pi\sqrt{-1}}$$
is a right inverse of the operator
$\left(\begin{smallmatrix}
\Bar{\partial} & 0\\
0 & \partial
\end{smallmatrix}\right).$
Due to the
Hardy--Littlewood--Sobolev theorem \cite[Chapter~V. \S1.2 Theorem~1]{St}
for all $1<p<2$ and $2<q<\infty$ with $\frac{1}{p}=\frac{1}{q}+\frac{1}{2}$
this is a bounded operator from
$\banach{p}(\mathbb{C},\mathbb{H})$
into $\banach{q}(\mathbb{C},\mathbb{H})$.
Moreover, the restriction $\Op{I}_{\Omega}(0)$ of
$\Op{I}_{\mathbb{C}}(0)$ to a bounded open domain $\Omega$ is a
bounded operator from $\banach{p}(\Omega,\mathbb{H})$
into $\banach{q}(\Omega,\mathbb{H})$.
On the other hand H\"older's inequality \cite[Theorem~III.1~(c)]{RS1}
implies that the multiplication operators with
$U\in\banach{2}(\Omega)$ are
bounded operators from $\banach{q}(\Omega)$ into
$\banach{p}(\Omega)$.
Hence the operator
$\unity+\Op{I}_{\Omega}(0)\comp\left(\begin{smallmatrix}
0 & -\Bar{U}\\
U & 0
\end{smallmatrix}\right)$
is a bounded operator on
$\banach{q}(\Omega,\mathbb{H})$.
For smooth $U$ a spinor $\psi$ belongs to the kernel of
$\left(\begin{smallmatrix}
\Bar{\partial} & -\Bar{U}\\
U & \partial
\end{smallmatrix}\right)$, if and only if
$\unity+\Op{I}_{\Omega}(0)\comp\left(\begin{smallmatrix}
0 & -\Bar{U}\\
U & 0
\end{smallmatrix}\right)$
maps $\psi$ into the kernel of
$\left(\begin{smallmatrix}
\Bar{\partial} & 0\\
0 & \partial
\end{smallmatrix}\right)$.
Due to Weyl's Lemma \cite[Theorem~IX.25]{RS2} all elements
in the kernel of this differential operator are smooth functions.
Consequently, for all
$U\in\banach{2}_{\text{\scriptsize\rm loc}}(\Omega)$ the kernel of
$\left(\begin{smallmatrix}
\Bar{\partial} & -\Bar{U}\\
U & \partial
\end{smallmatrix}\right)$ is defined as all spinors
$\psi\in\banach{q}_{\text{\scriptsize\rm loc}}(\Omega,\mathbb{H})$,
which are mapped by
$\unity+\Op{I}_{\Omega}(0)\comp\left(\begin{smallmatrix}
0 & -\Bar{U}\\
U & 0
\end{smallmatrix}\right)$
into the kernel of
$\left(\begin{smallmatrix}
\Bar{\partial} & 0\\
0 & \partial
\end{smallmatrix}\right)$.
Finally, we remark that we may always cover $\Omega$ by small
sets $\Omega'$ such that the von Neumann series
$$\left(\unity+\Op{I}_{\Omega}(0)\comp\begin{pmatrix}
0 & -\Bar{U}\\
U & 0
\end{pmatrix}\right)^{-1}=
\sum\limits_{l=0}^{\infty}
\left(\Op{I}_{\Omega}(0)\comp\begin{pmatrix}
0 & \Bar{U}\\
-U & 0
\end{pmatrix}\right)^{l}$$
converges as an operator on
$\banach{q}(\Omega',\mathbb{H})$, which maps the
closed subspace of bounded elements in the kernel of
$\left(\begin{smallmatrix}
\Bar{\partial} & 0\\
0 & \partial
\end{smallmatrix}\right)$ onto the closed subspace of bounded
elements in the kernel of $\left(\begin{smallmatrix}
\Bar{\partial} & -\Bar{U}\\
U & \partial
\end{smallmatrix}\right)$. Therefore the latter kernel is contained in
$\bigcap\limits_{1<p<2}
\sobolev{1,p}_{\text{\scriptsize\rm loc}}(\Omega,\mathbb{H})
\subset\bigcap\limits_{q<\infty}
\banach{q}_{\text{\scriptsize\rm loc}}(\Omega,\mathbb{H})$.

Moreover, on small domains $\Omega\subset\mathbb{C}$ the operator
$$\Op{I}_{\Omega}(U)=\Op{I}_{\Omega}(0)\comp\left(\unity+\begin{pmatrix}
0 & -\Bar{U}\\
U & 0
\end{pmatrix}\comp\Op{I}_{\Omega}(0)\right)^{-1}=
\left(\unity+\Op{I}_{\Omega}(0)\comp\begin{pmatrix}
0 & -\Bar{U}\\
U & 0
\end{pmatrix}\right)^{-1}\comp\Op{I}_{\Omega}(0)$$
is a right inverse of the operator $\left(\begin{smallmatrix}
\Bar{\partial} & -\Bar{U}\\
U & \partial
\end{smallmatrix}\right)$.
If $\Func{K}_{\Omega}(U,z,z')\frac{d\Bar{z}\wedge dz}{2\pi\sqrt{-1}}$
denotes the integral kernel of this operator $\Op{I}_{\Omega}(U)$,
then we have
\begin{align*}
\begin{pmatrix}
\Bar{\partial} & -\Bar{U}\\
U & \partial
\end{pmatrix}\Func{K}_{\Omega}(U,z,z')&=
\pi\delta(z-z')\unity &
\begin{pmatrix}
-\Bar{\partial} & U\\
-\Bar{U} & -\partial\\
\end{pmatrix}\Func{K}_{\Omega}^{t}(U,z',z)&=
\pi\delta(z-z')\unity.
\end{align*}
Here the differential operator and his transposed acts on the integral kernel
as a function depending on $z$ for fixed $z'$.
If $\psi$ is an element of the kernel of $\left(\begin{smallmatrix}
\Bar{\partial} & -\Bar{U}\\
U & \partial
\end{smallmatrix}\right)$ and $\phi$ an element of the kernel of
$\left(\begin{smallmatrix}
-\Bar{\partial} & U\\
-\Bar{U} & -\partial\\
\end{smallmatrix}\right)$, then a direct calculation shows
\begin{eqnarray*}
d\left(\Func{K}_{\Omega}(U,z',z)
\begin{pmatrix}
dz & 0\\
0 & -d\Bar{z}
\end{pmatrix}\psi(z)\right)&=&\pi\delta(z-z')\psi(z)d\Bar{z}\wedge dz\\
d\left(\phi(z)\begin{pmatrix}
dz & 0\\
0 & -d\Bar{z}
\end{pmatrix}\Func{K}_{\Omega}(U,z,z')\right)&=&
\pi\delta(z-z')\phi(z)dz\wedge d\Bar{z}.
\end{eqnarray*}
This implies a quaternionic version of

\newtheorem{Cauchy}[Lemma]{Cauchy's Integral Formula}
\index{Cauchy's integral formula}
\begin{Cauchy}\label{cauchy formula}
All elements $\psi$ and $\phi$ in the kernel of $\left(\begin{smallmatrix}
\Bar{\partial} & -\Bar{U}\\
U & \partial
\end{smallmatrix}\right)$ and $\left(\begin{smallmatrix}
-\Bar{\partial} & U\\
-\Bar{U} & -\partial
\end{smallmatrix}\right)$  on a small open set $\Omega$
obey the formula
\begin{equation*}\begin{split}
\psi(z')&=\frac{1}{2\pi\sqrt{-1}}\oint\Func{K}_{\Omega}(U,z',z)
\begin{pmatrix}
dz & 0\\
0 & -d\Bar{z}
\end{pmatrix}\psi(z)\\
\phi(z')&=\frac{-1}{2\pi\sqrt{-1}}\oint\phi(z)\begin{pmatrix}
dz & 0\\
0 & -d\Bar{z}
\end{pmatrix}\Func{K}_{\Omega}(U,z,z'),
\end{split}\end{equation*}
as long as the integration path surrounds $z'$
one times in the anti--clockwise--order, respectively.
\end{Cauchy}

\begin{Remark}\label{smeared closed path}
At a first look it is not clear, whether the integral along the closed
path is well defined. However, since on the complement of $\{z'\}$ the
corresponding one--forms are closed, we may extend the integration
over the closed path to an integration over a cylinder around $z'$.
More precisely, let $f$ be a quaternionic smooth function
with compact support in $\Omega$,
which is equal to $\unity$ on an open subset $\Omega'$
with $\Bar{\Omega}'\subset\Omega$. Then we have the following
equality of measurable functions on $z'\in\Omega':$
\begin{equation*}\begin{split}
\psi(z')&=\frac{1}{2\pi\sqrt{-1}}
\int\limits_{\Omega}df\wedge\Func{K}_{\Omega}(U,z',z)
\begin{pmatrix}
dz & 0\\
0 & -d\Bar{z}
\end{pmatrix}\psi(z)\\
\phi(z')&=\frac{-1}{2\pi\sqrt{-1}}\int\limits_{\Omega}\phi(z)\begin{pmatrix}
dz & 0\\
0 & -d\Bar{z}
\end{pmatrix}\Func{K}_{\Omega}(U,z,z')\wedge df.
\end{split}\end{equation*}
\end{Remark}

In particular, for all square integrable Hopf fields
(i.\ e.\ all potentials belong to $\banach{2}_{\text{\scriptsize\rm loc}}$)
the holomorphic sections of the corresponding
holomorphic quaternionic line bundle define a sheaf
on $\Spa{X}$. In the sequel we shall denote by $\Sh{Q}_{D}$
the sheaf of sections of a holomorphic quaternionic line bundle
over the complex line bundle corresponding to $\Sh{O}_{D}$.

\begin{Lemma}\label{integral kernel bound}
For small potentials $U\in\banach{2}(\Omega)$ on a bounded domain
$\Omega$ there exists positive functions
$\Func{A},\Func{B}\in\bigcap\limits_{q<\infty}\banach{q}(\Omega)$
such that the integral kernels
$\Func{K}_{\Omega}(U,z,z')$ may be estimated by
$$\left|\Func{K}_{\Omega}(U,z,z')\right|\leq
\frac{\Func{A}(z)\Func{B}(z')}{|z-z'|}.$$
\end{Lemma}

\begin{proof}
The equation
$\displaystyle{\frac{1}{(z-z')(z'-z'')}+
\frac{1}{(z'-z'')(z''-z)}+\frac{1}{(z''-z)(z-z')}=0}$
implies the estimate
$$\frac{1}{\left|z-z'\right|\left|z'-z''\right|}\leq
\frac{1}{\left|z'-z''\right|\left|z''-z\right|}+
\frac{1}{\left|z''-z\right|\left|z-z'\right|}.$$
Therefore the integral kernel of the operator
$\Op{I}_{\Omega}(0)\comp\left(\begin{smallmatrix}
0 & -\Bar{U}\\
U & 0
\end{smallmatrix}\right)\comp\Op{I}_{\Omega}(0)$ is bounded by
$\displaystyle{\frac{\Func{F}(z)+\Func{F}(z')}{\pi^2|z-z'|}}$,
where $\Func{F}(z)$ is the convolution of the
$\banach{2}$--function $|U|$ with the positive function
$\displaystyle{\frac{1}{|z|}}$. Due to
Young's inequality \cite[Section~IX.4 Example~1]{RS2}
this function $\Func{F}$ belongs to
$\bigcap\limits_{q<\infty}\banach{q}(\Omega)$.
An iterative application of this argument to all terms of the
von Neumann series of $\Op{I}_{\Omega}(U)$ yields a bound of the
integral kernel of the form
$\displaystyle{\sum\limits_l
\frac{\Func{A}_l(z)\Func{B}_l(z')}{|z-z'|}}$.
For small $\banach{2}$--norms of $U$ all $\banach{q}(\Omega)$--norms
of $\displaystyle{\sum\limits_l\Func{A}_l}$ and
$\displaystyle{\sum\limits_l\Func{B}_l}$ are bounded.
This completes the proof.
\end{proof}

The holomorphic sections of a holomorphic quaternionic line bundle are
in general not continuous. Nevertheless they share many properties
with the sheaves of holomorphic functions. In particular, they have the
\De{Strong unique continuation property}, which is proven by a
\De{Carleman inequality} (compare with \cite{Ca} and
\cite[Proposition~1.3]{Wo}).

\newtheorem{Carleman inequality}[Lemma]{Carleman inequality}
\begin{Carleman inequality}\label{carleman inequality}
There exists some constant $S_p$
depending only on $1<p<2$, such that for all $n\in\mathbb{Z}$
and all $\psi\in C^{\infty}_0(\mathbb{C}\setminus\{0\},\mathbb{H})$
the following inequality holds:
$$\left\| |z|^{-n}\psi\right\|_{\frac{2p}{2-p}}\leq S_p
\left\| |z|^{-n}\left(\begin{smallmatrix}
\Bar{\partial} & 0\\
0 & \partial
\end{smallmatrix}\right)\psi\right\|_{p}.$$
\end{Carleman inequality}

The literature \cite{Je,Ki1,Ma,Ki2} deals with the much more difficult
higher--dimensional case and does not treat our case.
David Jerison pointed out to the author, that the arguments of
\cite[Proposition~2.6]{Wo}, where the analogous but weaker statement
about the gradient term of the Laplace operator is treated,
carry over to the Dirac operator.

\begin{proof}
Dolbeault's Lemma \cite[Chapter~I Section~D 2.~Lemma]{GuRo}
implies for all smooth $\psi$ with compact support the equality
$$\psi(z)=\int\limits_{\mathbb{C}}
\left(\begin{smallmatrix}
(z-z')^{-1} & 0\\
0 & (\Bar{z}-\Bar{z}')^{-1}
\end{smallmatrix}\right)\comp\left(\begin{smallmatrix}
\Bar{\partial} & 0\\
0 & \partial
\end{smallmatrix}\right)\psi(z')
\frac{d\Bar{z}'\wedge dz'}{2\pi\sqrt{-1}}.$$
In fact, the components of the difference of the left hand side minus
the right hand side are holomorphic and anti--holomorphic functions on
$\mathbb{C}$, respectively, which vanish at $z=\infty$.
In particular, the integrals
$\int_\mathbb{C} z^n\Bar{\partial}\psi_1d\Bar{z}\wedge dz$ and
$\int_\mathbb{C}\Bar{z}^n\partial\psi_2d\Bar{z}\wedge dz$
with $n\in\mathbb{N}_0$ are proportional to the Taylor coefficients
of $\psi$ at $\infty$, which vanish.
Moreover, if the support of $\psi$
does not contain $0$ and therefore also a small neighbourhood of $0$,
then the integrals
$\int_\mathbb{C}z^{-n}\Bar{\partial}\psi_1d\Bar{z}\wedge dz$ and
$\int_\mathbb{C}\Bar{z}^{-n}\partial\psi_2d\Bar{z}\wedge dz$ with
$n\in\mathbb{N}$ are proportional to the Taylor coefficients of
$\psi$ at $0$, which in this case also
vanish. These cancellations follow also from partial integration.
We conclude that for all $n\in\mathbb{Z}$
$$\psi(z)=\int\limits_{\mathbb{C}}
\begin{pmatrix}
\left(\frac{z}{z'}\right)^n\frac{1}{z-z'} & 0\\
0 & \left(\frac{\Bar{z}}{\Bar{z}'}\right)^n\frac{1}{\Bar{z}-\Bar{z}'}\\
\end{pmatrix}\comp\left(\begin{smallmatrix}
\Bar{\partial} & 0\\
0 & \partial
\end{smallmatrix}\right)\psi(z')
\frac{d\Bar{z}'\wedge dz'}{2\pi\sqrt{-1}}.$$
In fact, for negative $n$ the left hand side
minus the left hand side
of the foregoing formula is equal to the Taylor polynomial of $\psi$
at $\infty$ up to order $|n|$, and for positive $n$ equal to the
Taylor polynomial of $\psi$ at $0$ up to order $n-1$.
Finally, the Hardy--Littlewood--Sobolev theorem
\cite[Chapter~V. \S1.2 Theorem~1]{St} implies that the operator with
integral kernel
$$\begin{pmatrix}
\left(\frac{|z'|z}{|z|z'}\right)^n\frac{1}{z-z'} & 0\\
0 & \left(\frac{|z'|\Bar{z}}{|z|\Bar{z}'}\right)^n\frac{1}{\Bar{z}-\Bar{z}'}
\end{pmatrix}\frac{d\Bar{z}'\wedge dz'}{2\pi\sqrt{-1}}$$
from $\banach{p}(\mathbb{C},\mathbb{H})$ into
$\banach{\frac{2p}{2-p}}(\mathbb{C},\mathbb{H})$
is bounded by some constant $S_p$ not depending on $n$, and maps
$|z|^{-n}\left(\begin{smallmatrix}
\Bar{\partial} & 0\\
0 & \partial
\end{smallmatrix}\right)\psi$ onto $|z|^{-n}\psi$.
\end{proof}

Due to a standard argument
(e.\ g.\ \cite[Proof of Theorem~5.1.4]{So} and
\cite[Section Carleman Method]{Wo})
this \De{Carleman inequality} implies the

\newtheorem{Unique continuation}[Lemma]{Strong unique continuation property}
\begin{Unique continuation}\label{strong unique continuation}
Let $U$ be a potential in
$\banach{2}_{\text{\scriptsize\rm loc}}(\Omega)$ on an
open connected set $0\ni\Omega\subset\mathbb{C}$ and
$\psi\in \sobolev{1,p}_{\text{\scriptsize\rm loc}}(\Omega,\mathbb{H})$
an element of the kernel of $\left(\begin{smallmatrix}
\Bar{\partial} & -\Bar{U}\\
U & \partial
\end{smallmatrix}\right)$ on $\Omega$ with $1<p<2$.
If the $\banach{\frac{2p}{2-p}}$--norm of the restriction of
$\psi$ to the balls $B(0,\varepsilon)$ converges in the limit
$\varepsilon\downarrow 0$ faster to zero than any power of $\varepsilon$:
$$\left(\int\limits_{B(0,\varepsilon)}
\left|\psi\right|^{\frac{2p}{2-p}}d^2x
\right)^{\frac{2-p}{2p}}\leq
\text{\bf{O}}(\varepsilon^n)\;\forall n\in\mathbb{N},$$
then $\psi$ vanishes identically on $\Omega$.
\end{Unique continuation}

\begin{proof}
The question is local so we may assume that $U$ is an element
of $\banach{2}$ rather than $\banach{2}_{\text{\scriptsize\rm loc}}$.
We fix $\varepsilon$ small enough that
$\left\|U\right\|_{\banach{2}(B(z,2\varepsilon))}\leq 1/(2S_p)$
for all $z$, where $S_p$ is the constant of the
\De{Carleman inequality}.
Let $\phi\in C^{\infty}$ be $1$ on $B(0,\varepsilon)$ and $0$ on
$\mathbb{C}\setminus B(0,2\varepsilon)$. A limiting argument using the
infinite order vanishing of $\psi$ and the equality
$\left(\begin{smallmatrix}
\Bar{\partial} & 0\\
0 & \partial
\end{smallmatrix}\right)\psi=\left(\begin{smallmatrix}
0 & \Bar{U}\\
-U & 0
\end{smallmatrix}\right)\psi$ shows that the proof of the
\De{Carleman inequality} is also true for $\phi\psi$. So
$$\left\| |z|^{-n}\phi\psi\right\|_{\frac{2p}{2-p}}\leq
S_p\left\| |z|^{-n}\left(\begin{smallmatrix}
\Bar{\partial} & 0\\
0 & \partial
\end{smallmatrix}\right)\phi\psi\right\|_p\leq
S_p\left\| |z|^{-n}\phi\left(\begin{smallmatrix}
\Bar{\partial} & 0\\
0 & \partial
\end{smallmatrix}\right)\psi\right\|_p+
S_p\left\| |z|^{-n}E\right\|_p.$$
Here $E$ (for error) $=\left(\begin{smallmatrix}
\psi_1\Bar{\partial}\phi_1\\
-\psi_2\partial\phi_2
\end{smallmatrix}\right)$ is an $\banach{p}$
function supported in
$\{x\mid \varepsilon\leq|z|\leq 2\varepsilon\}$. Using the equality
$\left(\begin{smallmatrix}
\Bar{\partial} & 0\\
0 & \partial
\end{smallmatrix}\right)\psi=\left(\begin{smallmatrix}
0 & \Bar{U}\\
-U & 0
\end{smallmatrix}\right)\psi$ and
H\"older's inequality \cite[Theorem~III.1~(c)]{RS1} yields
\begin{eqnarray*}
\left\| |z|^{-n}\phi\psi\right\|_{\frac{2p}{2-p}}&\leq&
S_p\left\| |z|^{-n}\phi\left(\begin{smallmatrix}
0 & \Bar{U}\\
-U & 0
\end{smallmatrix}\right)\psi\right\|_p+S_p\left\| |z|^{-n}E\right\|_p\\
&\leq& S_p\left\|U\right\|_{\banach{2}(B(0,2\varepsilon))}
\left\| |z|^{-n}\phi\psi\right\|_{\frac{2p}{2-p}}+
S_p\left\| |z|^{-n}E\right\|_p.
\end{eqnarray*}
By the choice of $\varepsilon$ the first term can be absorbed
into a factor $2$
$$\left\| |z|^{-n}\phi\psi\right\|_{\frac{2p}{2-p}}\leq
2S_p\left\| |z|^{-n}E\right\|_p.$$
Now comes the crucial observation: $E$ is supported in
$\{z\mid \varepsilon\leq|z|\leq 2\varepsilon\}$, so
\begin{align*}
\left\| |z|^{-n}\phi\psi\right\|_{\frac{2p}{2-p}}&
\leq 2S_p\varepsilon^{-n}\left\|E\right\|_p
&\text{ and }\left\|\left(\frac{\varepsilon}{|z|}
\right)^n\phi\psi\right\|_{\frac{2p}{2-p}}
&\leq 2S_p\left\|E\right\|_p.
\end{align*}
Using the limit $n\rightarrow\infty$ we conclude that $\phi\psi$
vanishes on $B(0,\varepsilon)$, and therefore also $\psi$.
In other words, the set
$\{z\mid\psi\text{ vanishes to infinite order at }z\}$
is open, and in fact contains a ball of fixed radius $\varepsilon$
centered at any of its points.
So this set must be all of $\Omega$ and the proof is complete.
\end{proof}

The definition of the order of a zero extends
from the complex case to the quaternionic case.

\newtheorem{Order of zeroes}[Lemma]{Order of zeroes}
\begin{Order of zeroes}\label{order of zeroes}
The order of a zero of $\psi$ in the kernel of
$\left(\begin{smallmatrix}
\Bar{\partial} & -\Bar{U}\\
U & \partial
\end{smallmatrix}\right)$ on an open neighbourhood
$0\in\Omega\subset\mathbb{C}$ at $z=0$ is defined
as the largest integer $m$, such that
$\left(\begin{smallmatrix}
z & 0\\
0 & \Bar{z}
\end{smallmatrix}\right)^{-m}\psi\in
\banach{q}_{\text{\scriptsize\rm loc}}(\Omega,\mathbb{H})$
with $2<q<\infty$ belongs to the kernel of
$\left(\begin{smallmatrix}
\Bar{\partial} & -\Bar{U}\left(\frac{\Bar{z}}{z}\right)^{m}\\
U\left(\frac{z}{\Bar{z}}\right)^m & \partial
\end{smallmatrix}\right)$.
Due to the
\De{Strong unique continuation property}~\ref{strong unique continuation}
this number is finite and denoted by $\ord\limits_{0} \psi$.
\end{Order of zeroes}

On small open domains $\Omega$ with small $\banach{2}$--norms of the
potential $U$ the elements of the kernel of
$\left(\begin{smallmatrix}
\Bar{\partial} & -\Bar{U}\\
U & \partial
\end{smallmatrix}\right)$ are small perturbations of the elements of
the kernel of $\left(\begin{smallmatrix}
\Bar{\partial} & 0\\
0 & \partial
\end{smallmatrix}\right)$. We conclude that the quotient of
the space of holomorphic spinors divided by the subspace of
holomorphic spinors vanishing at $z_0$ is a one--dimensional
quaternionic vector space. Hence all holomorphic spinors
are non--vanishing sections
of another holomorphic quaternionic line bundle.
Moreover, the zeroes of such $\psi$ are isolated.
Furthermore, the divisor of an holomorphic section
of a holomorphic quaternionic line bundle is well defined.
Furthermore, for any divisor $D$ and any sheaf $\Sh{Q}_{D}$ of
holomorphic sections of a holomorphic quaternionic line bundle
those holomorphic sections, whose divisors are larger than $-D$,
define the sheaf of holomorphic sections of
another holomorphic quaternionic line bundle,
which is denoted by $\Sh{Q}_{D}$.
In particular, any holomorphic section
of a holomorphic quaternionic line bundle is the non--vanishing
holomorphic section of another holomorphic quaternionic line bundle.

For an effective divisor $D$ (i.\ e.\ $D\geq 0$) the quotient sheaf
$\Sh{Q}_{D}/\Sh{Q}$ has the same support as the divisor $D$.
For the divisor of the function $z\mapsto z^l$ on
$0\ni\Omega\subset\mathbb{C}$ this
quotient is isomorphic to the codimension of the image of the operator
$$\left(\unity+\begin{pmatrix}
0 & -\left(\frac{\Bar{z}}{z}\right)^l\Bar{U}\\
\left(\frac{z}{\Bar{z}}\right)^lU
\end{pmatrix}\comp\Op{I}_{\Omega}(0)\right)
\begin{pmatrix}
z & 0\\
0 & \Bar{z}
\end{pmatrix}^l\comp\left(\unity+\begin{pmatrix}
0 & -\Bar{U}\\
U & 0
\end{pmatrix}\comp\Op{I}_{\Omega}(0)\right)^{-1}$$
considered as an operator on the kernel of the free Dirac operator.
For small $\Omega$ with small $\banach{2}(\Omega)$--norms of $U$,
this operator is a small perturbation of $\left(\begin{smallmatrix}
z & 0\\
0 & \Bar{z}
\end{smallmatrix}\right)^l$. This proves

\begin{Lemma}\label{quotient dimension}
For any pair of divisors $D'\geq D$ on a Riemann surface $\Spa{X}$ and any
sheaf $\Sh{Q}_{D}$ of holomorphic sections of
a holomorphic quaternionic line bundle the quaternionic dimension of
$H^0\left(\Spa{X},\Sh{Q}_{D'}/\Sh{Q}_{D}\right)$ is equal to $\deg(D'-D)$
and $H^1\left(\Spa{X},\Sh{Q}_{D'}/\Sh{Q}_{D}\right)$ is trivial.\qed
\end{Lemma}

\section{Spectral theory of Dirac operators}\label{section spectral theory}
The holomorphic structures of quaternionic line bundles may be
described by first order differential operators similar to
Dirac operators. 
In this section we develop in six steps the spectral theory of
Dirac operators on compact Riemann surfaces.
Due to the Sobolev Embedding \cite[5.4 Theorem]{Ad}, these Dirac
operators are for all $1<p<2$ bounded operators
from the Sobolev spaces of $\sobolev{1,p}$--spinors into the
$\banach{p}$--spinors. Furthermore, their resolvents turn out to be
bounded operators from the $\banach{p}$--spinors onto the
$\sobolev{1,p}$--spinors. Hence we shall define the domains of these
Dirac operators, considered as unbounded closed operators on the
Hilbert space of $\banach{2}$--spinors, as the images
in the $\sobolev{1,p}$--spinors of the
$\banach{2}$--spinors under the resolvents.
The corresponding closed unbounded operators are defined
as the restrictions of the Dirac operators from the
$\sobolev{1,p}$--spinors into the $\banach{p}$--spinors.

\noindent{\bf 1. Uniformization of compact Riemann surfaces.}
We choose a holomorphic complex line bundle, which is a square root of
the canonical bundle. The corresponding holomorphic structures are
given by Dirac operators with potentials. In order to develop the
spectral theory of these Dirac operators we represent the compact
Riemann surfaces of genus larger than one as quotients $\mathbb{D}/\Gamma$
of the hyperbolic disk $\mathbb{D}$ modulo a Fuchsian group
and elliptic curves as quotients $\mathbb{C}/\lattice$ of $\mathbb{C}$
modulo a lattice $\lattice$ \cite[Chapter~IV.5.]{FK}.
The corresponding group action has a fundamental domain denoted by
$\Delta$ \cite[Chapter~IV.9.]{FK}. The corresponding spin bundle is
an induced bundle of an representation of $\Gamma$ and $\lattice$,
respectively. Therefore we shall consider the Dirac operators and
their resolvents on the Riemann sphere $\mathbb{P}$
with the elliptic metric $\frac{dz d\Bar{z}}{(1+z\Bar{z})^{2}}$,
on the complex plane $\mathbb{C}$ with flat metric $dz d\Bar{z}$
and on the hyperbolic disk $\mathbb{D}$
with hyperbolic metric $\frac{dz d\Bar{z}}{(1-z\Bar{z})^{2}}$.
The corresponding Dirac operators are of the form \cite[Chapter 3.4.]{Fr1}
\begin{align*}
\begin{pmatrix}
0 & (1+z\Bar{z})\partial\\
-(1+z\Bar{z})\Bar{\partial} & 0
\end{pmatrix} &\text{ on }\mathbb{P},& \begin{pmatrix}
0 & \partial\\
-\Bar{\partial} & 0
\end{pmatrix} &\text{ on }\mathbb{C} \text{ and}& \begin{pmatrix}
0 & (1-z\Bar{z})\partial\\
-(1-z\Bar{z})\Bar{\partial} & 0
\end{pmatrix}\text{ on }\mathbb{D}.
\end{align*}

\noindent{\bf 2. Green's functions.}
We shall calculate the integral kernels of the corresponding resolvents,
which are the inverse of the operators
\begin{equation*}\begin{array}{ccc}
\begin{pmatrix}
\sqrt{-1}\lambda & -(1+z\Bar{z})\partial\\
(1+z\Bar{z})\Bar{\partial} & \sqrt{-1}\lambda
\end{pmatrix},& \begin{pmatrix}
\sqrt{-1}\lambda & -\partial\\
\Bar{\partial} & \sqrt{-1}\lambda
\end{pmatrix}\text{ and}& \begin{pmatrix}
\sqrt{-1}\lambda & -(1-z\Bar{z})\partial\\
(1-z\Bar{z})\Bar{\partial} & \sqrt{-1}\lambda
\end{pmatrix}.
\end{array}\end{equation*}
The composition of these operators with the operators
\begin{equation*}\begin{array}{ccc}
\begin{pmatrix}
\sqrt{-1}\lambda & (1+z\Bar{z})\partial\\
-(1+z\Bar{z})\Bar{\partial}& \sqrt{-1}\lambda
\end{pmatrix},& \begin{pmatrix}
\sqrt{-1}\lambda & \partial\\
-\Bar{\partial} & \sqrt{-1}\lambda
\end{pmatrix}\text{ and}& \begin{pmatrix}
\sqrt{-1}\lambda & (1-z\Bar{z})\partial\\
-(1-z\Bar{z})\Bar{\partial} & \sqrt{-1}\lambda
\end{pmatrix}
\end{array}\end{equation*}
are equal to the operators
\begin{multline*}
\begin{pmatrix}
-\lambda^2+(1+z\Bar{z})^2\partial\Bar{\partial}+
(1+z\Bar{z})\Bar{z}\Bar{\partial} & 0\\
0 & -\lambda^2+(1+z\Bar{z})^2\Bar{\partial}\partial+
(1+z\Bar{z})\partial
\end{pmatrix},\\\begin{pmatrix}
-\lambda^2+\partial\Bar{\partial} & 0\\
0 & -\lambda^2+\Bar{\partial}\partial
\end{pmatrix}\text{ and}\\
\begin{pmatrix}
-\lambda^2+(1-z\Bar{z})^2\partial\Bar{\partial}-
(1-z\Bar{z})\Bar{z}\Bar{\partial} & 0\\
0 & -\lambda^2+(1-z\Bar{z})^2\Bar{\partial}\partial-
(1-z\Bar{z})z\partial
\end{pmatrix}.
\end{multline*}
On functions, which depend only on $r=|z|$ these operators
act as diagonal matrices, whose entries are the operators
\begin{multline*}
-\lambda^2+\left(\frac{1+r^2}{2}\right)^2
\left(\frac{d^2}{dr^2}+\frac{1}{r}\frac{d}{dr}\right)+
\frac{1+r^2}{2}r\frac{d}{dr},\\
-\lambda^2+\frac{1}{4}\left(\frac{d^2}{dr^2}+\frac{1}{r}\frac{d}{dr}\right)
\text{ and}\\
-\lambda^2+\left(\frac{1-r^2}{2}\right)^2
\left(\frac{d^2}{dr^2}+\frac{1}{r}\frac{d}{dr}\right)-
\frac{1-r^2}{2}r\frac{d}{dr}\text{, respectively.}
\end{multline*}
the substitutions $r=\tan(y/2)$, $r=y/2$  and $r=\tanh(y/2)$ transforms
these operators into
\begin{equation*}\begin{array}{ccc}
\displaystyle{-\lambda^2+
\frac{d^2}{dy^2}+\frac{1}{\sin(y)}\frac{d}{dy}},&
\displaystyle{-\lambda^2+
\frac{d^2}{dy^2}+\frac{1}{y}\frac{d}{dy}}\text{ and}&
\displaystyle{-\lambda^2+
\frac{d^2}{dy^2}+\frac{1}{\sinh(y)}\frac{d}{dy}}.
\end{array}\end{equation*}
We remark that in all three cases $y$ is twice the distance from the
origin. They define self--adjoint operators on the Hilbert spaces
corresponding to the measure spaces
\begin{align*}
\frac{\pi\sin(y)dy}{\cos(y)+1}&\text{ on }y\in[0,\pi],&
\frac{\pi ydy}{2}&\text{ on }y\in[0,\infty),&
\frac{\pi\sinh(y)dy}{\cosh(y)+1}&\text{ on }y\in[0,\infty).
\end{align*}
Let $\Func{G}_{\mathbb{P},\lambda}$,
$\Func{G}_{\mathbb{C},\lambda}$ and $\Func{G}_{\mathbb{D},\lambda}$
denote the corresponding Green's functions, i.\ e.\ the applications
of the three operators above on these functions yields the
$\delta$--function with respect to the corresponding measures
(which are equal to the usual two--dimensional $\delta$--function on
the complex plane).
Due to \cite[Chapter~V. \S3.1 and \S6.5]{St} and \cite[Section~7.2]{GJ}
the second function $\Func{G}_{\mathbb{C},\lambda}$ is for $y>0$ given by
\begin{eqnarray*}
\Func{G}_{\mathbb{C},\lambda}(y)&=&
-\int\limits_{-\infty}^{\infty}\int\limits_{-\infty}^{\infty}
\frac{\exp\left(-\pi\sqrt{-1}yk_1\right)}
     {\lambda^2+\pi^2(k_1^2+k_2^2)}dk_1dk_2\\
&=&-\int\limits_{-\infty}^{\infty}\int\limits_{-\infty}^{\infty}
\frac{\exp\left(-\sqrt{-1}|\lambda|yk_1\right)}
     {\pi^2(1+k_1^2+k_2^2)}dk_1dk_2\\
&=&-\int\limits_{-\infty}^{\infty}
\frac{\exp\left(-|\lambda|y\sqrt{1+k^2}\right)}
     {\pi\sqrt{1+k^2}}dk\\
&=&-\frac{2}{\pi}\int\limits_{1}^{\infty}
\frac{\exp\left(-|\lambda|yx\right)}
     {\sqrt{x^2-1}}dx\\
&=&-\frac{2}{\pi}\int\limits_{|\lambda|y}^{\infty}
\frac{\exp\left(-x\right)}
     {\sqrt{x^2-\lambda^2y^2}}dx.
\end{eqnarray*}
This implies that this function $\Func{G}_{\mathbb{C},\lambda}$
has the following properties
\begin{description}
\item[(i)] $\displaystyle{0<-\Func{G}_{\mathbb{C},\lambda}(y)\leq
  \text{\bf{O}}(1)\exp\left((\varepsilon-|\lambda|)y\right)}$
  with an appropriate $\varepsilon>0$ and large $|\lambda|y$.
\item[(ii)] $\displaystyle{0<\Func{G}'_{\mathbb{C},\lambda}(y)\leq
  \text{\bf{O}}(1)\exp\left((\varepsilon-|\lambda|)y\right)}$
  with an appropriate $\varepsilon>0$ and large $|\lambda|y$.
\item[(iii)] $\displaystyle{0<-\Func{G}_{\mathbb{C},\lambda}(y)\leq
  -\frac{2}{\pi}\ln(|\lambda|y)+\text{\bf{O}}(1)}$
  for small $y$.
\item[(iv)] $\displaystyle{0<\Func{G}'_{\mathbb{C},\lambda}(y)\leq
  \frac{2}{\pi}\frac{1}{y}+\text{\bf{O}}(1)}$
  for small $y$.
\end{description}
The first and the third operator may be transformed into the operators
\begin{eqnarray*}
\cos^{-1}\left(\frac{y}{2}\right)\left(\lambda^2-\frac{d^2}{dy^2}-
\frac{1}{\sin(y)}\frac{d}{dy}\right)\cos\left(\frac{y}{2}\right)&=&
\lambda^2-\frac{d^2}{dy^2}-\frac{\cos(y)}{\sin(y)}\frac{d}{dy}+
\frac{\sin^2\left(\frac{y}{2}\right)}{4\cos^2\left(\frac{y}{2}\right)}\\
\cosh^{-1}\left(\frac{y}{2}\right)\left(\lambda^2-\frac{d^2}{dy^2}-
\frac{1}{\sinh(y)}\frac{d}{dy}\right)\cosh\left(\frac{y}{2}\right)&=&
\lambda^2-\frac{d^2}{dy^2}-\frac{\cosh(y)}{\sinh(y)}\frac{d}{dy}-
\frac{\cosh^2\left(\frac{y}{2}\right)+1}{4\cosh^2\left(\frac{y}{2}\right)}
\end{eqnarray*}
Let $\Tilde{\Func{G}}_{\mathbb{P},\lambda}$,
and $\Tilde{\Func{G}}_{\mathbb{D},\lambda}$
denote the Green's functions of the operators
\begin{align*}
&-\lambda^2+\frac{d^2}{dy^2}+\frac{\cos(y)}{\sin(y)}\frac{d}{dy}&\text{and }
&-\lambda^2+\frac{d^2}{dy^2}+\frac{\cosh(y)}{\sinh(y)}\frac{d}{dy}
\end{align*}
on the measure spaces $\displaystyle{\frac{\pi\sin(y)dy}{2}}$
with $y\in[0,\pi]$ and
$\displaystyle{\frac{\pi\sinh(y)dy}{2}}$ with $y\in[0,\infty)$.
They describe the Laplace operators
acting on functions \cite[Chapter~VII \S5.]{Cha}.
The corresponding Green's functions
have representations analogous to the representation of
$\Func{G}_{\mathbb{C},\lambda}$ (compare \cite[Chapter~5]{Da}).
In fact, the substitution $x=-\cos(y)$ transforms the former operator
into $\displaystyle{-\lambda^2+(1-x^2)\frac{d^2}{dx^2}-2x\frac{d}{dx}}$,
whose eigenfunctions are the Legendre polynomials \cite[\S10.10.]{MOT}.
We apply a variant of Mehlers integral \cite[\S10.10. (43)]{MOT}:
\begin{eqnarray*}
P_n(-\cos(y))&=&P_n(\cos(\pi-y))
=\frac{1}{\pi}\int\limits_{y-\pi}^{\pi-y}
\frac{\exp\left(\sqrt{-1}(n+\frac{1}{2})x\right)dx}
     {\sqrt{2\cos(x)-2\cos(\pi-y)}}\\
&=&\frac{1}{\pi}\int\limits_{y}^{2\pi-y}
\frac{\exp\left(\sqrt{-1}(n+\frac{1}{2})(\pi-x)\right)dx}
     {\sqrt{2\cos(\pi-x)-2\cos(\pi-y)}}\\
&=&\int\limits_{y}^{\pi}
\frac{\exp\left(\sqrt{-1}(n+\frac{1}{2})(\pi-x)\right)dx}
     {\pi\sqrt{2\cos(y)-2\cos(x)}}-
\int\limits_{-\pi}^{-y}
\frac{\exp\left(\sqrt{-1}(n+\frac{1}{2})(\pi-x)\right)dx}
     {\pi\sqrt{2\cos(y)-2\cos(x)}}\\
&=&\frac{2(-1)^{n}}{\pi}\int\limits_{y}^{\pi}
\frac{\sin\left((n+\frac{1}{2})(x)\right)dx}
     {\sqrt{2\cos(y)-2\cos(x)}}.
\end{eqnarray*}
Hence we obtain (compare \cite[\S10.10. (2),(4) and (18)]{MOT})
\begin{eqnarray*}
\Tilde{\Func{G}}_{\mathbb{P},\sqrt{\lambda^2+\frac{1}{4}}}(y)
&=&-\frac{2}{\pi}\sum\limits_{n=0}^{\infty}
P_n(-\cos(y))\frac{(-1)^n(n+\frac{1}{2})}{\lambda^2+(n+\frac{1}{2})^2}\\
&=&\frac{2}{\pi^2}\int\limits_{y}^{\pi}
\sum\limits_{n=0}^{\infty}
\left(\frac{\sqrt{-1}}{\lambda-\sqrt{-1}(n+\frac{1}{2})}-
     \frac{\sqrt{-1}}{\lambda+\sqrt{-1}(n+\frac{1}{2})}\right)
\frac{\sin\left((n+\frac{1}{2})x\right)dx}
     {\sqrt{2\cos(y)-2\cos(x)}}\\
&=&\frac{1}{\pi^2}\int\limits_{y}^{\pi}
\sum\limits_{n\in\mathbb{Z}+\frac{1}{2}}
\frac{\exp\left(\sqrt{-1}nx\right)}{\lambda-\sqrt{-1}n}-
\frac{\exp\left(\sqrt{-1}nx\right)}{\lambda+\sqrt{-1}n}
\frac{dx}{\sqrt{2\cos(y)-2\cos(x)}}\\
&=&-\frac{2\pi}{\pi^2}\int\limits_{y}^{\pi}
\frac{\cosh\left(\lambda(\pi-x)\right)dx}
     {\cosh(\lambda\pi)\sqrt{2\cos(y)-2\cos(x)}}\\
&=&-\frac{2}{\pi}\int\limits_{|\lambda|y}^{|\lambda|\pi}
\frac{\cosh\left(|\lambda|\pi-x\right)dx}
     {\cosh(\lambda\pi)|\lambda|\sqrt{2\cos(y)-2\cos(\frac{x}{\lambda})}}.
\end{eqnarray*}
On the other hand, the heat kernel of the hyperbolic plane yields the
following representation
\begin{eqnarray*}
\Tilde{\Func{G}}_{\mathbb{D},\sqrt{\lambda^2-\frac{1}{4}}}(y)
&=&-\frac{1}{\pi^{3/2}}\int\limits_{0}^{\infty}\int\limits_{y}^{\infty}
\frac{x\exp\left(-\frac{t}{4}-\frac{x^2}{4t}-\lambda^2t+\frac{t}{4}\right)}
     {t^{3/2}\sqrt{2\cosh(x)-2\cosh(y)}}dx dt\\
&=&-\frac{2}{\pi}\int\limits_{y}^{\infty}
\frac{\exp\left(-|\lambda|x\right)}
     {\sqrt{2\cosh(x)-2\cosh(y)}}dx\\
&=&-\frac{2}{\pi}\int\limits_{|\lambda|y}^{\infty}
\frac{\exp\left(-x\right)}
     {|\lambda|\sqrt{2\cosh(\frac{x}{\lambda})-2\cosh(y)}}dx.
\end{eqnarray*}
We conclude that the Green's functions
$\Tilde{\Func{G}}_{\mathbb{P},\lambda}$, and
$\Tilde{\Func{G}}_{\mathbb{D},\lambda}$ have also the properties (i)--(iv).
Due to \cite[Chapter~1.3 and Chapter~1.8]{Da}
the resolvents of the Laplace operators on the three simply connected
Riemann surfaces $\mathbb{P}$,$\mathbb{C}$ and $\mathbb{D}$
are positivity preserving. Moreover, if $\Op{H}_0$ is an
elliptic second order differential operator
and $V$ a non--negative potential, the difference of the resolvents
$$\left(\lambda^2+\Op{H}_0\right)^{-1}-
\left(\lambda^2+\Op{H}_0+V\right)^{-1}=
\left(\lambda^2+\Op{H}_0\right)^{-1}
V\left(\lambda^2+\Op{H}_0+V\right)^{-1}$$
is positivity preserving. Moreover, if in addition
$-\frac{d}{dy}\left(\lambda^2+\Op{H}_0\right)^{-1}$
is positivity preserving, then also the difference
$$\frac{d}{dy}\left(\lambda^2+\Op{H}_0+V\right)^{-1}-
\frac{d}{dy}\left(\lambda^2+\Op{H}_0\right)^{-1}=
-\frac{d}{dy}\left(\lambda^2+\Op{H}_0\right)^{-1}
V\left(\lambda^2+\Op{H}_0+V\right)^{-1}$$
is positivity preserving. Hence we may estimate the positive
Green's functions
\begin{align*}
0<-\Func{G}_{\mathbb{P},\lambda}(y)&\leq
-\Tilde{\Func{G}}_{\mathbb{P},\sqrt{\lambda^2+\frac{1}{4}}}(y)&
0<\Func{G}'_{\mathbb{P},\lambda}(y)&\leq
\Tilde{\Func{G}}'_{\mathbb{P},\sqrt{\lambda^2+\frac{1}{4}}}(y)\\
0<-\Func{G}_{\mathbb{D},\lambda}(y)&\leq
-\Tilde{\Func{G}}_{\mathbb{D},\sqrt{\lambda^2+\frac{1}{4}}}(y)&
0<\Func{G}'_{\mathbb{D},\lambda}(y)&\leq
\Tilde{\Func{G}}'_{\mathbb{D},\sqrt{\lambda^2+\frac{1}{4}}}(y)
\end{align*}
This proves

\begin{Lemma}\label{free resolvents}
The Green's functions $\Func{G}_{\mathbb{P},\lambda}$,
$\Func{G}_{\mathbb{D},\lambda}$ and $\Func{G}_{\mathbb{D},\lambda}$
have the properties
\begin{description}
\item[(i)] $\displaystyle{0<-\Func{G}_{\cdot,\lambda}(y)\leq
  \text{\bf{O}}(1)\exp\left((\varepsilon-|\lambda|)y\right)}$
  with an appropriate $\varepsilon>0$ and large $|\lambda|y$.
\item[(ii)] $\displaystyle{0<\Func{G}'_{\cdot,\lambda}(y)\leq
  \text{\bf{O}}(1)\exp\left((\varepsilon-|\lambda|)y\right)}$
  with an appropriate $\varepsilon>0$ and large $|\lambda|y$.
\item[(iii)] $\displaystyle{0<-\Func{G}_{\cdot,\lambda}(y)\leq
  -\frac{2}{\pi}\ln(|\lambda|y)+\text{\bf{O}}(1)}$
  for small $y$.
\item[(iv)] $\displaystyle{0<\Func{G}'_{\cdot,\lambda}(y)\leq
  \frac{2}{\pi}\frac{1}{y}+\text{\bf{O}}(1)}$
  for small $y$.\qed
\end{description}
\end{Lemma}

\noindent{\bf 3. Integral kernels of the resolvents of Dirac operators
on simply connected Riemann surfaces.}
The free Dirac operator on $\mathbb{C}$ is translation invariant.
Moreover, the free Dirac operators on $\mathbb{P}$ and $\mathbb{D}$
are invariant under group actions of the subgroups
$SU(2)$ and $SU(1,1)$ of the M\"o\-bius group, respectively.
More precisely, for $\left(\begin{smallmatrix}
a & b\\
c & d
\end{smallmatrix}\right)\in SU(2)$ and $SU(1,1)$ the transformation
\begin{align*}
z'&=\frac{az+b}{cz+d}\text{ implies}&
\partial'&=(cz+d)^2\partial,&
\Bar{\partial}'&=(\Bar{c}\Bar{z}+\Bar{d})^2\Bar{\partial},&
1\pm z'\Bar{z}'&=\frac{1\pm z\Bar{z}}{|cz+d|^2}\text{ and}
\end{align*}
\begin{multline*}
\begin{pmatrix}
\sqrt{-1}\lambda & -(1\pm z'\Bar{z}')\partial'\\
(1\pm z'\Bar{z}')\Bar{\partial}' & \sqrt{-1}\lambda
\end{pmatrix}=\\=
\begin{pmatrix}
cz+d & 0\\
0 & \Bar{c}\Bar{z}+\Bar{d}
\end{pmatrix}\comp\begin{pmatrix}
\sqrt{-1}\lambda & -(1\pm z\Bar{z})\partial\\
(1\pm z\Bar{z})\Bar{\partial} & \sqrt{-1}\lambda
\end{pmatrix}\comp\begin{pmatrix}
cz+d & 0\\
0 & \Bar{c}\Bar{z}+\Bar{d}
\end{pmatrix}^{-1},
\end{multline*}
\begin{multline*}
\begin{pmatrix}
\sqrt{-1}\lambda & (1\pm z'\Bar{z}')\partial'\\
-(1\pm z'\Bar{z}')\Bar{\partial}' & \sqrt{-1}\lambda
\end{pmatrix}=\\=
\begin{pmatrix}
cz+d & 0\\
0 & \Bar{c}\Bar{z}+\Bar{d}
\end{pmatrix}\comp\begin{pmatrix}
\sqrt{-1}\lambda & (1\pm z\Bar{z})\partial\\
-(1\pm z\Bar{z})\Bar{\partial} & \sqrt{-1}\lambda
\end{pmatrix}\comp\begin{pmatrix}
cz+d & 0\\
0 & \Bar{c}\Bar{z}+\Bar{d}
\end{pmatrix}^{-1},
\end{multline*}
respectively. Therefore the spin bundle of the compact Riemann surface
$\mathbb{P}$ is the trivial $\mathbb{C}^2$-- bundle on the
two members of the covering
$\mathbb{P}=\{z\in\mathbb{C}\}\cup\{z'\in\mathbb{C}\}$ with the
transformation $z'=-1/z$ and the transition matrix
$\left(\begin{smallmatrix}
z & 0\\
0 & \Bar{z}
\end{smallmatrix}\right)$, which transforms the spinors on
$\{z\in\mathbb{C}\}$ into the spinors on $\{z'\in\mathbb{C}\}$.
The spin bundles of $\mathbb{C}$ and $\mathbb{D}$ are the trivial
$\mathbb{C}^2$--bundles over these non--compact Riemann surfaces.
The translation invariance of the free Dirac operator on $\mathbb{C}$
implies that the resolvent
$\Op{R}_{\mathbb{C}}(0,0,\sqrt{-1}\lambda)=\left(\begin{smallmatrix}
\sqrt{-1}\lambda & -\partial\\
\Bar{\partial} & \sqrt{-1}\lambda
\end{smallmatrix}\right)^{-1}$ has the integral kernel
$\Func{K}_{\mathbb{C},\lambda}(z,z')
\frac{d\Bar{z}'\wedge dz'}{2\sqrt{-1}}$ with
$$\Func{K}_{\mathbb{C},\lambda}(z,z')=\begin{pmatrix}
\sqrt{-1}\lambda & \partial\\
-\Bar{\partial} & \sqrt{-1}\lambda
\end{pmatrix}\begin{pmatrix}
\Func{G}_{\mathbb{C},\lambda}(2|z-z'|)& 0\\
0 & \Func{G}_{\mathbb{C},\lambda}(2|z-z'|)
\end{pmatrix}\frac{d\Bar{z}'\wedge dz'}{2\sqrt{-1}}.$$
On $\mathbb{P}$ and $\mathbb{D}$ we use the invariance under
$SU(2)$ and $SU(1,1)$. The transformed coordinate under the
           M\"obius transformation
$\left(\begin{smallmatrix}
1 & -z'\\
\pm\Bar{z}' & 1
\end{smallmatrix}\right)/\sqrt{1\pm z'\Bar{z}'}\in SU(2)$ and $SU(1,1)$
vanishes at $z'\in\mathbb{P}$ and $\mathbb{D}$, respectively.
Therefore the integral kernels of
the resolvents $\left(\begin{smallmatrix}
\sqrt{-1}\lambda & -(1\pm z\Bar{z})\partial\\
(1\pm z\Bar{z})\Bar{\partial} & \sqrt{-1}\lambda
\end{smallmatrix}\right)^{-1}$
on $\mathbb{P}$ and $\mathbb{D}$ have the integral kernels
$\Func{K}_{\mathbb{P},\lambda}(z,z')
\frac{d\Bar{z}'\wedge dz'}{2\sqrt{-1}(1+z'\Bar{z}')^2}$ and
$\Func{K}_{\mathbb{D},\lambda}(z,z')
\frac{d\Bar{z}'\wedge dz'}{2\sqrt{-1}(1-z'\Bar{z}')^2}$ with
\begin{align*}
\Func{K}_{\mathbb{P},\lambda}(z,z')&=
\begin{pmatrix}
\sqrt{-1}\lambda & (1+z\Bar{z})\partial\\
-(1+z\Bar{z})\Bar{\partial} & \sqrt{-1}\lambda
\end{pmatrix}\begin{pmatrix}
\frac{1+z'\Bar{z}'}{1+z\Bar{z}'}
\Func{G}_{\mathbb{P},\lambda}\left(2d_{\mathbb{P}}(z,z')\right) & 0\\
0 &\frac{1+z'\Bar{z}'}{1+\Bar{z}z'}
\Func{G}_{\mathbb{P},\lambda}\left(2d_{\mathbb{P}}(z,z')\right)
\end{pmatrix}\\
\Func{K}_{\mathbb{D},\lambda}(z,z')&=
\begin{pmatrix}
\sqrt{-1}\lambda & (1-z\Bar{z})\partial\\
-(1-z\Bar{z})\Bar{\partial} & \sqrt{-1}\lambda
\end{pmatrix}\begin{pmatrix}
\frac{1-z'\Bar{z}'}{1-z\Bar{z}'}
\Func{G}_{\mathbb{D},\lambda}\left(2d_{\mathbb{D}}(z,z')\right) & 0\\
0 &\frac{1-z'\Bar{z}'}{1-\Bar{z}z'}
\Func{G}_{\mathbb{D},\lambda}\left(2d_{\mathbb{D}}(z,z')\right)
\end{pmatrix}.
\end{align*}
Here $d_{\mathbb{P}}(z,z')$ and $d_{\mathbb{D}}(z,z')$ denote the
distance between $z$ and $z'$ with respect to the invariant metrics
$\frac{dz d\Bar{z}}{(1\pm z\Bar{z})^2}$ on $\mathbb{P}$ and
$\mathbb{D}$, respectively.

\noindent{\bf 4. Integral kernels of the resolvents of
Dirac operators on compact Riemann surfaces.}
A spin bundle of the elliptic Riemann surface $\mathbb{C}/\lattice$
is the trivial $\mathbb{C}^2$ bundle. Finally, a spin bundle
of the hyperbolic compact Riemann surface $\mathbb{D}/\Gamma$ is
the induced bundle of the discrete Fuchsian subgroup
$\Gamma\subset SU(1,1)$ of the following action
on the sections of the trivial spin bundle on $\mathbb{D}$:
\begin{align*}
\begin{pmatrix}
a & b\\
c & d
\end{pmatrix}&\text{ acts on spinors as }
\psi\mapsto\psi'\text{ with}&
\psi'(z)&=\begin{pmatrix}
a-cz & 0\\
0 & \Bar{a}-\Bar{c}\Bar{z}
\end{pmatrix}^{-1}\psi\left(\frac{dz-b}{a-cz}\right).
\end{align*}
Consequently, the resolvent
$\Op{R}_{\mathbb{C}/\lattice}(0,0,\sqrt{-1}\lambda)=\left(\begin{smallmatrix}
\sqrt{-1}\lambda & -\partial\\
\Bar{\partial} & \sqrt{-1}\lambda
\end{smallmatrix}\right)^{-1}$
of the free Dirac operator on $\mathbb{C}/\lattice$ has
the integral kernel
$$\sum\limits_{\gamma\in\lattice}
\Func{K}_{\mathbb{C},\lambda}(z,z'+\gamma)
\frac{d\Bar{z}'\wedge dz'}{2\sqrt{-1}}$$
with $z,z'\in\Delta$.
Due to property~(i) of Lemma~\ref{free resolvents} these sum converges
for all non--vanishing real $\lambda$.
Analogously, the resolvent 
$\Op{R}_{\mathbb{D}/\Gamma}(0,0,\sqrt{-1}\lambda)=\left(\begin{smallmatrix}
\sqrt{-1}\lambda & -(1-z\Bar{z})\partial\\
(1-z\Bar{z})\Bar{\partial} & \sqrt{-1}\lambda
\end{smallmatrix}\right)^{-1}$
of the free Dirac operator on $\mathbb{D}/\Gamma$ has the integral kernel
$$\sum\limits_{\left(\begin{smallmatrix}
               a & b\\
               c & d
               \end{smallmatrix}\right)\in\Gamma}
\Func{K}_{\mathbb{D},\lambda}\left(z,\frac{dz'-c}{a-cz'}\right)
\comp\begin{pmatrix}
a-cz' & 0\\
0 & \Bar{a}-\Bar{c}\Bar{z}'
\end{pmatrix}\frac{d\Bar{z}'\wedge dz'}{2\sqrt{-1}(1-z'\Bar{z}')^2}.$$

\noindent{\bf 5. Banach spaces of spinors}
The $\banach{p}$--spinors on $\mathbb{C}$ belong to
$\banach{p}(\mathbb{C},\mathbb{H})$.
Moreover, on $\mathbb{P}$ and $\mathbb{D}$ the $\banach{p}$--spinors
have finite norms
\begin{align*}
\|f\|&=\left(\int\limits_{\mathbb{P}}
|f(z)|^p(1+z\Bar{z})^{\left(\frac{p}{2}-2\right)}
\frac{d\Bar{z}\wedge dz}{2\sqrt{-1}}\right)^{\frac{1}{p}}&
\|f\|&=\left(\int\limits_{\mathbb{D}}
|f(z)|^p(1-z\Bar{z})^{\left(\frac{p}{2}-2\right)}
\frac{d\Bar{z}\wedge dz}{2\sqrt{-1}}\right)^{\frac{1}{p}},
\end{align*}
which are invariant under the actions of $SU(2)$ and $SU(1,1)$.
Moreover, the $\banach{p}$--spinors on $\mathbb{C}/\lattice$ are
defined as sections of the spin bundle on $\mathbb{C}/\lattice$,
with finite norm
$$\|f\|=\left(\int\limits_{\Delta}
|f|^p\frac{d\Bar{z}\wedge dz}{2\sqrt{-1}}\right)^{\frac{1}{p}}.$$
Finally, the $\banach{p}$--spinors on $\mathbb{D}/\Gamma$ are
defined as section of the spin bundle on $\mathbb{D}/\Gamma$,
with finite norm
$$\|f\|=\left(\int\limits_{\Delta}
|f|^p(1-z\Bar{z})^{\frac{p}{2}}
\frac{d\Bar{z}\wedge dz}{2\sqrt{-1}(1-z\Bar{z})^{2}}\right)^{\frac{1}{p}}.$$
Let $\Op{R}_{\mathbb{P}}(V,W,\lambda)$ and
$\Op{R}_{\mathbb{D}}(V,W,\lambda)$ denote the resolvents
\begin{align*}
\begin{pmatrix}
\lambda-V & -(1+z\Bar{z})\partial\\
(1+z\Bar{z})\Bar{\partial} & \lambda-W
\end{pmatrix}^{-1}&&
\begin{pmatrix}
\lambda-V & -(1-z\Bar{z})\partial\\
(1-z\Bar{z})\Bar{\partial} & \lambda-W
\end{pmatrix}^{-1}
\end{align*}
considered as operators from
the $\banach{p}$--spinors into
the $\banach{q}$--spinors. The corresponding free resolvents
$\Op{R}_{\mathbb{P}}(0,0,\sqrt{-1}\lambda)$ and
$\Op{R}_{\mathbb{D}}(0,0,\sqrt{-1}\lambda)$ have the integral kernels
\begin{align*}
\Func{K}_{\mathbb{P},\lambda}(z,z')&
\frac{d\Bar{z}\wedge dz}{2\sqrt{-1}(1+z'\Bar{z}')^{2}}&
\Func{K}_{\mathbb{D},\lambda}(z,z')&
\frac{d\Bar{z}\wedge dz}{2\sqrt{-1}(1-z\Bar{z})^{2}}.
\end{align*}
Analogously let $\Op{R}_{\mathbb{C}/\lattice}(V,W,\lambda)$
and $\Op{R}_{\mathbb{D}/\Gamma}(V,W,\lambda)$ denote the resolvents of
Dirac operators on with potentials $V$ and $W$
on the compact Riemann surfaces
$\mathbb{C}/\lattice$ and $\mathbb{D}/\Gamma$.
Therefore for all compact Riemann surfaces
$\Spa{X}=\mathbb{P},\mathbb{C}/\lattice,\mathbb{D}/\Gamma$
the resolvents of the Dirac operators with potentials
$V$ and $W$ are equal to
\begin{eqnarray*}
\Op{R}_{\Spa{X}}(V,W,\sqrt{-1}\lambda)&=&
\left(\unity-\Op{R}_{\Spa{X}}(0,0,\sqrt{-1}\lambda)\comp
\begin{pmatrix}
V & 0\\
0 & W
\end{pmatrix}\right)^{-1}\comp
\Op{R}_{\Spa{X}}(0,0,\sqrt{-1}\lambda)\\
&=&\Op{R}_{\Spa{X}}(0,0,\sqrt{-1}\lambda)\comp
\left(\unity-\begin{pmatrix}
V & 0\\
0 & W
\end{pmatrix}\comp
\Op{R}_{\Spa{X}}(0,0,\sqrt{-1}\lambda)\right)^{-1}.
\end{eqnarray*}

\noindent{\bf 6. The resolvents of Dirac operators with
$\banach{2}$--potentials on compact Riemann surfaces.}

\begin{Theorem}\label{weakly continuous resolvent}
For all $1\leq p,q<\infty$ with $\frac{1}{p}<\frac{1}{q}+\frac{1}{2}$,
there exists a constant $C_p>0$, with the following property:
For all compact Riemann surfaces $\Spa{X}$ and all $\varepsilon>0$ there
exists a $\delta>0$ such that for all real
$\lambda\in(-\infty,-\delta)\cup(\delta,\infty)$ the mapping
$(V,W)\mapsto\Op{R}_{\Spa{X}}(V,W,\sqrt{-1}\lambda)$ is holomorphic
and weakly continuous from the weakly compact space of all potentials,
whose restrictions to all $\varepsilon$--balls of $\Spa{X}$
(with respect to the elliptic, flat or hyperbolic metric)
have $\banach{2}$--norm not greater than $C_p$
(with respect to the induced measure), into the compact operators from
the $\banach{p}$--spinors into the $\banach{q}$--spinors.
\end{Theorem}

\begin{proof}
Due to \cite[Chapter~V. \S.3.4 Lemma~3.]{St} the operators
$\partial\comp
\left(\unity-\partial\Bar{\partial}\right)^{-\frac{1}{2}}$
and $\Bar{\partial}\comp
\left(\unity-\partial\Bar{\partial}\right)^{-\frac{1}{2}}$
are bounded operators on $\banach{p}(\mathbb{C})$ with $1<p<2$.
Therefore any function $f\in\banach{p}(\mathbb{C})$ with
either $\partial f\in\banach{p}(\mathbb{C})$ or
$\Bar{\partial}f\in\banach{p}(\mathbb{C})$ belongs to the Sobolev space
$\sobolev{1,p}(\mathbb{C})$ \cite[Chapter~V. \S.3.4 Theorem~3.]{St}.
Since $\mathbb{C}$ and $\mathbb{D}$ are homogeneous spaces, they obey
the assumptions of the Sobolev Embedding \cite[Theorem~2.21]{Au} on
Riemannian manifolds. We conclude that for all $1<p<2$ the resolvent
$\Op{R}_{\Spa{X}}(0,0,\sqrt{-1})$ considered as an operator from
the space of $\banach{p}$--spinors into the space of
$\banach{\frac{2p}{2-p}}$--spinors are bounded.
Moreover, due to Lemma~\ref{free resolvents} there exists a constant
$C_p>S_p^{-1}$ such that
$\|\Op{R}_{\Spa{X}}(0,0,\sqrt{-1}\lambda)\|\leq S_p$
for all $\lambda\in(-\infty,-1)\cup(1,\infty)$.
Now we decompose this resolvent into the sum
$$\Op{R}_{\Spa{X}}(0,0,\sqrt{-1}\lambda)=
\Op{R}_{\Spa{X},\varepsilon'\text{\scriptsize\rm--near}}
(0,0,\sqrt{-1}\lambda)+
\Op{R}_{\Spa{X},\varepsilon'\text{\scriptsize\rm--distant}}
(0,0,\sqrt{-1}\lambda),$$
whose integral kernel either vanish or are equal to the integral
kernel of $\Op{R}_{\Spa{X}}(0,0,\sqrt{-1}\lambda)$,
in cases that $z$ and $z'$ have distance larger than $\varepsilon$
or smaller than $\epsilon$ and vice versa, respectively.
Obviously the norm of the first term is smaller than $S_p$.
If the potentials $V$ and $W$ belong to the set described in the Lemma,
the operator $\left(\begin{smallmatrix}
V & 0\\
0 & W
\end{smallmatrix}\right)\comp
\Op{R}_{\Spa{X},\varepsilon'\text{\scriptsize\rm--near}}
(0,0,\sqrt{-1}\lambda)$ has for small $\varepsilon'$
norm smaller than $\left(
\frac{\vol(B(0,\varepsilon+\varepsilon'))}{\vol(B(0,\varepsilon))}
\right)^{1/p}$. In fact, for all $x\in\Spa{X}$ the restriction
of $\Op{R}_{\Spa{X},\varepsilon'\text{\scriptsize\rm--near}}
(0,0,\sqrt{-1}\lambda)\psi$ to $B(x,\varepsilon)$ is
smaller than the norm of the restriction of $\psi$ to
$B(x,\varepsilon+\varepsilon'))$. Since $\Spa{X}$ is either the
homogeneous space $\mathbb{P}$ or a quotient of the homogeneous spaces
$\mathbb{C}$ or $\mathbb{D}$ by a discrete group, for all small
$\varepsilon$ and all $\banach{p}$--functions on
$\Spa{X}$, the $\banach{p}$--norm of the function
$$x\mapsto\left\|f\mid_{B(x,\varepsilon)}\right\|_p$$
is equal to $\vol^{\frac{1}{p}}(B(0,\varepsilon))$ times the
$\banach{p}$--norm $\|f\|_p$ of $f$.
We conclude that for small $\varepsilon'$ the operator
$\left(\begin{smallmatrix}
V & 0\\
0 & W
\end{smallmatrix}\right)\comp
\Op{R}_{\Spa{X},\varepsilon'\text{\scriptsize\rm--near}}
(0,0,\sqrt{-1}\lambda)$ has norm smaller than $1$.
Due to Lemma~\ref{free resolvents}~(ii),
in the limit $|\lambda|\rightarrow\infty$ the norm of the second
term converge to zero
$\lim\limits_{|\lambda|\rightarrow\infty}
\|\Op{R}_{\Spa{X},\varepsilon'\text{\scriptsize\rm--distant}}
(0,0,\sqrt{-1}\lambda)\|=0.$
Hence there exists a $\delta>0$,
such that the operator $\left(\begin{smallmatrix}
V & 0\\
0 & W
\end{smallmatrix}\right)\comp
\Op{R}_{\Spa{X}}(0,0,\sqrt{-1}\lambda)$ on the $\banach{p}$--spinors
has for all $\lambda\in(-\infty,-\delta)\cup(\delta,\infty)$
norm smaller than $1$. Consequently the von Neumann series
$$\Op{R}_{\Spa{X}}(V,W,\sqrt{-1}\lambda)=
\sum\limits_{l=0}^{\infty}\Op{R}_{\Spa{X}}(V,W,\sqrt{-1}\lambda)\comp
\left(\left(\begin{smallmatrix}
V & 0\\
0 & W
\end{smallmatrix}\right)\comp
\Op{R}_{\Spa{X}}(0,0,\sqrt{-1}\lambda)\right)^{l}$$
converges to a holomorphic function with values in the operators form
the $\banach{p}$--spinors into the $\banach{\frac{2p}{2-p}}$--spinors.
Moreover, due to Kondrakov's Theorem \cite[Theorem~2.34]{Au} the resolvent
$\Op{R}_{\Spa{X}}(0,0,\sqrt{-1}\lambda)$ considered as an operator from
the $\banach{p}$--spinors into the $\banach{q}$--spinors with
$1\leq q<\frac{2p}{2-p}$ is compact.
Due to \cite[Theorem~II.5.11]{LT}, all Banach spaces of
$\banach{q}$--spinors have a \Em{Schauder basis}. Consequently they
have the \Em{approximation property}, and all compact operators into
one of these Banach spaces of $\banach{q}$--spinors are norm--limits
of finite rank operators (compare \cite[Section~I.1.a]{LT}).
Hence all terms in the von Neumann series are norm--limits of weakly
continuous functions from the set described in the Lemma into the
compact operators from the $\banach{p}$--spinors into the
$\banach{q}$ spinors. Since this set is weakly compact, the uniform
limit of weakly continuous functions is again
a weakly continuous function \cite[Theorem~IV.8]{RS1}.
\end{proof}

We shall explain the relation of these Dirac operators and
holomorphic structures. In the introduction we mentioned already,
that Dirac operators on $\mathbb{C}$ are the composition of
holomorphic structures with an invertible operator
$$\begin{pmatrix}
U & \partial\\
-\Bar{\partial} & \Bar{U}
\end{pmatrix}=\begin{pmatrix}
0 & \unity\\
-\unity & 0
\end{pmatrix}\comp\begin{pmatrix}
\Bar{\partial} & -\Bar{U}\\
U & \partial
\end{pmatrix}.$$
The Dirac operator on $\mathbb{P}$ and $\mathbb{D}$ with
potentials $(1\pm z\Bar{z})U$ and $(1\pm z\Bar{z})\Bar{U}$ are equal to
$$\begin{pmatrix}
(1\pm z\Bar{z})U & (1\pm z\Bar{z})\partial\\
-(1\pm z\Bar{z})\Bar{\partial} & (1\pm z\Bar{z})\Bar{U}
\end{pmatrix}=(1\pm z\Bar{z})\begin{pmatrix}
0 & \unity\\
-\unity & 0
\end{pmatrix}\comp\begin{pmatrix}
\Bar{\partial} & -\Bar{U}\\
U & \partial
\end{pmatrix}.$$
Hence these Dirac operators are also compositions of
holomorphic structures with invertible operators. We remark that the
$\banach{2}$--norms of the potentials $(1\pm z\Bar{z})U$ with respect
to the induced measures
$\frac{d\Bar{z}\wedge dz}{2\sqrt{-1}(1\pm z\Bar{z})^2}$ coincides with
the integrals over $\frac{1}{2\sqrt{-1}}Q\wedge\Bar{Q}$
with the corresponding Hopf fields $Q=-\Bar{U}d\Bar{z}$.

Finally let us deduce a simple criterion for a bounded sequence of
square integrable Hopf fields on a compact Riemann surface $\Spa{X}$,
whether they contain subsequences in the sets of the form
described in Theorem~\ref{weakly continuous resolvent} or not.
Due to the Banach--Alaoglu theorem \cite[Theorem~IV.21]{RS1} and the
Riesz Representation theorem \cite[Chapter~13 Section~5]{Ro2}
any bounded sequence $\left(Q_n\right)_{n\in\mathbb{N}}$ of square
integrable Hopf fields has a subsequence,
with the property that the corresponding
sequence of measures $\frac{1}{2\sqrt{-1}}Q_n\wedge\Bar{Q}_n$
converge weakly to a finite Baire measure on $\Spa{X}$.
This subsequence is contained in a set of the form described in
Theorem~\ref{weakly continuous resolvent}, if the limit
of the measure does not contain point measures of mass
larger or equal to the constant $S_p^{-2}$.
In fact, if the weak limit of the measures
does not contain point measures of mass larger or equal to $S_p^{-2}$,
we may cover $\Spa{X}$ by open sets, whose measures with
respect to the limit of the measures is smaller than $S_p^{-2}$.
Due to the compactness of $\Spa{X}$ this open covering has a finite
subcovering. For any finite open covering, the function on
$\Spa{X}$, which associates to each $x$ the radius of the
maximal open disk around $x$, which is entirely contained in one member
of the covering, is continuous. We exclude the trivial case,
where one member of the covering contains the whole of $\Spa{X}$ and
therefore all disks. So this function is the maximum of the
distances of the corresponding point to all complements of the
members of the covering. Therefore, there exists a small
$\varepsilon>0$, such that all disks with radius $2\varepsilon$
are contained in one member of the finite subcovering.
Obviously, for any member of the subcovering
there exists a continuous $[0,1]$--valued function,
whose support is contained in this member of the subcovering,
and which is equal to $1$ on those disks $B(x,\varepsilon)$,
whose extensions $B(x,2\varepsilon)$ are contained
in this member of the subcovering.
Since the sequence of measures converges weakly,
the integrals of these functions with respect to
the measures $\frac{1}{2\sqrt{-1}}Q_n(x)\wedge\Bar{Q}_n(x)$
corresponding to the sequence are also smaller than $S_p^{-2}$,
with the exception of finitely many elements of the sequence.
This shows that with the exception of
finitely many elements of the sequence,
the $\banach{2}$--norms of the restrictions of $Q_n$ to all
$\varepsilon$-balls are smaller than some $C_p<S_p^{-1}$.

\begin{Lemma}\label{bounded point measures}
If a weak limit of the sequence of finite Baire measures
$\frac{1}{2\sqrt{-1}}Q_n\wedge\Bar{Q}_nd^2x$
on the compact Riemann surface $\Spa{X}$ does not contain point
measures with mass larger or equal to $S_p^{-2}$,
then there exists a $C_p<S_p^{-1}$, an $\varepsilon>0$
and a subsequence of the bounded sequence $Q_n$ of square integrable
Hopf fields, whose $\banach{2}$--norms
of the restrictions of $Q_n$ to all $\varepsilon$--balls is smaller
than $C_p$.\qed
\end{Lemma}

\section{The Riemann--Roch Theorem}\label{section riemann roch}
In this section we shall prove that all
holomorphic quaternionic line bundles with square integrable
Hopf fields obey S\'{e}rre Duality and the Riemann--Roch Theorem.
In general, the holomorphic sections of a
holomorphic quaternionic line bundle with square integrable
Hopf field are not continuous.
Therefore we cannot use a non--trivial meromorphic section in order to
determine the Chern class of the bundle (compare \cite[\S2.3]{FLPP}).
Hence we use sheaf theory.

Let $U$ be a potential in
$\banach{2}_{\text{\scriptsize\rm loc}}(\Omega)$
over an open subset $\Omega\subset\mathbb{C}$.
For all $1<p<2$ the operator $\left(\begin{smallmatrix}
\Bar{\partial} & -\Bar{U}\\
U & \partial
\end{smallmatrix}\right)$ defines a linear operator from
$\sobolev{1,p}_{\text{\scriptsize\rm loc}}(\Omega,\mathbb{H})$
onto $\banach{p}_{\text{\scriptsize\rm loc}}(\Omega,\mathbb{H})$.
Due to \cite[Chapter~V. \S.3.4 Lemma~3.]{St} the operators
$\partial\comp
\left(\unity-\partial\Bar{\partial}\right)^{-\frac{1}{2}}$
and $\Bar{\partial}\comp
\left(\unity-\partial\Bar{\partial}\right)^{-\frac{1}{2}}$
are bounded operators on $\banach{p}(\mathbb{C})$ with $1<p<2$.
Therefore the operator $\left(\begin{smallmatrix}
\Bar{\partial} & 0\\
0 & \partial
\end{smallmatrix}\right)$ defines an isomorphism from
$\sobolev{1,p}(\Omega,\mathbb{H})$ onto
$\banach{p}(\Omega,\mathbb{H})$.
For any holomorphic quaternionic line bundle on a
Riemann surface $\Spa{X}$, whose holomorphic complex line bundle
corresponds to $\Sh{O}_{D}$,
let $\Sh{Q}_{D}$ denote the sheaf of holomorphic sections,
and $\Sh{W}^{1,p}_{D}$ the corresponding sheaf of
$\sobolev{1,p}_{\text{\scriptsize\rm loc}}$--sections.
Moreover, let $\Sh{L}^p_{D-K}$ denote the sheaf of
$\banach{p}_{\text{\scriptsize\rm loc}}$--sections of the
quaternionic line bundle corresponding to $\Sh{Q}_{D}$
tensored with the inverse of the canonical line bundle.
If $g dz d\Bar{z}$ denotes a hermitian metric with respect to local
coordinates $z$ on the compact Riemann surface $\Spa{X}$,
then the local operators
$$\frac{1}{g}\begin{pmatrix}
\Bar{\partial} & -\Bar{U}\\
U & \partial
\end{pmatrix}$$
fit together to a global operator from
$H^0\left(\Spa{X},\Sh{W}^{1,p}_{D}\right)$ into
$H^0\left(\Spa{X},\Sh{L}^p_{D-K}\right)$.
In fact, under the transformation $z\mapsto z'$ this operator
transforms to
$$\frac{1}{g'}\begin{pmatrix}
\Bar{\partial'} & -\Bar{U}'\\
U' & \partial'
\end{pmatrix}=
\frac{1}{g}\left|\frac{dz'}{dz}\right|^2\begin{pmatrix}
\overline{\frac{dz}{dz'}} & 0\\
0 &\frac{dz}{dz'}
\end{pmatrix}\begin{pmatrix}
\Bar{\partial} & -\Bar{U}\\
U & \partial
\end{pmatrix}=\frac{1}{g}\begin{pmatrix}
\frac{dz'}{dz} & 0\\
0 & \overline{\frac{dz'}{dz}}
\end{pmatrix}\begin{pmatrix}
\Bar{\partial} & -\Bar{U}\\
U & \partial
\end{pmatrix}.$$
Moreover, the holomorphic cocycle
of the underlying holomorphic complex line bundle
does only change the Hopf field $Q=-\Bar{U}d\Bar{z}$.
Consequently, the holomorphic structure
of the quaternionic line bundle, which is locally given by
operators of the form $\frac{1}{g}\left(\begin{smallmatrix}
\Bar{\partial} & -\Bar{U}\\
U & \partial
\end{smallmatrix}\right)$,
induces a morphism $\Sh{W}^{1,p}_{D}\rightarrow\Sh{L}^p_{D-K}$ which
fits to the following exact sequence of sheaves \cite[\S2.2]{FLPP}
$$0\rightarrow\Sh{Q}_{D}\hookrightarrow\Sh{W}_{D}^{1,p}\rightarrow
\Sh{L}^p_{D-K}\rightarrow 0.$$
Standard arguments \cite[Theorem~12.6.]{Fo}
show that the first cohomology group of the sheaf $\Sh{W}^{1,p}_{D}$ vanish.
Consequently, the corresponding long exact cohomology sequence
\cite[\S15.]{Fo} shows that the cokernel of the holomorphic structure,
considered as a Fredholm operator from
$H^0\left(\Spa{X},\Sh{W}^{1,p}_{D}\right)$ into
$H^0\left(\Spa{X},\Sh{L}^p_{D-K}\right)$ is naturally isomorphic
to the first cohomology group $H^1\left(\Spa{X},\Sh{Q}_{D}\right)$
of the sheaf $\Sh{Q}_{D}$. On the other hand, this cokernel is dual to
the kernel of the transposed operator acting on the dual space of
$H^0\left(\Spa{X},\Sh{L}^p_{D-K}\right)$, which is equal to
$H^0\left(\Spa{X},\Sh{L}^\frac{p}{p-1}_{K-D}\right)$.
This transposed operator defines a natural holomorphic structure
on the quaternionic line bundle over $\Sh{O}_{K-D}$ \cite[\S2.3.]{FLPP}.
The corresponding sheaf of holomorphic sections is denoted by
$\Sh{Q}_{K-D}$. Hence we have proven (compare \cite[\S8.--\S9.]{Na})

\newtheorem{Serre Duality}[Lemma]{S\'{e}rre Duality}
\begin{Serre Duality}\label{serre duality}
Let $\Spa{X}$ be a compact Riemann surface and $\Sh{Q}_{D}$ the sheaf of
holomorphic sections of a holomorphic structure
with a square integrable Hopf field $Q$
(i.\ e.\ $\frac{1}{2\sqrt{-1}}\int\limits_{\Spa{X}}Q\wedge\Bar{Q}<\infty$)
on the quaternionic line bundle over the complex line bundle
corresponding to $\Sh{O}_{D}$. Then the \u{C}ech cohomology groups
$H^1\left(\Spa{X},\Sh{Q}_{D}\right)$ and
$H^0\left(\Spa{X},\Sh{Q}_{K-D}\right)$ are naturally dual to each
other.\qed
\end{Serre Duality}

\newtheorem{Riemann--Roch Theorem}[Lemma]{Riemann--Roch Theorem}
\begin{Riemann--Roch Theorem}\label{riemann roch}
Let $\Spa{X}$ be a compact Riemann surface and $\Sh{Q}_{D}$ the sheaf of
holomorphic sections of a holomorphic quaternionic line bundle
with square integrable Hopf fields
(i.\ e.\ $\frac{1}{2\sqrt{-1}}\int\limits_{\Spa{X}}Q\wedge\Bar{Q}<\infty$).
over the complex line bundle corresponding to $\Sh{O}_{D}$.
Then the quaternionic dimensions of the corresponding
\u{C}ech cohomology groups are finite and obey the formula
$$\dim_{\mathbb{Q}}H^0\left(\Spa{X},\Sh{Q}_{D}\right)-
\dim_{\mathbb{Q}}H^1\left(\Spa{X},\Sh{Q}_{D}\right)
=1-g+\deg(D).$$
\end{Riemann--Roch Theorem}

\begin{proof}
Due to the long exact cohomology sequence corresponding
to the exact sequence of sheaves \cite[\S15.]{Fo}
$$0\rightarrow\Sh{Q}_{D}\rightarrow\Sh{Q}_{D'}\rightarrow
\Sh{Q}_{D'}/\Sh{Q}_{D}\rightarrow 0$$ and Lemma~\ref{quotient dimension}
the \De{Riemann--Roch Theorem} for the sheaf
$\Sh{Q}_{D}$ is equivalent to the
\De{Riemann--Roch Theorem} for the sheaf
$\Sh{Q}_{D'}$ with $D\leq D'$. Since for all pairs of divisors $D$ and
$D'$ there exists a divisor $D''$ with $D\leq D''$ and $D'\leq D''$,
this equivalence holds also for arbitrary $D$ and $D'$.
Consequently, it suffices to proof the \De{Riemann--Roch Theorem}
for the holomorphic quaternionic line bundles with one fixed
underlying holomorphic complex line bundle.
Theorem~\ref{weakly continuous resolvent} shows that
holomorphic structures with square integrable Hopf fields on the spin bundle,
considered as a quaternionic line bundle, are Fredholm operators
of index zero from $H^0\left(\Spa{X},\Sh{W}^{1,p}_{D}\right)$ into
$H^0\left(\Spa{X},\Sh{L}^p_{D}\right)$, where $D$ is the
corresponding square root of the canonical divisor, i.\ e.\ $2D=K$.
This implies $\deg(D)=g-1$ and the claim follows from the proof of
\De{S\'{e}rre Duality}~\ref{serre duality}.
\end{proof}

\section{A B\"acklund transformation}
\label{section baecklund}

In the following discussion concerning this transformation we make use
of the \Em{Lorentz spaces} \index{Lorentz spaces} $\banach{p,q}$.
These rearrangement invariant Banach spaces are an extension of the
family of the usual Banach spaces $\banach{p}$ indexed by an additional
parameter $1\leq q\leq\infty$ for $1<p<\infty$. For $p=1$ or
$p=\infty$ we consider only the \Em{Lorentz spaces} $\banach{p,\infty}$,
which in these cases are isomorphic to $\banach{p}$ 
(\cite[Chapter~V. \S3.]{SW},\cite[Chapter~4 Section~4.]{BS} 
and \cite[Chapter~1. Section~8.]{Zi}).
We recall some properties of these Banach spaces:
\begin{description}
\item[(i)] For $1<p\leq\infty$ the \Em{Lorentz spaces} $\banach{p,p}$
  coincide with the usual $\banach{p}$--spaces. Moreover, the
  \Em{Lorentz space} $\banach{1,\infty}$ coincides with the usual Banach
  space $\banach{1}$.
\item[(ii)] On a finite measure space the \Em{Lorentz space}
  $\banach{p,q}$ is contained in $\banach{p',q'}$
  either if $p>p'$ or if $p=p'$ and $q\leq q'$.
\end{description}
In \cite{O} H\"older's inequality and Young's inequality are generalized
to these \Em{Lorentz spaces}
(\cite[Chapter~4 Section~7.]{BS} and \cite[Chapter~2. Section~10.]{Zi}):

\newtheorem{Generalized Hoelder}[Lemma]{Generalized H\"older's inequality}
\index{generalized!H\"older's inequality}
\begin{Generalized Hoelder}\label{generalized hoelder}
Either for $1/p_1+1/p_2=1/p_3<1$ and $1/q_1+1/q_2\geq 1/q_3$ or for
$1/p_1+1/p_2=1$, $1/q_1+1/q_2\geq 1$ and $(p_3,q_3)=(1,\infty)$
there exists some constant $C>0$ with
$$\|fg\|_{(p_3,q_3)}\leq C\|f\|_{(p_1,q_1)}\|g\|_{(p_2,q_2)}.$$
\end{Generalized Hoelder}
\newtheorem{Generalized Young}[Lemma]{Generalized Young's inequality}
\index{generalized!Young's inequality}
\begin{Generalized Young}\label{generalized young}
Either for $1/p_1+1/p_2-1=1/p_3>0$ and $1/q_1+1/q_2\geq 1/q_3$ or for
$1/p_1+1/p_2=1$, $1/q_1+1/q_2\geq 1$ and $(p_3,q_3)=(\infty,\infty)$
there exists some constant $C>0$ with
$$\|f\ast g\|_{(p_3,q_3)}\leq C\|f\|_{(p_1,q_1)}\|g\|_{(p_2,q_2)}.$$
\end{Generalized Young}

Therefore, the resolvent of the Dirac operators
on $\mathbb{C}$ is a bounded  operators
from the $\banach{1}$--spinors into the $\banach{2,\infty}$--spinors,
from the $\banach{2,1}$--spinors into the continuous spinors,
from the $\banach{p}$--spinors into the $\banach{q,p}$--spinors,
and finally from the $\banach{p,q}$--spinors into the $\banach{q}$--spinors,
with $1<p<2$ and $q=2p/(2-p)$.
Moreover, the Sobolev constant $S_p$
(compare with Lemma~\ref{free resolvents}
and Theorem~\ref{weakly continuous resolvent})
is equal to the corresponding norm $\|f\|_{2,\infty}$
times the corresponding constant of the
\De{Generalized Young's inequality}~\ref{generalized young}.

Let $\xi$ and $\chi$ be two elements in the kernel of
$\left(\begin{smallmatrix}
\Bar{\partial} & -\Bar{A}\\
A & \partial
\end{smallmatrix}\right)$ on an open domain $\Omega\subset\mathbb{C}$.
If $\chi$ does not vanish, then the quotient of these two holomorphic
sections of the corresponding holomorphic quaternionic line bundle is
equal to
$$\begin{pmatrix}
\chi_1 & -\Bar{\chi}_2\\
\chi_2 & \Bar{\chi}_1
\end{pmatrix}^{-1}\begin{pmatrix}
\xi_1 & -\Bar{\xi}_2\\
\xi_2 & \Bar{\xi}_1
\end{pmatrix}=
\frac{1}{\chi_1\Bar{\chi}_1+\chi_2\Bar{\chi}_2}\begin{pmatrix}
\Bar{\chi}_1 & \Bar{\chi}_2\\
-\chi_2 & \chi_1
\end{pmatrix}\begin{pmatrix}
\xi_1 & -\Bar{\xi}_2\\
\xi_2 & \Bar{\xi}_1
\end{pmatrix}.$$
The derivatives of these quaternionic--valued functions are equal to
\begin{eqnarray*}
d\begin{pmatrix}
\xi_1 & -\Bar{\xi}_2\\
\xi_2 & \Bar{\xi}_1
\end{pmatrix}&=&\left(\begin{pmatrix}
d\Bar{z} & 0\\
0 & dz
\end{pmatrix}\begin{pmatrix}
0 & \Bar{A}\\
-A & 0
\end{pmatrix}+\begin{pmatrix}
dz & 0\\
0 & d\Bar{z}
\end{pmatrix}\begin{pmatrix}
\partial & 0\\
0 & \Bar{\partial}
\end{pmatrix}\right)\begin{pmatrix}
\xi_1 & -\Bar{\xi}_2\\
\xi_2 & \Bar{\xi}_1
\end{pmatrix}\\
d\begin{pmatrix}
\chi_1 & -\Bar{\chi}_2\\
\chi_2 & \Bar{\chi}_1
\end{pmatrix}&=&\left(\begin{pmatrix}
d\Bar{z} & 0\\
0 & dz
\end{pmatrix}\begin{pmatrix}
0 & \Bar{A}\\
-A & 0
\end{pmatrix}-\begin{pmatrix}
dz & 0\\
0 & d\Bar{z}
\end{pmatrix}\begin{pmatrix}
B & \Bar{U}\\
-U & \Bar{B}
\end{pmatrix}\right)
\begin{pmatrix}
\chi_1 & -\Bar{\chi}_2\\
\chi_2 & \Bar{\chi}_1
\end{pmatrix}\text{, with}\\
\begin{pmatrix}
B & \Bar{U}\\
-U & \Bar{B}
\end{pmatrix}&=&-\begin{pmatrix}
\partial\chi_1 & -\partial\Bar{\chi}_2\\
\Bar{\partial}\chi_2 & \Bar{\partial}\Bar{\chi}_1
\end{pmatrix}\begin{pmatrix}
\chi_1 & -\Bar{\chi}_2\\
\chi_2 & \Bar{\chi}_1
\end{pmatrix}^{-1}.
\end{eqnarray*}
Therefore, the derivative of the foregoing quotient is equal to
\begin{multline*}
d\begin{pmatrix}
\chi_1 & -\Bar{\chi}_2\\
\chi_2 & \Bar{\chi}_1
\end{pmatrix}^{-1}\begin{pmatrix}
\xi_1 & -\Bar{\xi}_2\\
\xi_2 & \Bar{\xi}_1
\end{pmatrix}=\\
\begin{aligned}
&=\begin{pmatrix}
\chi_1 & -\Bar{\chi}_2\\
\chi_2 & \Bar{\chi}_1
\end{pmatrix}^{-1}\left(d\begin{pmatrix}
\xi_1 & -\Bar{\xi}_2\\
\xi_2 & \Bar{\xi}_1
\end{pmatrix}-\left(d\begin{pmatrix}
\chi_1 & -\Bar{\chi}_2\\
\chi_2 & \Bar{\chi}_1
\end{pmatrix}\right)\begin{pmatrix}
\chi_1 & -\Bar{\chi}_2\\
\chi_2 & \Bar{\chi}_1
\end{pmatrix}^{-1}\begin{pmatrix}
\xi_1 & -\Bar{\xi}_2\\
\xi_2 & \Bar{\xi}_1
\end{pmatrix}\right)\\
&=\begin{pmatrix}
\chi_1 & -\Bar{\chi}_2\\
\chi_2 & \Bar{\chi}_1
\end{pmatrix}^{-1}\begin{pmatrix}
dz & 0\\
0 & d\Bar{z}
\end{pmatrix}\left(\begin{pmatrix}
\partial & 0\\
0 & \Bar{\partial}
\end{pmatrix}+\begin{pmatrix}
B & \Bar{U}\\
-U & \Bar{B}
\end{pmatrix}\right)\begin{pmatrix}
\xi_1 & -\Bar{\xi}_2\\
\xi_2 & \Bar{\xi}_1
\end{pmatrix}.
\end{aligned}\end{multline*}
The non--vanishing section $\left(\begin{smallmatrix}
\chi_1 & -\Bar{\chi}_2\\
\chi_2 & \Bar{\chi}_1
\end{smallmatrix}\right)$ induces a flat connection on the
quaternionic line bundle. The zero curvature equation takes the from
$$\left[\begin{pmatrix}
\partial+B & \Bar{U}\\
A & \partial
\end{pmatrix},\begin{pmatrix}
\Bar{\partial} & -\Bar{A}\\
-U & \Bar{\partial}+\Bar{B}
\end{pmatrix}\right]=0.$$
In the framework of `quaternionic function theory' \cite{FLPP}
this equation takes the form 
$$\begin{pmatrix}
\Bar{\partial} & -\Bar{U}\\
U & \partial
\end{pmatrix}\begin{pmatrix}
\partial+B & \Bar{U}\\
-U & \Bar{\partial}+\Bar{B}\end{pmatrix}=\begin{pmatrix}
\partial+B & \Bar{A}\\
-A & \Bar{\partial}+\Bar{B}
\end{pmatrix}\begin{pmatrix}
\Bar{\partial} & -\Bar{A}\\
A & \partial
\end{pmatrix}.$$
This implies the equation
$$\begin{pmatrix}
\Bar{\partial} & -\Bar{U}\\
U & \partial
\end{pmatrix}\begin{pmatrix}
\partial+B & \Bar{U}\\
-U & \Bar{\partial}+\Bar{B}\end{pmatrix}\begin{pmatrix}
\xi_1 & -\Bar{\xi}_2\\
\xi_2 & \Bar{\xi}_1
\end{pmatrix}=0.$$
Therefore the quaternionic--valued function
$\left(\begin{smallmatrix}
\psi_1 & \Bar{\psi}_2\\
\psi_2 & \Bar{\psi}_1
\end{smallmatrix}\right)=\left(\begin{smallmatrix}
\partial+B & \Bar{U}\\
-U & \Bar{\partial}+\Bar{B}
\end{smallmatrix}\right)\left(\begin{smallmatrix}
\xi_1 & -\Bar{\xi}_2\\
\xi_2 & \Bar{\xi}_1
\end{smallmatrix}\right)$ belongs to the kernel of
$\left(\begin{smallmatrix}
\Bar{\partial} & -\Bar{U}\\
U & \partial
\end{smallmatrix}\right)$.
On the other hand the quaternionic--valued function
$\left(\begin{smallmatrix}
\chi_1 & \chi_2\\
-\Bar{\chi}_2 & \Bar{\chi}_1
\end{smallmatrix}\right)^{-1}$ obeys the differential equation
\begin{multline*}
\begin{pmatrix}
\Bar{\partial} & 0\\
0 & \partial
\end{pmatrix}\begin{pmatrix}
\chi_1 & \chi_2\\
-\Bar{\chi}_2 & \Bar{\chi}_1
\end{pmatrix}^{-1}=
\begin{pmatrix}
\Bar{\partial} & 0\\
0 & \partial
\end{pmatrix}\begin{pmatrix}
\Bar{\chi}_1 & -\chi_2\\
\Bar{\chi}_2 & \chi_1
\end{pmatrix}\frac{1}{\chi_1\Bar{\chi}_1+\chi_2\Bar{\chi}_2}\\
\begin{aligned}
&&&=\left(\begin{pmatrix}
-\Bar{B} & -U\\
\Bar{U} & -B
\end{pmatrix}-\begin{pmatrix}
\Bar{\partial}\ln(\chi_1\Bar{\chi}_1+\chi_2\Bar{\chi}_2) & 0\\
0 & \partial\ln(\chi_1\Bar{\chi}_1+\chi_2\Bar{\chi}_2)
\end{pmatrix}\right)\begin{pmatrix}
\chi_1 & \chi_2\\
-\Bar{\chi}_2 & \Bar{\chi}_1
\end{pmatrix}^{-1}\\
&&&=\begin{pmatrix}
0 & -U\\
\Bar{U} & 0
\end{pmatrix}\begin{pmatrix}
\chi_1 & \chi_2\\
-\Bar{\chi}_2 & \Bar{\chi}_1
\end{pmatrix}^{-1}.
\end{aligned}\end{multline*}
Here we used
\begin{align*}
\Bar{\partial}\ln(\chi_1\Bar{\chi}_1+\chi_2\Bar{\chi}_2)&=
\frac{\chi_1\Bar{\partial}\Bar{\chi}_1+\Bar{\chi}_2\Bar{\partial}\chi_2}
     {\chi_1\Bar{\chi}_1+\chi_2\Bar{\chi}_2}=-\Bar{B}&
\partial\ln(\chi_1\Bar{\chi}_1+\chi_2\Bar{\chi}_2)&=
\frac{\Bar{\chi}_1\partial\chi_1+\chi_2\partial\Bar{\chi}_2}
     {\chi_1\Bar{\chi}_1+\chi_2\Bar{\chi}_2}=-B.
\end{align*}
Therefore this function belongs to the kernel of
$\left(\begin{smallmatrix}
\Bar{\partial} & U\\
-\Bar{U} & \partial
\end{smallmatrix}\right)$.

\newtheorem{Baecklund}[Lemma]{B\"acklund transformation}
\begin{Baecklund}\label{baecklund}
Let $\left(\begin{smallmatrix}
\chi_1\\
\chi_2
\end{smallmatrix}\right)$ and $\left(\begin{smallmatrix}
\xi_1\\
\xi_2
\end{smallmatrix}\right)$ be two spinors in the kernel of
$\left(\begin{smallmatrix}
\Bar{\partial} & -\Bar{A}\\
A & \partial
\end{smallmatrix}\right)$
on an open domain
$\Omega\subset\mathbb{C}$ with square integrable potential
$A\in\banach{2}_{\text{\scriptsize\rm loc}}(\Omega)$.
Moreover, let $\chi$ have no zeroes on $\Omega$
(in the sense of \De{Order of zeroes}~\ref{order of zeroes}).
Then the  components of the quaternionic--valued functions
$\left(\begin{smallmatrix}
B & \Bar{U}\\
-U & \Bar{B}
\end{smallmatrix}\right)=-\left(\begin{smallmatrix}
\partial\chi_1 & -\partial\Bar{\chi}_2\\
\Bar{\partial}\chi_2 & \Bar{\partial}\Bar{\chi}_1
\end{smallmatrix}\right)\left(\begin{smallmatrix}
\chi_1 & -\Bar{\chi}_2\\
\chi_2 & \Bar{\chi}_1
\end{smallmatrix}\right)^{-1}$ belong to
$U\in\banach{2}_{\text{\scriptsize\rm loc}}(\Omega)$ and
$B\in\banach{2,\infty}_{\text{\scriptsize\rm loc}}(\Omega)$.
More precisely, the function $\Bar{\partial}B$ is a measure on
$\Omega$ without point measures. Moreover, the derivative of the
quotient $\left(\begin{smallmatrix}
\chi_1 & -\Bar{\chi}_2\\
\chi_2 & \Bar{\chi}_1
\end{smallmatrix}\right)^{-1}\left(\begin{smallmatrix}
\xi_1 & -\Bar{\xi}_2\\
\xi_2 & \Bar{\xi}_1
\end{smallmatrix}\right)$ is equal to
$$d\begin{pmatrix}
\chi_1 & -\Bar{\chi}_2\\
\chi_2 & \Bar{\chi}_1
\end{pmatrix}^{-1}\begin{pmatrix}
\xi_1 & -\Bar{\xi}_2\\
\xi_2 & \Bar{\xi}_1
\end{pmatrix}=\begin{pmatrix}
\chi_1 & -\Bar{\chi}_2\\
\chi_2 & \Bar{\chi}_1
\end{pmatrix}^{-1}\begin{pmatrix}
dz & 0\\
0 & d\Bar{z}
\end{pmatrix}\begin{pmatrix}
\partial+B & \Bar{U}\\
-U & \Bar{\partial}+\Bar{B}
\end{pmatrix}\begin{pmatrix}
\xi_1 & -\Bar{\xi}_2\\
\xi_2 & \Bar{\xi}_1
\end{pmatrix}.$$
Furthermore, the function $\left(\begin{smallmatrix}
\psi_1 & -\Bar{\psi}_2\\
\psi_2 & \Bar{\psi}_1
\end{smallmatrix}\right)=\left(\begin{smallmatrix}
\partial+B & \Bar{U}\\
-U & \Bar{\partial}+\Bar{B}
\end{smallmatrix}\right)\left(\begin{smallmatrix}
\xi_1 & -\Bar{\xi}_2\\
\xi_2 & \Bar{\xi}_1
\end{smallmatrix}\right)$ belongs to the kernel of
$\left(\begin{smallmatrix}
\Bar{\partial} & -\Bar{U}\\
U & \partial
\end{smallmatrix}\right)$ and the function
$\left(\begin{smallmatrix}
\phi_1 & -\Bar{\phi}_2\\
\phi_2 & \Bar{\phi}_1
\end{smallmatrix}\right)=\left(\begin{smallmatrix}
\chi_1 & \chi_2\\
-\Bar{\chi}_2 & \Bar{\chi}_1
\end{smallmatrix}\right)^{-1}$ belongs to the kernel of
$\left(\begin{smallmatrix}
\Bar{\partial} & U\\
-\Bar{U} & \partial
\end{smallmatrix}\right)$. In particular, the quotient
$\left(\begin{smallmatrix}
\chi_1 & -\Bar{\chi}_2\\
\chi_2 & \Bar{\chi}_1
\end{smallmatrix}\right)^{-1}\left(\begin{smallmatrix}
\xi_1 & -\Bar{\xi}_2\\
\xi_2 & \Bar{\xi}_1
\end{smallmatrix}\right)$ belongs to
$\bigcap\limits_{1<p<2}
\sobolev{2,p}_{\text{\scriptsize\rm loc}}(\Omega,\mathbb{H})$.

Conversely, if $\left(\begin{smallmatrix}
\psi_1\\
\psi_2
\end{smallmatrix}\right)$ belong on $\Omega$ to the kernel of
$\left(\begin{smallmatrix}
\Bar{\partial} & -\Bar{U}\\
U & \partial
\end{smallmatrix}\right)$ and
$\left(\begin{smallmatrix}
\phi_1\\
\phi_2
\end{smallmatrix}\right)$ to the kernel of
$\left(\begin{smallmatrix}
\Bar{\partial} & U\\
-\Bar{U} & \partial
\end{smallmatrix}\right)$, then
$$d\begin{pmatrix}
f_1 & -\Bar{f}_2\\
f_2 & \Bar{f}_1
\end{pmatrix}=\begin{pmatrix}
\phi_1 & \phi_2\\
-\Bar{\phi}_2 & \Bar{\phi}_1
\end{pmatrix}\begin{pmatrix}
dz & 0\\
0 & d\Bar{z}
\end{pmatrix}\begin{pmatrix}
\psi_1 & -\Bar{\psi}_2\\
\psi_2 & \Bar{\psi}_1
\end{pmatrix}$$
is a closed quaternionic--valued form on $\Omega$.
If in addition $\phi$ has no zeroes on $\Omega$, then the two spinors 
$\left(\begin{smallmatrix}
\chi_1 & -\Bar{\chi}_2\\
\chi_2 & \Bar{\chi}_1
\end{smallmatrix}\right)=
\left(\begin{smallmatrix}
\phi_1 & \phi_2\\
-\Bar{\phi}_2 & \Bar{\phi}_1
\end{smallmatrix}\right)^{-1}$ and $\left(\begin{smallmatrix}
\xi_1 & -\Bar{\xi}_2\\
\xi_2 & \Bar{\xi}_1
\end{smallmatrix}\right)=\left(\begin{smallmatrix}
\phi_1 & \phi_2\\
-\Bar{\phi}_2 & \Bar{\phi}_1
\end{smallmatrix}\right)^{-1}\left(\begin{smallmatrix}
f_1 & -\Bar{f}_2\\
f_2 & \Bar{f}_1
\end{smallmatrix}\right)$ belong on $\Omega$ to the kernel of
$\left(\begin{smallmatrix}
\Bar{\partial} & -\Bar{A}\\
A & \partial
\end{smallmatrix}\right)$ with potential
$A=\frac{\Bar{\phi}_1\partial\phi_2-\phi_2\partial\Bar{\phi}_1}
       {\phi_1\Bar{\phi}_1+\phi_2\Bar{\phi}_2}
\in\banach{2}_{\text{\scriptsize\rm loc}}(\Omega)$.
\end{Baecklund}

\begin{proof} This proposition establishes a
one--to--one correspondence between two holomorphic sections $\xi$ and
$\chi$ of a holomorphic quaternionic line bundle,
and two holomorphic sections $\psi$ and $\phi$
of two paired quaternionic holomorphic line bundles.
We prove this proposition in four steps.
In steps~{1}--{3} we proof the statements concerning the
mapping from two spinors $\chi$ and $\xi$ in the kernel of
$\left(\begin{smallmatrix}
\Bar{\partial} & -\Bar{A}\\
A & \partial\\
\end{smallmatrix}\right)$ to two `paired' spinors $\psi$ and $\phi$.
In the final step we prove the statements concerning the
inverse transformation from two `paired' spinors
$\phi$ and $\psi$ to two `holomorphic' spinors $\xi$ and $\chi$
of one holomorphic quaternionic bundle.

\noindent{\bf 1. For potentials
$A\in\banach{2,1}_{\text{\scriptsize\rm loc}}(\Omega)$.}
If the Hopf field belongs locally to $\banach{2,1}(\Omega)$,
then the \De{Generalized H\"older's inequality}~\ref{generalized hoelder}
and \De{Generalized Young's inequality}~\ref{generalized young}
together with the arguments in section~\ref{section local} imply that the
holomorphic sections $\chi$ and $\xi$ are continuous and belong to the
Sobolev space
$\sobolev{1,2}_{\text{\scriptsize\rm loc}}(\Omega,\mathbb{H})$.
Hence for non--vanishing $\chi$ the potentials $U$ and $B$ belong to
$\banach{2}_{\text{\scriptsize\rm loc}}(\Omega)$ and
$\left(\begin{smallmatrix}
\chi_1 & -\Bar{\chi}_2\\
\chi_2 & \Bar{\chi}_1
\end{smallmatrix}\right)^{-1}$ is continuous and belong to
$\sobolev{1,2}_{\text{\scriptsize\rm loc}}(\Omega,\mathbb{H})$.
In this case the statements concerning $\left(\begin{smallmatrix}
\psi_1 & -\Bar{\psi}_2\\
\psi_2 & \Bar{\psi}_1
\end{smallmatrix}\right)=\left(\begin{smallmatrix}
\partial+B & \Bar{U}\\
-U & \Bar{\partial}+\Bar{B}
\end{smallmatrix}\right)\left(\begin{smallmatrix}
\xi_1 & -\Bar{\xi}_2\\
\xi_2 & \Bar{\xi}_1
\end{smallmatrix}\right)$ and $\left(\begin{smallmatrix}
\phi_1 & -\Bar{\phi}_2\\
\phi_2 & \Bar{\phi}_1
\end{smallmatrix}\right)=\left(\begin{smallmatrix}
\chi_1 & \chi_2\\
-\Bar{\chi}_2 & \Bar{\chi}_1
\end{smallmatrix}\right)^{-1}$ follow from the foregoing calculations.

\noindent{\bf 2. For $\chi=\left(\unity+
\Op{I}_{\Omega}(0)\comp\left(\begin{smallmatrix}
0 & -\Bar{A}\\
A & 0
\end{smallmatrix}\right)\right)^{-1}\left(\begin{smallmatrix}
a\\
b
\end{smallmatrix}\right)$ with $(a,b)\in\mathbb{P}$.}
We shall extend the arguments of step~{1} with a limiting argument
to small potentials $A\in\banach{2}(\Omega)$.
In fact, for any small $A\in\banach{2}(\Omega)$ we choose a
sequence $A_n$ of smooth potentials in $\banach{2}(\Omega)$ with limit $A$.
We extend all potentials to a slightly larger open domain $\Omega'$
which contains the closure $\Bar{\Omega}$, so that they vanish on the
relative complement of $\Omega$ in $\Omega'$.
Obviously, the corresponding sequence of spinors $\chi_n$ defined
above extend to $\Omega'$.
By definition these spinors are smooth on
$\Omega'\setminus\Bar{\Omega}$.
Furthermore, the sequence of integrals of the corresponding
one--forms $B_ndz$ along a closed path in
$\Omega'\setminus\Bar{\Omega}$ around $\Omega$ converges.
Since the sequence of measures
$\frac{1}{2\sqrt{-1}}A_n\Bar{A}_nd\Bar{z}\wedge dz$ converges,
this implies that the sequence $U_n$ is a bounded sequence in
$\banach{2}(\Omega)$. 
Due to the Banach--Alaoglu theorem \cite[Theorem~IV.21]{RS1},
this sequence $U_n$ has a weakly convergent subsequence
with limit $U$. Also the sequence of real signed measures
$\frac{1}{2\sqrt{-1}}(A_n\Bar{A}_n-U_n\Bar{U}_n)d\Bar{z}\wedge dz$
on $\Omega$ has a weakly convergent subsequence.
Finally, due to the equations $\Bar{\partial}B=A\Bar{A}-U\Bar{U}$,
the sequence of functions $B_n$ is bounded
in the \Em{Lorentz space} $\banach{2,\infty}(\Omega)$.
Due to \cite[Chapter~2 Theorem~2.7. and Chapter~4 Corollary~4.8.]{BS}
this \Em{Lorentz space} is the dual space
of the corresponding \Em{Lorentz space} $\banach{2,1}(\Omega)$.
The sequence $B_n$ has also a weakly convergent subsequence with limit $B$
and $\Bar{\partial}B$ considered as a measure is equal to the limit of
the measures $\frac{1}{2\sqrt{-1}}
\left(A_n\Bar{A}_n-U_n\Bar{U}_n\right)d\Bar{z}\wedge dz$.
Since the sequence of spinors $\chi_n$
converges in $\banach{q}(\Omega,\mathbb{H})$,
and since the sequences $U_n$ and $B_n$ both converge weakly,
the limit $\chi$ is anti--holomorphic with respect to
the anti--holomorphic structure defined by the limits $U$ and $B$.

Next we prove that the function
$\Bar{\partial}B=-\Bar{\partial}\partial\ln
\left(\chi_1\Bar{\chi}_1+\chi_2\Bar{\chi}_2\right)$
considered as a measure contains no point measures.

\begin{Lemma}\label{zygmund estimate}
If $\Omega$ denotes a bounded open subset of $\mathbb{C}$,
then for all finite signed Baire measures $d\mu$ on $\Omega$
\cite[Chapter~13 Section~5]{Ro2}
there exists a function $h$ in the Zygmund space
$\banach{}_{\exp}(\Omega)$
\cite[Chapter~4 Section~6.]{BS} such that
$-\Bar{\partial}\partial h=d\mu$ (in the sense of distributions).
Moreover, if for a suitable $\varepsilon>0$
all $\varepsilon$--balls of $\Omega$
have measure smaller than $\pi/q$
with respect to the positive part $d\mu^+$ of the
Hahn decomposition of the finite signed Baire measure
$d\mu$ on $\Omega$ \cite[Chapter~11 Section~5]{Ro2},
then the exponentials $\exp(h)$
of all $h\in \banach{}_{\exp}(\Omega)$ with
$-\Bar{\partial}\partial h=d\mu$ belong to
$\banach{q}_{\text{\scriptsize\rm loc}}(\Omega)$.
Conversely, if the positive part $d\mu^+$ contains a point measure
with mass $\pi/q$, then the corresponding functions
$h\in \banach{}_{\exp}(\Omega)$ with
$-\Bar{\partial}\partial h=d\mu$ do not belong to
$\banach{q}_{\text{\scriptsize\rm loc}}(\Omega)$.
\end{Lemma}

\begin{proof}
Due to Dolbeault's Lemma \cite[Chapter~I Section~D 2.~Lemma]{GuRo}
the convolution with the function $-\frac{2}{\pi}\ln|z|$
defines a right inverse of the operator $-\Bar{\partial}\partial$.
Now we claim that the restriction of this convolution operator
defines a bounded operator from $\banach{1}(\Omega)$ into the
Zygmund space $\banach{}_{\exp}(\Omega)$.
Since the domain $\Omega$ is bounded, the claim
is equivalent to the analogous statement about the restriction
to $\Omega$  of the convolution with the non--negative function
$$f(z)=\begin{cases}-\frac{2}{\pi}\ln|z| & \text{if }|z|<1\\
                    0 & \text{if }1\leq |z|
       \end{cases}.$$
Associated to this function $f$ is its distribution function $\mu_f$
and its non--increasing rearrangement $f^{\ast}$
(\cite[Chapter~II \S3. Chapter~V \S3.]{SW},
\cite[Chapter~2 Section~1.]{BS} and \cite[Chapter~1. Section~8.]{Zi}):
\begin{align*}
\mu_f(s)&=\pi\exp\left(-\pi s\right) &
f^{\ast}(t)&=\begin{cases}-\frac{\ln\left(t/\pi\right)}{\pi} &
                        \text{if }0\leq t\leq\pi\\
                       0 & \text{if } \pi\leq t
          \end{cases}.
\end{align*}
If $g\in\banach{1}(\Omega)$,
then $g^{\ast\ast}(t)=\frac{1}{t}\int_{0}^{t}g^{\ast}(s)ds$
is bounded by $\|g\|_1/t$, since the $\banach{1,\infty}$--norm
$\|g\|_{(1,\infty)}=\sup\{tg^{\ast\ast}(t)\mid t>0\}=
\int_0^{\infty}g^{\ast}(t)dt$
coincides with the $\banach{1}$--norm \cite[Chapter~V (3.9)]{SW}.
Therefore, \cite[(1.8.14) and (1.8.15)]{Zi} in the proof of
\cite[1.8.8.~Lemma]{Zi} (borrowed from \cite[Lemma~1.5.]{O})
implies that the non--increasing rearrangement $h^{\ast}(t)$
of the convolution $h=f\ast g$ is bounded by
\begin{eqnarray*}
h^{\ast}(t)\leq h^{\ast\ast}(t)\leq h_2^{\ast\ast}(t)+h_1^{\ast\ast}(t)&\leq &
g^{\ast\ast}(t)\int\limits_{f^{\ast}(t)}^{\infty}\mu_{f}(s)ds-
\int\limits_{t}^{\infty}sg^{\ast\ast}(s)df^{\ast}(s)\\
&\leq &\frac{\|g\|_1}{t}\exp\left(-\pi f^{\ast}(t)\right)-
\|g\|_1\int\limits_{t}^{\infty}df^{\ast}(s)\\
&\leq &\frac{\|g\|_1}{\pi}+\|g\|_1f^{\ast}(t).
\end{eqnarray*}
Since by definition the non--increasing rearrangement $h^{\ast}(t)$
vanishes for all arguments,
which are larger than the Lebesgue measure of $\Omega$,
we conclude the validity of the following estimate:
$$\int\limits_{0}^{\infty}
\left(\exp\left(q h^{\ast}(t)\right)-1\right)dt\leq
\int\limits_{0}^{|\Omega|}
\exp\left(q h^{\ast}(t)\right)dt\leq
|\Omega|\exp\left(q\|g\|_1/\pi\right)\int\limits_{0}^{|\Omega|}
\left(t/\pi\right)^{-\|g\|_1q/\pi}dt,$$
with an obvious modification when $\pi<|\Omega|$.
Due to a standard argument \cite[Chapter~2 Exercise~3.]{BS}
the finiteness of this integral is equivalent to the statement that
$\exp|h|$ belongs to $\banach{q}(\Omega)$.
To sum up, the exponential $\exp(h)$ of the convolution $h=f\ast g$
belongs to $\banach{q}(\Omega)$, if $q<\frac{\pi}{\|g\|_1}$.
This proves the claim.
In particular, for all $g\in \banach{1}(\Omega)$ there exists
an element $h\in \banach{}_{\exp}(\Omega)$ with
$-\Bar{\partial}\partial h=g$.

Due to \cite[Chapter~4 Theorem~6.5]{BS}
$\banach{}_{\exp}(\Omega)$ is the dual space of
the Zygmund space $\banach{}\log\banach{}(\Omega)$.
Hence we shall improve the
previous estimate and show that the convolution with
$-\frac{2}{\pi}\ln|z|$ defines a bounded operator from
$\banach{}\log\banach{}(\Omega)\subset \banach{1}(\Omega)$ into
$C(\Omega)\subset \banach{}_{\exp}(\Omega)$.
By definition of the norm \cite[Chapter~4 Definition~6.3.]{BS}
$$\|g\|_{\banach{}\log\banach{}}=-\frac{1}{|\Omega|}
\int\limits_{0}^{|\Omega|}g^{\ast}(t)\ln(t/|\Omega|)dt=
\int\limits_{0}^{|\Omega|}g^{\ast\ast}(t)dt$$
we may improve the previous estimate to
\cite[(1.8.14) and (1.8.15)]{Zi}
$$h^{\ast\ast}(t)\leq
g^{\ast\ast}(t)\int\limits_{f^{\ast}(t)}^{\infty}\mu_{f}(s)ds-
\int\limits_{t}^{\infty}sg^{\ast\ast}(s)df^{\ast}(s)\leq
\frac{1}{\pi}\int\limits_{0}^{t}g^{\ast}(s)ds+
\int\limits_{t}^{\pi}g^{\ast\ast}(s)ds\leq \|g\|_{\banach{}\log\banach{}}.$$
This implies that in this case $h^{\ast\ast}(t)$ is bounded,
and consequently $h\in \banach{\infty}(\Omega)$.
Furthermore, since the function $\ln|z|$ is continuous for $z\neq 0$,
the convolution with $-\frac{2}{\pi}\ln|z|$ is a bounded operator from
$\banach{}\log\banach{}(\Omega)$ into the Banach space $C(\Omega)$.
Finally, the dual of this operator yields a bounded operator from
the Banach space of finite signed Baire measures on $\Omega$
\cite[Chapter~13 Section~5 25.~Riesz Representation Theorem]{Ro2}
into $\banach{}_{\exp}(\Omega)$.
More precisely, if the measure of $\Omega$ with respect to a
finite positive measure $d\mu$ is smaller than $\pi/q$, then the
exponential $\exp(h)$ of the corresponding function $h=f\ast d\mu$
belongs to $\banach{q}(\Omega)$.

Due to Weyl's Lemma \cite[Theorem~IX.25]{RS2} the difference of
two arbitrary functions $h_1$ and $h_2$
with $-\Bar{\partial}\partial h_1=-\Bar{\partial}\partial h_2=d\mu$
is analytic. Therefore, it suffices to show the second and third
statement of the lemma for the convolution of $-\frac{2}{\pi}\ln|z|$
with $d\mu$. Due to the boundedness of $\Omega$
we may neglect that part of this convolution,
where the former function is negative. Therefore, we may neglect the
negative part of $d\mu$ in order to bound the exponential $\exp(h)$.
The decomposition of the convolution into
an $\varepsilon$--near and an $\varepsilon$--distant part
analogous to the decomposition in the proof of
Lemma~\ref{weakly continuous resolvent} completes the proof.
\end{proof}

If the function $\Bar{\partial}B$ considered as a finite Baire measure
contains a point measure at $z=z'$
of negative mass smaller or equal to $-n\pi$,
then, due to Lemma~\ref{zygmund estimate},
the spinor $\Tilde{\chi}=\left(\begin{smallmatrix}
z-z' & 0\\
0 & \Bar{z}-\Bar{z}'
\end{smallmatrix}\right)^{-n}\chi$ belongs to
$\bigcap\limits_{q<\infty}
\banach{q}_{\text{\scriptsize\rm loc}}(\Omega,\mathbb{H})$.
This implies that $\chi$ has a zero of order $n$ at $z'$.
Hence, due to our assumptions, the masses of all point measures are
larger than $-\pi$. Again the following Lemma~\ref{zygmund estimate}
implies that $\left(\begin{smallmatrix}
\chi_1 & \chi_2\\
-\Bar{\chi}_2 & \Bar{\chi}_1
\end{smallmatrix}\right)^{-1}$ is a
$\banach{2}_{\text{\scriptsize\rm loc}}$--spinor in the kernel of
$\left(\begin{smallmatrix}
\Bar{\partial} & U\\
-\Bar{U} & \partial
\end{smallmatrix}\right)$. Since these kernels are contained in
$\bigcap\limits_{q<\infty}
\banach{q}_{\text{\scriptsize\rm loc}}(\Omega,\mathbb{H})$,
Lemma~\ref{zygmund estimate} implies
that this measure contains no point measures.

Finally we show that $\left(\begin{smallmatrix}
\psi_1 & -\Bar{\psi}_2\\
\psi_2 & \Bar{\psi}_1
\end{smallmatrix}\right)=\left(\begin{smallmatrix}
\partial+B & \Bar{U}\\
-U & \Bar{\partial}+\Bar{B}
\end{smallmatrix}\right)\left(\begin{smallmatrix}
\xi_1 & -\Bar{\xi}_2\\
\xi_2 & \Bar{\xi}_1
\end{smallmatrix}\right)$ belongs to the kernel of
$\left(\begin{smallmatrix}
\Bar{\partial} & -\Bar{U}\\
U & \partial
\end{smallmatrix}\right)$.
For this purpose we use again the sequence of smooth potentials $A_n$
in $\banach{2}(\Omega)$ with limits $A$ and the corresponding
sequence of spinors $\chi_n$ on $\Omega'$. We choose $\Omega$ small
enough such that the corresponding sequences of potentials
$U_n$ belong to the subsets described in
Theorem~\ref{weakly continuous resolvent}, on which the resolvents are
weakly continuous.
The arguments of Theorem~\ref{weakly continuous resolvent}
imply also that the sequence $\Op{I}_{\Omega}(U_n)$
considered as operators from
$\banach{p}(\Omega,\mathbb{H})$ into
$\banach{q}(\Omega,\mathbb{H})$ with
$\frac{1}{p}<\frac{1}{q}+\frac{1}{2}$
converges to $\Op{I}_{\Omega}(U)$.
Now for any quaternionic function $f$ in $\banach{q}(\Omega,\mathbb{H})$,
the sequence of quaternionic functions $\left(\begin{smallmatrix}
\partial+B_n & \Bar{U}_n\\
-U_n & \Bar{\partial}+\Bar{B}_n
\end{smallmatrix}\right)\comp\Op{I}_{\Omega}(A_n)f$
belong on the complement of the support of $f$ in $\Omega$ to the
kernel of $\left(\begin{smallmatrix}
\Bar{\partial} & -\Bar{U}_n\\
U_n & \partial
\end{smallmatrix}\right)$. Therefore it satisfies on this complement
the corresponding quaternionic version of
\De{Cauchy's Integral Formula}~\ref{cauchy formula}.
Due to continuity this implies that the limits obeys
the quaternionic version of
\De{Cauchy's Integral Formula}~\ref{cauchy formula}
in the sense of distributions
on the complement of the support of $f$ in $\Omega$.
Hence the limit belongs to the kernel of
$\left(\begin{smallmatrix}
\Bar{\partial} & -\Bar{U}\\
U & \partial
\end{smallmatrix}\right)$.
Since the spinor $\xi$ in the kernel of
$\left(\begin{smallmatrix}
\Bar{\partial} & -\Bar{A}\\
A & \partial
\end{smallmatrix}\right)$
obey the corresponding quaternionic version of
\De{Cauchy's Integral Formula}~\ref{cauchy formula},
we may represent it on any open subset, whose closure is contained in
$\Omega$, as $\xi=\Op{I}_{\Omega}(A)f$
with an appropriate $f$, whose support is disjoint
from the open subset in $\Omega$.
Furthermore, $\xi$ is the limit of $\Op{I}_{\Omega}(A_n)f$.
This implies that $\psi$ belongs on $\Omega$ to the kernel of
$\left(\begin{smallmatrix}
\Bar{\partial} & -\Bar{U}\\
U & \partial
\end{smallmatrix}\right)$.

\noindent{\bf 3. For general $\chi$.}
Due to step~{2}, the quotient $\left(\begin{smallmatrix}
\Tilde{\chi}_1 & -\Bar{\Tilde{\chi}}_2\\
\Tilde{\chi}_2 & \Bar{\Tilde{\chi}}_1
\end{smallmatrix}\right)^{-1}\left(\begin{smallmatrix}
\chi_1 & -\Bar{\chi}_2\\
\chi_2 & \Bar{\chi}_1
\end{smallmatrix}\right)$ of $\chi$ divided by the inverse of the
$\Tilde{\chi}$, which was considered in step~{2} is continuous and
belongs to $\sobolev{2,p}_{\text{\scriptsize\rm loc}}(\Omega)$.
This implies that all components of the difference
$\left(\begin{smallmatrix}
\partial\chi_1 & -\partial\Bar{\chi}_2\\
\Bar{\partial}\chi_2 & \Bar{\partial}\Bar{\chi}_1
\end{smallmatrix}\right)\left(\begin{smallmatrix}
\chi_1 & -\Bar{\chi}_2\\
\chi_2 & \Bar{\chi}_1
\end{smallmatrix}\right)^{-1}-\left(\begin{smallmatrix}
\partial\Tilde{\chi}_1 & -\partial\Bar{\Tilde{\chi}}_2\\
\Bar{\partial}\Tilde{\chi}_2 & \Bar{\partial}\Bar{\Tilde{\chi}}_1
\end{smallmatrix}\right)\left(\begin{smallmatrix}
\Tilde{\chi}_1 & -\Bar{\Tilde{\chi}}_2\\
\Tilde{\chi}_2 & \Bar{\Tilde{\chi}}_1
\end{smallmatrix}\right)^{-1}$ belong to
$\bigcap\limits_{1<p<2}
\sobolev{1,p}_{\text{\scriptsize\rm loc}}(\Omega)\times
\sobolev{1,p}_{\text{\scriptsize\rm loc}}(\Omega)$.
Now the arguments of step~{2} carry over to all $\chi$ in the kernel
of $\left(\begin{smallmatrix}
\Bar{\partial} & -\Bar{A}\\
A & \partial
\end{smallmatrix}\right)$ without zeroes on $\Omega$.

\noindent{\bf 4. Inverse transformation.}
The arguments of steps~{1}--{3} carry over and show,
that $\left(\begin{smallmatrix}
\phi_1 & \phi_2\\
-\Bar{\phi}_2 & \Bar{\phi}_1
\end{smallmatrix}\right)^{-1}$
belong on $\Omega$ to the kernel of
$\left(\begin{smallmatrix}
\Bar{\partial} & -\Bar{A}\\
A & \partial
\end{smallmatrix}\right)$.
All other statements follow from direct calculations.
\end{proof}

Actually we proved the following quaternionic version of

\newtheorem{Weyls Lemma}[Lemma]{Weyl's Lemma}
\begin{Weyls Lemma}\label{weyls lemma}
Let $\left(\begin{smallmatrix}
\phi_1 & -\Bar{\phi}_2\\
\phi_2 & \Bar{\phi}_1
\end{smallmatrix}\right)$
be spinor without zeros in the kernel of
$\left(\begin{smallmatrix}
\Bar{\partial} & U\\
-\Bar{U} & \partial
\end{smallmatrix}\right)$
with potential $U\in\banach{2}_{\text{\scriptsize\rm loc}}(\Omega)$
on a domain $\Omega\subset\mathbb{C}$.
Then a function $\left(\begin{smallmatrix}
\psi_1 & -\Bar{\psi}_2\\
\psi_2 & \Bar{\psi}_1
\end{smallmatrix}\right)
\in\banach{p}_{\text{\scriptsize\rm loc}}(\Omega,\mathbb{H})$ with $1<p<2$
belongs to the kernel of $\left(\begin{smallmatrix}
\Bar{\partial} & -\Bar{U}\\
U & \partial
\end{smallmatrix}\right)$ if
$\left(\begin{smallmatrix}
\phi_1 & \phi_2\\
-\Bar{\phi}_2 & \Bar{\phi}_2
\end{smallmatrix}\right)\left(\begin{smallmatrix}
dz & 0\\
0 & d\Bar{z}
\end{smallmatrix}\right)\left(\begin{smallmatrix}
\psi_1 & -\Bar{\psi}_2\\
\psi_2 & \Bar{\psi}_1
\end{smallmatrix}\right)$ is a closed current on $\Omega$.
\end{Weyls Lemma}

\begin{proof}
Due to the assumptions there exists a function
$f\in\bigcap\limits_{r<p}\sobolev{1,r}(\Omega,\mathbb{H})$ with
$$d\begin{pmatrix}
f_1 & -\Bar{f}_2\\
f_2 & \Bar{f}_1
\end{pmatrix}=\begin{pmatrix}
\phi_1 & \phi_2\\
-\Bar{\phi}_2 & \Bar{\phi}_1
\end{pmatrix}\begin{pmatrix}
dz & 0\\
0 & d\Bar{z}
\end{pmatrix}\begin{pmatrix}
\psi_1 & -\Bar{\psi}_2\\
\psi_2 & \Bar{\psi}_1
\end{pmatrix}.$$
Now the \De{B\"acklund transformation}~\ref{baecklund}
implies that $\chi=\left(\begin{smallmatrix}
\phi_1 & \phi_2\\
-\Bar{\phi}_2 & \Bar{\phi}_1
\end{smallmatrix}\right)^{-1}$ and $\xi=\left(\begin{smallmatrix}
\phi_1 & \phi_2\\
-\Bar{\phi}_2 & \Bar{\phi}_1
\end{smallmatrix}\right)^{-1}f$ belong to the kernel of
$\left(\begin{smallmatrix}
\Bar{\partial} & -\Bar{A}\\
A & \partial
\end{smallmatrix}\right)$ with an appropriate
$A\in\banach{2}_{\text{\scriptsize\rm loc}}(\Omega)$.
Finally, again due to the
\De{B\"acklund transformation}~\ref{baecklund}, $\psi$ belongs to
the kernel of $\left(\begin{smallmatrix}
\Bar{\partial} & -\Bar{U}\\
U & \partial
\end{smallmatrix}\right)$.
\end{proof}

\section{The Pl\"ucker formula}\label{section pluecker}

Let $H\subset H^0\left(\Spa{X},\Sh{Q}_{D}\right)$
be a quaternionic linear system in the space of holomorphic sections
of a holomorphic quaternionic line bundle
on a compact Riemann surface $\Spa{X}$.
At any point $x\in\Spa{X}$ we have a sequence
$\ord\limits_{x}\!_{1} H<\ldots<\ord\limits_{x}\!_{\dim H} H$
of \De{Orders of zeroes}~\ref{order of zeroes} of elements of $H$,
which differ only at finitely many points form the sequence
$\ord\limits_{x}\!_{1} H=1,\ldots,\ord\limits_{x}\!_{\dim H} H=\dim H$.
The order of $H$ is defined as \cite[Definition~4.2.]{FLPP}:
$$\ord H=\sum\limits_{x\in\Spa{X}}
\left(\ord\limits_{x}\!_{1} H-1\right)+\ldots+
\left(\ord\limits_{x}\!_{\dim H} H-\dim H\right).$$
For smooth Hopf fields the following estimate is proven in
\cite[Corollary~4.8.]{FLPP}:
$$\frac{1}{2\pi\sqrt{-1}}\int\limits_{\Spa{X}}
Q\wedge\Bar{Q}\geq\dim H
\left((1-g)\left(\dim H-1\right)-\deg(D)\right)+\ord H.$$
We shall show that this inequality holds
for all square integrable Hopf fields.
For this purpose we fit together the local
\De{B\"acklund transformation}~\ref{baecklund}
to a global transformation.

\begin{Corollary}\label{global baecklund}
Let $\xi,\chi\in H^0\left(\Spa{X},\Sh{Q}_{D}\right)$ be two
holomorphic spinors of a quaternionic holomorphic line bundle with
Hopf field $Q$ over the complex holomorphic line bundle
corresponding to $\Sh{O}_D$ on a compact Riemann surface $\Spa{X}$.
If $\chi$ has no zeroes, then the local
\De{B\"acklund transformation}~\ref{baecklund} induces a
global B\"acklund Transformation $Q\mapsto\Tilde{Q}$
from the Hopf field $Q$ to an Hopf field $\Tilde{Q}$ of a
quaternionic holomorphic line bundle over the
complex holomorphic line bundle corresponding to $\Sh{O}_{D+K}$ and a
paired quaternionic holomorphic line bundle over the
complex holomorphic line bundle corresponding to $\Sh{O}_{-D}$ with two
holomorphic sections, respectively. The Willmore functionals of these
Hopf fields obey the equation

\noindent\hspace{\fill}
$\displaystyle{\frac{1}{2\pi\sqrt{-1}}\int\limits_{\Spa{X}}
\left(Q\wedge\Bar{Q}-\Tilde{Q}\wedge\Bar{\Tilde{Q}}\right)=-\deg(D)}$.\qed
\end{Corollary}

An $\dim(H)$--fold application of this Corollary immediately implies
the Pl\"ucker formula. Indeed, first we choose a member $\chi$
of the linear system $H$ of lowest vanishing order at all points of
$\Spa{X}$. Since the Riemann surface has complex dimension one,
such sections of the quaternionic vector space $H$ always exists.
We change the divisor of the quaternionic holomorphic line bundle,
such that $\chi$ is a section without zeroes of
$H^0\left(\Spa{X},\Sh{Q}_D\right)$.
An application of Corollary~\ref{global baecklund} with this $\chi$
transforms the linear system
$H\subset H^0\left(\Spa{X},\Sh{Q}_D\right)$
into  a linear system
$\Tilde{H}\subset H^0\left(\Spa{X},\Sh{Q}_{D+K}\right)$
of (quaternionic) dimension $\dim H-1$.
We may repeat such an application of Corollary~\ref{global baecklund}
until we end with a trivial linear system with Hopf field $A$.
We remark that the sum over the degrees of the corresponding sequence
of quaternionic holomorphic line bundles is equal to
$\displaystyle{\deg(D)\dim H-\ord H+\sum\limits_{j=0}^{\dim H-1}j\deg(K)}$.
Consequently, these Hopf fields obey the formula
\begin{eqnarray*}
\frac{1}{2\pi\sqrt{-1}}\int\limits_{\Spa{X}}
\left(Q\wedge\Bar{Q}-A\wedge\Bar{A}\right)&=&
-\deg(D)\dim H+\ord H-\sum_{j=0}^{\dim H-1}j\deg(K)\\
&=&\dim H
\left((1-g)\left(\dim H-1\right)-\deg(D)\right)+\ord H.
\end{eqnarray*}
This implies the general

\newtheorem{Pluecker formula}[Lemma]{Pl\"ucker formula}
\begin{Pluecker formula}\label{pluecker formula}
Let $\Spa{X}$ be a compact Riemann surface and $\Sh{Q}_{D}$ the sheaf of
holomorphic sections of a holomorphic structure
with a square integrable Hopf field $Q$
(i.\ e.\ $\frac{1}{2\sqrt{-1}}\int\limits_{\Spa{X}}Q\wedge\Bar{Q}<\infty$)
on the quaternionic line bundle over the complex line bundle
corresponding to $\Sh{O}_{D}$. Then all linear systems
$H\subset H^0\left(\Spa{X},\Sh{Q}_{K-D}\right)$ obey

\noindent\hspace{\fill}
$\displaystyle{\frac{1}{2\pi\sqrt{-1}}\int\limits_{\Spa{X}}
Q\wedge\Bar{Q}\geq\dim H
\left((1-g)\left(\dim H-1\right)-\deg(D)\right)+\ord H}$.\qed
\end{Pluecker formula}

\section{Weak limits of Hopf fields}\label{section weak limits}

In this section we consider sequences of non--trivial sections of
sequences of holomorphic quaternionic line bundles over a compact
Riemann surface $\Spa{X}$. If the degrees of the underlying complex
line bundles and the Hopf fields are bounded, then these sequences
have convergent subsequences.

\begin{Theorem}\label{convergent subsequences}
Let $\psi_n$ be a sequence of non--trivial holomorphic sections of a
sequence of quaternionic line bundles over the holomorphic complex
line bundles corresponding to $\Sh{O}_{D_n}$ with Hopf fields $Q_n$.
If the sequence of degrees $\deg(D_n)$ is bounded, then the sequence
of underlying holomorphic complex line bundles has a convergent subsequence.
If in addition the sequence of Hopf fields is bounded
(i.\ e.\ $\frac{1}{2\sqrt{-1}}\int\limits_{\Spa{X}}
Q_n\wedge\Bar{Q}_n\leq C<\infty$),
then the appropriate renormalized sequence $\psi_n$ has a subsequence,
which converges to a non--trivial holomorphic section of
a holomorphic quaternionic line bundle over a
holomorphic complex line bundle corresponding to $\Sh{O}_{D}$,
where $D-D_n$ converges to an effective divisor $D'$.
Moreover, the Hopf fields is a weak limit of the Hopf fields
of the holomorphic structures corresponding to $\Sh{Q}_{D_n+D'}$.
\end{Theorem}

\begin{proof}
The proof precedes in five steps.

\noindent{\bf 1. The decomposition of the sequence of Hopf fields.}
Due to the Banach--Alaoglu theorem \cite[Theorem~IV.21]{RS1} and the
Riesz Representation theorem \cite[Chapter~13 Section~5]{Ro2}
the sequence of bounded finite Baire measures
$\frac{1}{2\sqrt{-1}}Q_n\wedge\Bar{Q}_n$ on $\Spa{X}$
has a convergent subsequence. The limit can have only
finitely many points $\{x_1,\ldots,x_L\}$,
whose mass is larger than the constant
of Theorem~\ref{weakly continuous resolvent}.
We shall decompose the sequence of Hopf fields $Q_n$ into a sum
$$Q_n=Q_{\text{\scriptsize\rm reg},n}+
\sum\limits_{l=1}^L Q_{\text{\scriptsize\rm sing},n,l}$$
of Hopf fields with disjoint support.
Here $Q_{\text{\scriptsize\rm sing},n,1},\ldots,
Q_{\text{\scriptsize\rm sing},n,L}$
are the restrictions of $Q_n$ to small disjoint balls
$B(x_1,\varepsilon_{n,l}),\ldots,B(x_L,\varepsilon_{n,L})$,
whose radii $\varepsilon_{n,l}$ tend to zero.
Consequently, $Q_{\text{\scriptsize\rm reg},n}$
are the restrictions of $Q_n$ to the complements
of the union of these balls. More precisely, we assume
\begin{description}
\item[Decomposition (i)] For all $l=1,\ldots,L$ the weak limit
  of the sequence of finite Baire measures
  $\frac{1}{2\sqrt{-1}}Q_{\text{\scriptsize\rm sing},n,l}\wedge
  \Bar{Q}_{\text{\scriptsize\rm sing},n,l}$
  \cite[Chapter~13]{Ro2}
  is equal to the point measures of the weak limit of
  $\frac{1}{2\sqrt{-1}}Q_n\wedge\Bar{Q}_n$ at $x_l$.
  Consequently, the weak limit of the measures
  $\frac{1}{2\sqrt{-1}}Q_{\text{\scriptsize\rm reg},n}\wedge
  \Bar{Q}_{\text{\scriptsize\rm reg},n}$
  is equal to the weak limit of the measures
  $\frac{1}{2\sqrt{-1}}Q_n\wedge\Bar{Q}_n$
  minus the corresponding point measures at $x_1\ldots,x_L$.
\end{description}
Obviously there are many possible choices of the sequences
$\varepsilon_{n,l}$ with this property
(e.\ g.\ for a unique choice of $\varepsilon_{n,l}$ the square of
the  $\banach{2}$--norm of $Q_{\text{\scriptsize\rm sing},n,l}$ is equal
to the mass of the point measure at $x_l$ of the weak limit of
$\frac{1}{2\sqrt{-1}}Q_n\wedge\Bar{Q}_n$). Locally we may consider the
Hopf fields $Q_{\text{\scriptsize\rm sing},n,l}$ as Hopf fields on
$\mathbb{P}$. We want to transform each of these $L$ sequences
of Hopf fields by M\"obius transformations on $\mathbb{P}$,
such that the transformed Hopf fields belong to a set of the form
described in Theorem~\ref{weakly continuous resolvent}.
The pullbacks under the inverse of the action of the M\"obius group
$SL(2,\mathbb{C})/\mathbb{Z}_2$ on $\mathbb{P}$
$$SL(2,\mathbb{C})\ni\begin{pmatrix}
a & b\\
c & d
\end{pmatrix}:\mathbb{P}\rightarrow\mathbb{P},
z\mapsto\frac{az+b}{cz+d}$$
yields a unitary representation of the M\"obius group
on the Hilbert space of square integrable Hopf fields.
In doing so we transform each of these $L$ sequences of Hopf fields
$Q_{\text{\scriptsize\rm sing},n,l}$
(considered as Hopf fields on $\mathbb{P})$
by a sequence of M\"obius transformations $g_{n,l}$ in such a way
that the transformed sequence of Hopf fields has the following property:
\begin{description}
\item[Decomposition (ii)] There exists some $\varepsilon>0$, such that
  the $\banach{2}$--norm of the restrictions of the transformed Hopf fields
  $\left(g_{n,l}^{-1}\right)^{\ast}Q_{\text{\scriptsize\rm sing},n,l}$
  to all $\varepsilon$--balls
  (with respect to the usual metric of $\mathbb{P}$) is bounded by
  the constant $C_p<S_p^{-1}$.
\end{description}
Such decompositions do not always exist. But we shall see now that,
if all points of $\mathbb{P}$ have measure smaller than $2S_p^{-2}$
with respect to the weak limit of the finite Baire measures
$\frac{1}{2\sqrt{-1}}Q_n\wedge\Bar{Q}_n$,
then these decompositions indeed exist.
The free Dirac operator on $\mathbb{P}$ is invariant
under the compact subgroup
$SU(2,\mathbb{C})/\mathbb{Z}_2$ of the M\"obius group
($\simeq SL(2,\mathbb{C})/\mathbb{Z}_2$) as well as
the usual metric on $\mathbb{P}$.
Therefore, due to the global Iwasawa decomposition
\cite[Chapter~VI Theorem~5.1]{He},
it suffices to consider in place of the whole M\"obius group only the
semidirect product of the scaling transformations
($z\mapsto\exp(t)z$ with $t\in\mathbb{R}$)
with the translations ($z\mapsto z+z_0$ with $z_0\in\mathbb{C}$).
In the sequel we assume that all $g_{n,l}$
belong to this subgroup of the M\"obius group.

\begin{Lemma}\label{optimal Moebius}
If the square of the $\banach{2}$--norm of a Hopf field $Q$ on
$\mathbb{P}$ is smaller than $2S_p^{-2}$,
then there exists a constant $C_p<S_p^{-1}$ and a M\"obius transformation
such that the $\banach{2}$--norm of the restrictions
of the transformed Hopf fields to all balls of radius $\pi/6$
is not larger than $C_p$.
\end{Lemma}

\begin{proof}
Let $r_{\max}(Q)$ be the maximum of the set
$$\left\{r\mid
\text{the $\banach{2}$--norms of the restrictions of $Q$
      to all balls of radius $r$ are not larger than $C_p$}\right\}.$$
Since the $\banach{2}$--norm of the restriction of $Q$ to a ball
depends continuously on the center and the diameter of the ball,
this set has indeed a maximum.
Moreover, there exist balls with radius $r_{\max}(Q)$,
on which the restricted Hopf field has $\banach{2}$--norm equal to $C_p$.

We claim that there exists a M\"obius transformation $h$,
such that $r_{\max}\left(\left(h^{-1}\right)^{\ast}Q\right)$
is the supremum of the set of all
$r_{\max}\left(\left(g^{-1}\right)^{\ast}Q\right)$,
where $g$ runs through the M\"obius group.
Let $g_n$ be a maximizing sequence of this set,
i.\ e.\ the limit of the sequence
$r_{\max}\left(\left(g_n^{-1}\right)^{\ast}Q\right)$
is equal to the supremum of the former set.
Since $r_{\max}\left(\left(g^{-1}\right)^{\ast}Q\right)$ is
equal to $r_{\max}(Q)$, if $g$ belongs to the subgroup
$SU(2,\mathbb{C})/\mathbb{Z}_2$ of isometries of the M\"obius group,
and due to the global Iwasawa decomposition
\cite[Chapter~VI Theorem~5.1]{He},
the sequence $g_n$ may be chosen in the semidirect product of the
scaling transformations $z\mapsto\exp(t)z$
with the translations $z\mapsto z+z_0$.
It is quite easy to see, that if the values of $t$ and $z_0$
corresponding to a sequence $g_n$ of such M\"obius transformations
are not bounded, then there exist arbitrary small balls,
on which the $\banach{2}$--norms of the restrictions of
$\left(g_n^{-1}\right)^{\ast}Q$
have subsequences converging to the $\banach{2}$--norm of $Q$.
Hence if the $\banach{2}$--norm of $Q$ is larger than $C_p$,
then the maximizing sequence of M\"obius transformations
can be chosen to be bounded and therefore has a convergent subsequence.
In this case the continuity implies the claim.
If the $\banach{2}$--norm of $Q$ is not larger than $C_p$,
then $r_{\max}\left(\left(g^{-1}\right)^{\ast}Q\right)$
does not depend on $g$ and the claim is obvious.

If the $\banach{2}$--norm of $Q$ is smaller than $\sqrt{2}C_p$,
then all $r_{\max}(Q)$--balls,
on which the restriction of $Q$ has $\banach{2}$--norm equal to $C_p$,
have pairwise non--empty intersection.
In particular, all of them have non--empty intersection
with one of these balls.
Consequently, if $r_{\max}(Q)$ is smaller than $\pi/6$,
then these $r_{\max}(Q)$--balls are contained in one hemisphere.
In this case there exists a M\"obius transformation $g$,
which enlarges $r_{\max}(Q)$
(i.\ e.\ $r_{\max}(Q)<r_{\max}\left(\left(g^{-1}\right)^{\ast}Q\right)$).
We conclude that there exist a M\"obius transformation $g$,
such that $r_{\max}\left(\left(g^{-1}\right)^{\ast}Q\right)$
is smaller than $\pi/6$.
\end{proof}

The upper bound $2S_p^{-2}$ is sharp,
because for a sequence of $\banach{2}$--Hopf fields on $\mathbb{P}$,
whose square of the absolute values considered as a sequence
of finite Baire measures converges weakly to the sum
of two point measures of mass $S_p^{-2}$ at opposite points,
the corresponding sequence of maxima of
$r_{\max}\left(\left(h^{-1}\right)^{\ast}\cdot\right)$
converges to zero.
But the lower bound $\pi/6$ is of course not optimal.

If the $\banach{2}$--norms of the Hopf fields
$Q_{\text{\scriptsize\rm sing},n,l}$ are smaller than $\sqrt{2}C_p$,
then this lemma ensures the existence of M\"obius transformations
$g_{n,l}$ with the property \Em{Decomposition~(ii)}.
In general we showed the existence of a sequence of
M\"obius transformations $h_{n,l}$, which maximizes
$r_{\max}\left(\left(g^{-1}\right)^{\ast}
Q_{\text{\scriptsize\rm sing},n,l}\right)$.
If for one $l=1,\ldots,L$ the corresponding sequences
$\left(h_{n,l}^{-1}\right)^{\ast}Q_{\text{\scriptsize\rm sing},n,l}$
do not obey condition \Em{Decomposition~(ii)},
then we apply this procedure of decomposition to
the corresponding sequence of Hopf fields
$\left(h_{n,l}^{-1}\right)^{\ast}Q_{\text{\scriptsize\rm sing},n,l}$
on $\mathbb{P}$. Consequently, we decompose the sequence
$Q_{\text{\scriptsize\rm sing},n,l}$ into a finite sum
of Hopf fields with disjoint support, such that the corresponding
Hopf fields
$\left(g_{n,l}^{-1}\right)^{\ast}Q_{\text{\scriptsize\rm sing},n,l}$
on $\mathbb{P}$ obey
the analogous conditions \Em{Decomposition}~(i).
Due to Lemma~\ref{optimal Moebius} after finitely many iterations
of this procedure of decomposing the Hopf fields into finite
sums of Hopf fields with disjoint support, we arrive at a decomposition
$$Q_n=Q_{\text{\scriptsize\rm reg},n}+
\sum\limits_{l=1}^{L'} Q_{\text{\scriptsize\rm sing},n,l}$$
of Hopf fields with disjoint support.
More precisely, the Hopf fields
$Q_{\text{\scriptsize\rm sing},n,1},\ldots,
Q_{\text{\scriptsize\rm sing},n,L'}$ are restrictions of $Q_n$
either to small balls or to the relative complements
of finitely many small balls inside of small balls.
In particular, the domains of these Hopf fields are
excluded either from the domain of $Q_{\text{\scriptsize\rm reg},n}$,
or from the domain of another $Q_{\text{\scriptsize\rm sing},n,l}$.
The former Hopf fields obey condition \Em{Decomposition}~(i)
and the latter obey condition
\begin{description}
\item[Decomposition (i')] If the domain of
  $Q_{\text{\scriptsize\rm sing},n,l'}$ is excluded from the domain
  of $Q_{\text{\scriptsize\rm sing},n,l}$,
  then the weak limit of the sequence of finite Baire measures
  $\frac{1}{2\sqrt{-1}}\left(g_{n,l}^{-1}\right)^{\ast}
  Q_{\text{\scriptsize\rm sing},n,l'}\wedge
  \Bar{Q}_{\text{\scriptsize\rm sing},n,l'}$
  on $\mathbb{C}\subset\mathbb{P}$ converges weakly
  to the point measure of the weak limit of the sequence
  $\frac{1}{2\sqrt{-1}}\left(g_{n,l}^{-1}\right)^{\ast}
  Q_{\text{\scriptsize\rm sing},n,l}\wedge
  \Bar{Q}_{\text{\scriptsize\rm sing},n,l}$
  at some point of $\mathbb{C}$, whose measure with respect to the
  latter limit is not smaller than $S_p^{-2}$.
\end{description}
All these Hopf fields $Q_{\text{\scriptsize\rm sing},n,1},\ldots,
Q_{\text{\scriptsize\rm sing},n,L'}$ obey condition
\Em{Decomposition}~(ii).
We remark that if the weak limit of the finite Baire measures
$\frac{1}{2\sqrt{-1}}\left(h_{n,l}^{-1}\right)^{\ast}
Q_{\text{\scriptsize\rm sing},n,l}\wedge
\Bar{Q}_{\text{\scriptsize\rm sing},n,l}$
on $\mathbb{C}\subset\mathbb{P}$,
where $h_{n,l}$ denotes the sequence of M\"obius transformations
maximizing $r_{\max}\left(\left(g^{-1}\right)^{\ast}
Q_{\text{\scriptsize\rm sing},n,l}\right)$,
contains at $z=\infty$ a point measure,
whose mass is not smaller than $S_p^{-2}$,
then we decompose from the sequence
$\left(h_{n,l}^{-1}\right)^{\ast}Q_{\text{\scriptsize\rm sing},n,l}$
Hopf fields, whose domains are the complement of a large ball in the
domains of these Hopf fields.
In these cases the domains of the analog to the regular sequence
of the decomposition are excluded from the domains
of the analog to the singular sequence,
whose $\banach{2}$--norm accumulates at $z=\infty$.
Since the M\"obius transformations corresponding to the former are
faster divergent then the M\"obius transformations of the latter,
the latter should be considered as less singular than the former.
Therefore, also in this case the domains of the more singular sequences
are excluded from the domains of the less singular sequences.
To sum up, the sequence $Q_{\text{\scriptsize\rm reg},n}$
of Hopf fields on $\Spa{X}$ and the sequences
$\left(g_{n,1}^{-1}\right)^{\ast}Q_{\text{\scriptsize\rm sing},n,1},\ldots,
\left(g_{n,L'}^{-1}\right)^{\ast}Q_{\text{\scriptsize\rm sing},n,L'}$
of Hopf fields on $\mathbb{P}$ belong to a set of the form
described in Theorem~\ref{weakly continuous resolvent}.

\noindent{\bf 2. Limits of the sequence of underlying holomorphic
complex line bundles.}
Due to the Banach--Alaoglu theorem \cite[Theorem~IV.21]{RS1} and the
Riesz Representation theorem \cite[Chapter~13 Section~5]{Ro2}
the sequence of finite Baire measures
$\frac{1}{2\sqrt{-1}}Q_n\wedge\Bar{Q}_n$ on
$\Spa{X}$ has a convergent subsequence.
By passing to a subsequence we achieve that the sequence of
finite Baire measures $\frac{1}{2\sqrt{-1}}Q_n\wedge\Bar{Q}_n$
weakly converges. Since every divisor $D$ of bounded degree
is equivalent to the difference $D\sim D'-D''$ of two effective
divisors $D'$ and $D''$ of bounded degrees
(compare \cite[Theorem~21.7.]{Fo}) a subsequence of the sequence
of divisors $D_n$ is equivalent to a convergent sequence of divisors
with limit $D$. By passing to an equivalent subsequence we achieve
that the sequence of divisors $D_n$ converges to the divisor $D$.

We cover $\Spa{X}$ by open subsets
$$\Spa{X}=\Set{U}_0\cup\ldots\cup\Set{U}_L.$$
Here $\Set{U}_0$ is the complement of the union
of small neighbourhoods of the  support of the divisor $D$ with the
support of the divisor $D_{\text{\scriptsize\rm spin}}$
of the spin bundle used in Theorem~\ref{weakly continuous resolvent}
and all those points, whose mass with respect to the weak limit of the
measure $\frac{1}{2\sqrt{-1}}Q_n\wedge\Bar{Q}_n$ is greater or equal
than the constant $S_p^{-1}$. The other sets
$\Set{U}_1,\ldots,\Set{U}_L$ are small open disjoint disks,
which cover the connected components of the complement of $\Set{U}_0$.
Since the holomorphic structures of $\Sh{Q}_{D_{\text{\tiny\rm spin}}}$
are Dirac operators with potentials, whose resolvents are
investigated in Theorem~\ref{weakly continuous resolvent},
the restrictions of the holomorphic structures to $\Set{U}_0$ is also
of this form. Due to Theorem~\ref{weakly continuous resolvent} and
Lemma~\ref{bounded point measures} the resolvents of these
restrictions of the homomorphic structures to $\Set{U}_0$ converges.
By subtracting from $\Set{U}_0$ additional small closed disks
contained in additional open sets $\Set{U}_{L+1},\ldots,\Set{U}_{L'}$,
which are disjoint from $\Set{U}_1,\ldots,\Set{U}_L$ and form each other,
we may achieve that the corresponding limit of the sequence of
restrictions of the holomorphic structures to $\Set{U}_0$ has a
resolvent. Due to Theorem~\ref{weakly continuous resolvent} these
restrictions of holomorphic structures have always reduced resolvents
on the complement of a finite--dimensional subspace of holomorphic
sections. Our arguments in step~{5}, where we prove the existence
of convergent subsequences can be extended to this more general
situation, since all bounded subsets of these finite--dimensional
subsets are compact.

\noindent{\bf 3. Limits of the local resolvents near the singular
points with trivial kernels of the blown up holomorphic structures.}
In this step we consider the limits of the restrictions of the
holomorphic structures to $\Set{U}_1,\ldots,\Set{U}_L$.
We assume that local parameters maps these small open disks onto small
open domains in $\mathbb{C}$. Therefore the restrictions of the
holomorphic structures to $\Set{U}_1,\ldots,\Set{U}_L$ can be described
by Dirac operators with potentials $\left(\begin{smallmatrix}
U & \partial\\
-\Bar{\partial} & \Bar{U}
\end{smallmatrix}\right)$ on open sets of $\mathbb{C}$.
If $\Set{U}_l$ does not contain a point, whose mass with respect to
the weak limit of the measures
$\frac{1}{2\sqrt{-1}}Q_n\wedge\Bar{Q}_n$ is greater or equal than
$S_p^{-1}$, then due to Theorem~\ref{weakly continuous resolvent} and
Lemma~\ref{bounded point measures} the resolvents
of the restrictions of the holomorphic structures to $\Set{U}_0$
converges to the resolvent of the holomorphic structure,
whose Hopf field is the weak limit.

Let us assume that the support of the sequence of Hopf fields
$Q_{\text{\scriptsize\rm sing},n,l}=
-\Bar{U}_{\text{\scriptsize\rm sing},n,l}d\Bar{z}$
is contained in $\Set{U}_l$,
and that the sequence of holomorphic structures with Hopf fields
$\left(g_{n,1}^{-1}\right)^{\ast}Q_{\text{\scriptsize\rm sing},n,1}$
on $\mathbb{P}$ has a trivial kernel. We claim, that in this case
the corresponding sequence of resolvents of Dirac operators on
$\Set{U}_l$, whose Hopf fields are given by the restrictions of the
Hopf fields 
$$Q_{\text{\scriptsize\rm reg},n}+Q_{\text{\scriptsize\rm sing},n,l}=
-\Bar{U}_{\text{\scriptsize\rm reg},n}d\Bar{z}-
\Bar{U}_{\text{\scriptsize\rm sing},n,l}d\Bar{z}$$
to $\Set{U}_l$, considered as operators from
$\banach{p}(\Set{U}_l)\times\banach{p}(\Set{U}_l)$ into
$\banach{q}(\Set{U}_l)\times\banach{q}(\Set{U}_l)$ with $1<p<2$ and
$1<q<\frac{2p}{2-p}$ converges to the resolvent of the
Dirac operator, whose potential corresponds to the weak limit of the
sequence of Hopf fields.
The corresponding resolvents obey the relation
\begin{multline*}
\Op{R}_{\mathbb{C}}\left(U_{\text{\scriptsize\rm reg},n}+
U_{\text{\scriptsize\rm sing},n,l},\Bar{U}_{\text{\scriptsize\rm reg},n}+
\Bar{U}_{\text{\scriptsize\rm sing},n,l},0\right)=\\
=\Op{R}_{\mathbb{C}}\left(
U_{\text{\scriptsize\rm reg},n},\Bar{U}_{\text{\scriptsize\rm reg},n},
0\right)\comp\left(\unity-\begin{pmatrix}
U_{\text{\scriptsize\rm sing},n,l} & 0\\
0 & \Bar{U}_{\text{\scriptsize\rm sing},n,l}
\end{pmatrix}\comp\Op{R}_{\mathbb{C}}\left(
U_{\text{\scriptsize\rm reg},n},\Bar{U}_{\text{\scriptsize\rm reg},n},
0\right)\right)^{-1}\\
=\Op{R}_{\mathbb{C}}\left(
U_{\text{\scriptsize\rm reg},n},\Bar{U}_{\text{\scriptsize\rm reg},n},
0\right)+\Op{R}_{\mathbb{C}}\left(
U_{\text{\scriptsize\rm reg},n},\Bar{U}_{\text{\scriptsize\rm reg},n},
0\right)\comp\\
\comp\left(\unity-\begin{pmatrix}
U_{\text{\scriptsize\rm sing},n,l} & 0\\
0 & \Bar{U}_{\text{\scriptsize\rm sing},n,l}
\end{pmatrix}\comp\Op{R}_{\mathbb{C}}\left(
U_{\text{\scriptsize\rm reg},n},\Bar{U}_{\text{\scriptsize\rm reg},n},
0\right)\right)^{-1}\comp\begin{pmatrix}
U_{\text{\scriptsize\rm sing},n,l} & 0\\
0 & \Bar{U}_{\text{\scriptsize\rm sing},n,l}
\end{pmatrix}\comp\Op{R}_{\mathbb{C}}\left(
U_{\text{\scriptsize\rm reg},n},\Bar{U}_{\text{\scriptsize\rm reg},n},
0\right).
\end{multline*}
The operators
\begin{eqnarray*}
&&\left(\unity-\begin{pmatrix}
U_{\text{\scriptsize\rm sing},n,l} & 0\\
0 & \Bar{U}_{\text{\scriptsize\rm sing},n,l}
\end{pmatrix}\comp\Op{R}_{\mathbb{C}}\left(
U_{\text{\scriptsize\rm reg},n},\Bar{U}_{\text{\scriptsize\rm reg},n},
0\right)\right)^{-1}\comp\begin{pmatrix}
U_{\text{\scriptsize\rm sing},n,l} & 0\\
0 & \Bar{U}_{\text{\scriptsize\rm sing},n,l}
\end{pmatrix}\\
&=&\begin{pmatrix}
U_{\text{\scriptsize\rm sing},n,l} & 0\\
0 & \Bar{U}_{\text{\scriptsize\rm sing},n,l}
\end{pmatrix}\comp\left(\unity-\Op{R}_{\mathbb{C}}\left(
U_{\text{\scriptsize\rm reg},n},\Bar{U}_{\text{\scriptsize\rm reg},n},
0\right)\comp\begin{pmatrix}
U_{\text{\scriptsize\rm sing},n,l} & 0\\
0 & \Bar{U}_{\text{\scriptsize\rm sing},n,l}
\end{pmatrix}\right)^{-1}
\end{eqnarray*}
depend only on the restrictions of $\Op{R}_{\mathbb{C}}\left(
U_{\text{\scriptsize\rm reg},n},\Bar{U}_{\text{\scriptsize\rm reg},n},
0\right)$ to the support of $U_{\text{\scriptsize\rm sing},n,l}$.
We shall transform this sequence of operators under the corresponding
sequence of M\"obius transformations $g_{n,l}$. The small open sets
$\Set{U}_1,\ldots,\Set{U}_L$ are identified with bounded open sets of
$\mathbb{C}$. Therefore the restrictions of the holomorphic
structures may be described by Dirac operators with potentials on
bounded open sets of $\mathbb{C}$.
All M\"obius transformations $h$ induce isometries
\begin{align*}
\Op{I}_p(h)&:&\banach{p}(\mathbb{C})&\rightarrow\banach{p}(\mathbb{C})&
f&\mapsto \Tilde{f}&\tilde{f}(z)&=
f(h^{-1}z)\left|\frac{dh^{-1}z}{dz}\right|^{\frac{2}{p}}.
\end{align*}
A direct calculation shows that the resolvent
$\Op{R}_{\mathbb{C}}(0,0,0)$
of the free Dirac operator, considered as an operator from
$\banach{p}(\mathbb{C})\times\banach{p}(\mathbb{C})$ into
$\banach{\frac{2p}{2-p}}(\mathbb{C})\times
\banach{\frac{2p}{2-p}}(\mathbb{C})$ with $1<p<2$ is invariant under
the scaling transformations ($z\mapsto\exp(t)z$ with $t\in\mathbb{R}$)
and the translations ($z\mapsto z+z_0$ with $z_0\in\mathbb{C}$).
Since these sequences of M\"obius transformations belong to the
semidirect product of the scaling transformations
with the translations, the free resolvent is invariant under
these transformations $g_{n,l}$:
$$\begin{pmatrix}
\Op{I}_{\frac{2p}{2-p}}(g_{n,l}) & 0\\
0 & \Op{I}_{\frac{2p}{2-p}}(g_{n,l})
\end{pmatrix}\comp\Op{R}_{\mathbb{C}}(0,0,0)\comp\begin{pmatrix}
\Op{I}_p(g_{n,l}^{-1}) & 0\\
0 & \Op{I}_p(g_{n,l}^{-1})
\end{pmatrix}=\Op{R}_{\mathbb{C}}(0,0,0).$$

If the sets $\Set{U}_1,\ldots,\Set{U}_L$ are small, then the
restrictions of the sequence of transformed Hopf fields
$\left(g_{n,l}^{-1}\right)^{\ast}Q_n$ to the subset
$g_{n,l}^{-1}\left(\Set{U}_{l}\right)\subset\mathbb{P}$
still obey the condition of Lemma~\ref{bounded point measures}.
Therefore, due to Theorem~\ref{weakly continuous resolvent},
the corresponding sequence of resolvents on $\mathbb{P}$ converges.
We assume that the limit is the resolvents of a Dirac operator on
$\mathbb{P}$ without kernel. In this case the arguments of
Theorem~\ref{weakly continuous resolvent}
together with the first resolvent formula \cite[Theorem~VI.5]{RS1}:
$$(\lambda-\lambda')\Op{R}_{\lambda'}
=\Op{R}_{\lambda}\comp
\left(\frac{\unity}{\lambda-\lambda'}-\Op{R}_{\lambda}\right)^{-1}
=\left(\frac{\unity}{\lambda-\lambda'}-\Op{R}_{\lambda}\right)^{-1}
\comp\Op{R}_{\lambda},$$
imply that the corresponding sequence of resolvents considered as
operators from $\banach{p}(\mathbb{C})\times\banach{p}(\mathbb{C})$
into $\banach{\frac{2p}{2-p}}(\mathbb{C})\times
\banach{\frac{2p}{2-p}}(\mathbb{C})$ is bounded.
We conclude that the sequences of operators
\begin{eqnarray*}
&&\left(\unity-\begin{pmatrix}
U_{\text{\scriptsize\rm sing},n,l} & 0\\
0 & \Bar{U}_{\text{\scriptsize\rm sing},n,l}
\end{pmatrix}\comp\Op{R}_{\mathbb{C}}\left(
U_{\text{\scriptsize\rm reg},n},\Bar{U}_{\text{\scriptsize\rm reg},n},
0\right)\right)^{-1}\comp\begin{pmatrix}
U_{\text{\scriptsize\rm sing},n,l} & 0\\
0 & \Bar{U}_{\text{\scriptsize\rm sing},n,l}
\end{pmatrix}\\
&=&\begin{pmatrix}
U_{\text{\scriptsize\rm sing},n,l} & 0\\
0 & \Bar{U}_{\text{\scriptsize\rm sing},n,l}
\end{pmatrix}\comp\left(\unity-\Op{R}_{\mathbb{C}}\left(
U_{\text{\scriptsize\rm reg},n},\Bar{U}_{\text{\scriptsize\rm reg},n},
0\right)\comp\begin{pmatrix}
U_{\text{\scriptsize\rm sing},n,l} & 0\\
0 & \Bar{U}_{\text{\scriptsize\rm sing},n,l}
\end{pmatrix}\right)^{-1}
\end{eqnarray*}
are bounded.

Due to H\"older's inequality \cite[Theorem~III.1~(c)]{RS1}
for $1\leq q'<q\leq\infty$ the restriction of
$\banach{q}(\mathbb{C})$ into $\banach{q'}(B(0,\varepsilon))$
is bounded by
$\left(\pi\varepsilon^2\right)^{\frac{1}{q'}-\frac{1}{q}}$.
Since the radii of the supports of $U_{\text{\scriptsize\rm sing},n,l}$
tend to zero, the restrictions of the resolvents
$\Op{R}_{\mathbb{C}}(U_n,\Bar{U}_n,0)$ considered as operators from
$\banach{p}(\Set{U}_l)\times\banach{p}(\Set{U}_l)$ into
$\banach{q}(\Set{U}_l)\times\banach{q}(\Set{U}_l)$
with $1<p<2$ and $1<q<\frac{2p}{2-p}$ converge to the resolvent of the
weak limit of the sequences $U_n$.

\noindent{\bf 4. Limits of the local resolvents near the singular
points with non--trivial kernels of the blown up holomorphic structures.}
In this case we add to the sequence of divisors $D_n$ a sequence of effective
divisors $D'_n$ with support in the complements of $\Set{U}_0$,
such that the corresponding transformed sequences of holomorphic
structures corresponding with Hopf fields
$\left(g_{n,1}^{-1}\right)^{\ast}Q_{\text{\scriptsize\rm sing},n,1},\ldots,
\left(g_{n,L'}^{-1}\right)^{\ast}Q_{\text{\scriptsize\rm sing},n,L'}$
on $\mathbb{P}$ have trivial kernels.

\begin{Lemma}\label{transformed holomorphic structures}
For any holomorphic quaternionic line bundle on $\mathbb{P}$
with non--trivial kernel let $d$ be the unique natural number such that
\begin{align*}
\dim H^0\left(\mathbb{P},\Sh{Q}_{D-d\infty}\right)&=0&
\dim H^0\left(\mathbb{P},\Sh{Q}_{D-(d-1)\infty}\right)&=1
\end{align*}
Then there exists an effective divisor $D'$ of degree $d-\deg(D)-1$,
whose support is contained in $\mathbb{C}$,
such that
\begin{align*}
\dim H^0\left(\mathbb{P},\Sh{Q}_{D+D'+(n-d)\infty}\right)&=n
&\forall n&\in\mathbb{N}_0.
\end{align*}
\end{Lemma}

\begin{proof}
Due to the \De{Riemann--Roch Theorem}~\ref{riemann roch}
we have the inequality
$$\dim H^0\left(\mathbb{P},\Sh{Q}_{D-d\infty}\right)= \deg(D)+1-d+
\dim H^1\left(\mathbb{P},\Sh{Q}_{D-d\infty}\right)\geq \deg(D)+1-d.$$
By definition of $d$ this implies $\deg(D)\leq d-1$.
Moreover, the equality $\deg(D)=d-1$ is equivalent to
$\dim H^1\left(\mathbb{P},\Sh{Q}_{D-d\infty}\right)=0$.
Due to \De{S\'{e}rre Duality}~\ref{serre duality} this is equivalent to
\begin{align*}
\dim H^1\left(\mathbb{P},\Sh{Q}_{D+(n-d)\infty}\right)&=0&
\forall n&\in\mathbb{N}_0.
\end{align*}
Finally, due to the \De{Riemann--Roch Theorem}~\ref{riemann roch},
the equality $\deg(D)=d-1$ is equivalent to
\begin{align*}
\dim H^0\left(\mathbb{P},\Sh{Q}_{D+(n-d)\infty}\right)&=n
&\forall n&\in\mathbb{N}_0.
\end{align*}
Therefore it suffices to consider the cases $\deg(D)<d-1$.

We claim that in this case there exists an element
$z\in\mathbb{C}$, such that the analogous number $d$
corresponding to the holomorphic structure of $\Sh{Q}_{D+z}$
is equal to $d$. This is equivalent to
$\dim H^0\left(\mathbb{P},\Sh{Q}_{D+z-d\infty}\right)=0$.
Let us assume on the contrary that for all $z\in\mathbb{C}$ we have
$\dim H^0\left(\mathbb{P},\Sh{Q}_{D+z-d\infty}\right)=1$.
Consequently, for all pairwise different $z_1,\ldots,z_L\in\mathbb{C}$
the dimension of the linear system
$H^0\left(\mathbb{P},\Sh{Q}_{D+z_1+\ldots+z_L-d\infty}\right)$
is larger than $L$. For large $L$ due to
\De{S\'{e}rre Duality}~\ref{serre duality},
\De{Pl\"ucker formula}~\ref{pluecker formula}
the \u{C}ech cohomology group
$H^1\left(\mathbb{P},\Sh{Q}_{D+z_1+\ldots+z_L-d\infty}\right)$
is trivial. Consequently, due to
\De{Riemann--Roch Theorem}~\ref{riemann roch}, we obtain
$$L\leq H^0\left(\mathbb{P},\Sh{Q}_{D+z_1+\ldots+z_L-d\infty}\right)=
1+\deg(D)+L-d,$$
which contradicts to $\deg(D)<d-1$. This proves the claim.

By an iterated application of this claim we obtain
an effective divisor $D'$ with the desired properties.
\end{proof}

We apply this lemma to the holomorphic structure
corresponding to the weak limits of Hopf fields
$\left(g_{n,1}^{-1}\right)^{\ast}Q_{\text{\scriptsize\rm sing},n,1},\ldots,
\left(g_{n,L'}^{-1}\right)^{\ast}Q_{\text{\scriptsize\rm sing},n,L'}$
on $\mathbb{P}$. Since we are only interested in the restrictions of
the holomorphic structure to $\Set{U}_1,\ldots,\Set{U}_L$, we may
change the degree at $\infty$.
For all holomorphic quaternionic line bundles on $\mathbb{P}$,
with sheaf $\Sh{Q}_{D}$ of holomorphic sections,
the sheaf of holomorphic sections $\Sh{Q}_{D+D'-d\infty}$ of the corresponding
holomorphic structure on the spin bundle has a trivial kernel.
Here $D'$ denotes the divisor of degree $d-\deg(D)-1$ constructed in
Lemma~\ref{transformed holomorphic structures}.
Obviously, the sequence of divisors $D'_n=g_{n,l}(D')$ converge to the
divisor $\deg(D')x_l$ on $\Set{U}_l$. Hence the arguments of
step~{4} imply that the resolvents of the corresponding
Dirac operators on $\Set{U}_l$ converges to the resolvent of the
Dirac operator, whose potential is the weak limit of the sequence
of potentials.

\noindent{\bf 5. Limits of the sequence of holomorphic sections.}
In the preceding step we added to the sequence of divisors a sequence
of convergent effective divisors. Obviously, any sequence of sections
of the original sequence of holomorphic quaternionic line bundles are
also holomorphic sections of the latter sequence of holomorphic
quaternionic line bundles. We shall prove that this sequence converges
to a non--trivial section of the limit of the latter sequence of
holomorphic quaternionic line bundles. More precisely, the Hopf field
of the limit of the holomorphic structures is the weak limit of the
sequence of Hopf fields of the latter sequence of Hopf fields.

At the end of step~{2} we saw that the sequence of resolvents of
the restrictions of the holomorphic structures to $\Set{U}_0$
converged as an operator from
$H^0\left(\Set{U}_0,\Sh{L}^p_{D_{\text{\tiny\rm spin}}}\right)$ into
$H^0\left(\Set{U}_0,\Sh{L}^q_{D_{\text{\tiny\rm spin}}}\right)$ with
$1<p<2$ and $1<q<\frac{2p}{2-p}$. Moreover, in steps~{3}--{5}
we showed that for all $l=1,\ldots,L$ the resolvents of the
restriction of the holomorphic structures to $\Set{U}_l$
converged as an operator from
$H^0\left(\Set{U}_0,\Sh{L}^p_{D_{\text{\tiny\rm spin}}}\right)$ into
$H^0\left(\Set{U}_0,\Sh{L}^q_{D_{\text{\tiny\rm spin}}}\right)$ with
$1<p<2$ and $1<q<\frac{2p}{2-p}$.

Due to quaternionic version of
\De{Cauchy's Integral Formula}~\ref{cauchy formula}
the holomorphic sections are uniquely determined by their
restrictions to
$$\left(\Set{U}_1\cap\Set{U}_0\right)\cup\ldots\cup
\left(\Set{U}_l\cap\Set{U}_0\right).$$
Moreover, for all $1<p<2$ and $1<q<\frac{2p}{2-p}$
the $H^0\left(\Spa{X},\Sh{L}^q_{D_n}\right)$--norms are uniformly
bounded in terms of the
$H^0\left(\Spa{X},\Sh{L}^p_{D_n}\right)$--norms
of the restrictions to
$\left(\Set{U}_1\cap\Set{U}_0\right)\cup\ldots\cup
\left(\Set{U}_l\cap\Set{U}_0\right)$.
Since the sequence of Hopf fields is bounded, the
$H^0\left(\Spa{X},\Sh{W}^{1,p}_{D_n}\right)$--norms are bounded
uniformly in terms of the
$H^0\left(\Spa{X},\Sh{L}^{\frac{2p}{2-p}}_{D_n}\right)$--norms.
Now Kondrakov's Theorem \cite[Theorem~2.34]{Au} implies
that any sequence of non--trivial eigenfunctions, whose
$H^0\left(\Spa{X},\Sh{L}^q_{D_n}\right)$--norms are equal to one,
have a convergent subsequence, and that the limit is non--trivial.
Due to the convergence of the resolvents in steps~{2}--{4}
the limit is holomorphic with respect to the holomorphic structures,
whose Hopf field is the weak limit of the sequence of Hopf fields.
\end{proof}

\section{Existence of minimizers}\label{section minimizers}
In this section we prove the existence of minimizing surfaces in
$\mathbb{R}^3$ and $\mathbb{R}^4$ of the
Willmore functional inside all conformal classes.
More precisely, we show that any sequence of conformal mappings from a
compact Riemann surface $\Spa{X}$ into $\mathbb{R}^3$
(or $\mathbb{R}^4$), whose Willmore functionals is bounded,
may be transformed by a sequence of conformal mappings of
$\mathbb{R}^3\subset\Spa{S}^3$ (or $\mathbb{R}^4\subset\Spa{S}^4$)
into a sequence, which converges with respect to the
$\sobolev{2,p}(\Spa{X})$--topology for all $1<p<2$.\
Essentially this follows from the
\De{Quaternionic Weierstra{\ss} Representation}~\ref{weierstrass}
and Theorem~\ref{convergent subsequences}.
In \cite{PP} the global Weierstra{\ss} representation was generalized
to conformal mappings into $\mathbb{R}^4$.
In fact, the `quaternionic function theory' provides two version of
a global Weierstra{\ss} representation into $\mathbb{R}^4$. From our
point of view they are related by a
\De{B\"acklund transformation}~\ref{baecklund}.

\begin{Proposition}\label{compactness of conformal mappings}
For any sequence of mappings in one of the following classes
there exists a sequence of conformal transformations of the
target space, such that the transformed sequence
has a convergent subsequence with respect to the topologies of
$\bigcap\limits_{1<p<2}\sobolev{2,p}(\Spa{X})$:
\begin{description}
\item[(i)] Smooth conformal mappings from a compact Riemann surface
$\Spa{X}$ into $\mathbb{R}^3$ with bounded Willmore functional.
\item[(ii)] Smooth conformal mappings from a compact Riemann surface
$\Spa{X}$ into $\mathbb{R}^4$ with bounded Willmore functional.
\item[(iii)] Smooth conformal mappings from a compact Riemann surface
$\Spa{X}$ into $\mathbb{R}^4$ with a fixed complex holomorphic line bundle
underlying the quaternionic holomorphic line bundle
(compare \cite[Theorem~4.3]{PP}) and bounded Willmore functional.
\end{description}
\end{Proposition}

\begin{proof}
We use the
\De{Quaternionic Weierstra{\ss} Representation}~\ref{weierstrass}
\cite{PP,BFLPP} and its reduction to conformal mappings
into the pure imaginary quaternions $\simeq\mathbb{R}^3$
\cite{Ta1,Ta2,Fr2}.
Hence all immersion are represented by two non--trivial spinors of two
paired holomorphic quaternionic line bundles.
Let $\phi_n$ and $\psi_n$ be the sequences of paired spinors
corresponding to a minimizing sequence of the Willmore functional
on the space of conformal immersions of a compact Riemann surface
$\Spa{X}$ into $\mathbb{R}^3$ or $\mathbb{R}^4$.
The \De{Pl\"ucker formula}~\ref{pluecker formula} implies
that the degrees of the corresponding
quaternionic holomorphic line bundles are bounded from below.
Since these two line bundles are paired, the degrees are also bounded
from above. Therefore  Theorem~\ref{convergent subsequences} implies
that both sequences have convergent subsequences. But it might happen
that the corresponding limits of the holomorphic structures have
singularities. In this case the degrees of the limits are not the
limits of the degrees. The corresponding Weierstra{\ss}
representations describe immersions of $\Spa{X}$ into
$S^3\supset\mathbb{R}^3$ or $S^4\supset\mathbb{R}^4$.
But a conformal transformation of
$S^3\supset\mathbb{R}^3$ or $S^4\supset\mathbb{R}^4$
transforms these immersions into immersions into
$\mathbb{R}^3$ or $\mathbb{R}^4$.
We remark that the conformal transformations of
$S^4\supset\mathbb{R}^4$ are very easy to describe with the help of the
\De{B\"acklund transformation}~\ref{baecklund}.
In fact the conformal transformations are just equal to the action of
$GL(2,\mathbb{Q})$ on the corresponding quaternionic
two--dimensional subspace of holomorphic sections of the
B\"acklund transformed quaternionic holomorphic line bundle.
Observe that in Theorem~\ref{convergent subsequences} we implicitly
use translations and rotations of the immersions corresponding to
rescalings of the two holomorphic spinors of the two paired
quaternionic holomorphic line bundles.
If we use in addition some inversions, we may always achieve that the
limit stays inside of $\mathbb{R}^3$ or $\mathbb{R}^4$.
The corresponding two limits of $\phi_n$ and $\psi_n$ does not have poles.
Consequently the quaternionic holomorphic
line bundles have degrees equal to the limits of the corresponding
sequences of degrees and they are paired. In case the underlying
holomorphic complex line bundles are fixed, the limits of the
quaternionic holomorphic line bundles have also these underlying
complex holomorphic line bundles.
\end{proof}

We consider this Proposition as Montel's Theorem of
`quaternionic function theory'. It implies the existence of minimizers
of the Willmore functional.

\begin{Theorem}
The Willmore functional attains a minimum on the following classes:
\begin{description}
\item[(i)] Smooth conformal mappings from a compact Riemann surface
into $\mathbb{R}^3$.
\item[(ii)] Smooth conformal mappings from a compact Riemann surface
into $\mathbb{R}^4$.
\item[(iii)] Smooth conformal mappings from a compact Riemann surface
into $\mathbb{R}^4$ with a fixed complex holomorphic line bundle
underlying the quaternionic holomorphic line bundle
(compare \cite[Theorem~4.3]{PP}).
\end{description}
\end{Theorem}

\begin{proof}
Proposition~\ref{compactness of conformal mappings}
implies the convergence of a minimizing sequence
in the enlarged classes (i)--(iii)
of not necessarily smooth conformal mappings
with bounded Willmore functional.
It remains to ensure the smoothness of the minimizers.

\begin{Lemma}\label{local minimizer}
Let $\psi$ belong to the kernel of $\left(\begin{smallmatrix}
\Bar{\partial} & -\Bar{U}\\
U & \partial
\end{smallmatrix}\right)$ and $\phi$ to the kernel of
$\left(\begin{smallmatrix}
\Bar{\partial} & U\\
-\Bar{U} & \partial
\end{smallmatrix}\right)$ with potential $U\in\banach{2}(\Omega)$
on an open domain $\Omega\subset\mathbb{C}$.
If the Willmore functional $\willmore=4\int\limits_{\Omega}U\Bar{U}d^2x$
is minimal with respect to all $\banach{2}$--perturbations $\Var U$
with compact support in $\Omega$, which admit a perturbation of
$\psi$ and $\phi$ with \Em{compact support} in $\Omega$,
then there exists spinors $\Tilde{\psi}$ in the kernel of
$\left(\begin{smallmatrix}
\Bar{\partial} & -\Bar{U}\\
U & \partial
\end{smallmatrix}\right)$ and $\Tilde{\phi}$ in the kernel of
$\left(\begin{smallmatrix}
\Bar{\partial} & U\\
-\Bar{U} & \partial
\end{smallmatrix}\right)$, such that $(U,\Bar{U})$
is a complex linear combination of
$\left(\Tilde{\phi}_2\psi_1+\overline{\Tilde{\phi}_1\psi_2},
\Tilde{\phi}_1\psi_2+\overline{\Tilde{\phi}_2\psi_1}\right)$ and
$\left(\phi_2\Tilde{\psi}_1+\overline{\phi_1\Tilde{\psi}_2},
\phi_1\Tilde{\psi}_2+\overline{\phi_2\Tilde{\psi}_1}\right)$.
\end{Lemma}

\begin{proof}
If the support of $\Var U$ is contained in the open subdomain
$\Omega'\subset\overline{\Omega}'\subset\Omega$,
then, due to the quaternionic version of
\De{Cauchy's Integral Formula}~\ref{cauchy formula},
the restriction of $\phi$ and $\psi$ to the complement of
the domain of $\overline{\Omega'}$ in $\Omega$
are uniquely determined by their values on a cycle around $\Omega$
and on a cycle in $\Omega'$ around the support of $U$.
We conclude that $\psi$ and $\phi$ admit
perturbations with support contained in $\Omega'$, if and only if
the total residue with the corresponding integral kernels vanishes
on $\Omega'$. Again due to the quaternionic version of
\De{Cauchy's Integral Formula}~\ref{cauchy formula} this is equivalent
to the condition that for all elements $\tilde{\psi}$ in the kernel of
$\left(\begin{smallmatrix}
\Bar{\partial} & -\Bar{U}-\Var\Bar{U}\\
U+\Var U & \partial
\end{smallmatrix}\right)$ and all $\Tilde{\phi}$
in the kernel of $\left(\begin{smallmatrix}
\Bar{\partial} & U+\Var U\\
-\Bar{U}-\Var\Bar{U} & \partial
\end{smallmatrix}\right)$ the residues of the forms
$\Tilde{\phi}^{t}\left(\begin{smallmatrix}
d\Bar{z} & 0\\
0 & dz
\end{smallmatrix}\right)\psi$ and $\phi^{t}\left(\begin{smallmatrix}
d\Bar{z} & 0\\
0 & dz
\end{smallmatrix}\right)\Tilde{\psi}$ on $\Omega'$ vanish.
Due to the equations
\begin{eqnarray*}
d\left(\Tilde{\phi}^{t}\begin{pmatrix}
d\Bar{z} & 0\\
0 & dz
\end{pmatrix}\psi\right)&=&\Tilde{\phi}^{t}\begin{pmatrix}
0 & \Var\Bar{U}\\
\Var U & 0
\end{pmatrix}\psi dz\wedge d\Bar{z}\\
d\left(\phi^{t}\begin{pmatrix}
d\Bar{z} & 0\\
0 & dz
\end{pmatrix}\Tilde{\psi}\right)&=&\phi^{t}\begin{pmatrix}
0 & \Var\Bar{U}\\
\Var U & 0
\end{pmatrix}\Tilde{\psi}dz\wedge d\Bar{z}
\end{eqnarray*}
this is equivalent to the equations
\begin{align*}
\int\limits_{\Omega}
\Tilde{\phi}^{t}\begin{pmatrix}
0 & \sqrt{-1}\Var\Bar{U}\\
\sqrt{-1}\Var U & 0
\end{pmatrix}\psi d^2x&=0&
\int\limits_{\Omega}
\phi^{t}\begin{pmatrix}
0 & \sqrt{-1}\Var\Bar{U}\\
\sqrt{-1}\Var U & 0
\end{pmatrix}\Tilde{\psi}d^2x&=0.
\end{align*}
We shall apply the implicit function theorem and conclude that the
space of perturbations $\Var U$, which admit perturbations
of $\psi$ and $\phi$ with compact support are submanifolds.
Since the question is local, we may chose the domain $\Omega'$ to be
the unit disk $\mathbb{D}$. Indeed, appropriate small neighbourhoods
of any point are M\"obius transforms of $\mathbb{D}$.
On $\mathbb{D}$ we introduce the Banach spaces
$\banach{q}\left(\mathbb{D},(1-|z|^2)^{-s}d^2x,\mathbb{H}\right)$
of quaternionic valued $\banach{q}$--functions with respect
to the measure $\left(1-|z|^2\right)^{-s}d^2x$ on the unit disk
$\mathbb{D}\subset\mathbb{C}$ with $0\leq s<1$.

As a preparation we claim that the kernel of $\left(\begin{smallmatrix}
\Bar{\partial} & -\Bar{U}\\
U & \partial
\end{smallmatrix}\right)$ considered as a closed subspace of
$\banach{p}\left(\mathbb{D},(1-|z|^2)^{-s}d^2x,\mathbb{H}\right)$
is contained in
$\bigcap\limits_{q<\frac{2p}{2-s}}\banach{q}(\mathbb{D},\mathbb{H})$.
For the proof we apply the quaternionic version of
\De{Cauchy's Integral Formula}~\ref{cauchy formula}.
Due to Lemma~\ref{integral kernel bound} it suffices to
show that the integral along the boundary of $\mathbb{D}$ over
the integral kernel of $\Op{I}_{\mathbb{D}}(0)$
defines a bounded operator from
$\banach{p}\left(\mathbb{D},(1-|z|^2)^{-s}d^2x,\mathbb{H}\right)$,
into $\banach{\frac{2p}{2-s}}(\mathbb{D},\mathbb{H})$.
Due to Young's inequality \cite[Section~IX.4 Example~1]{RS2}
the convolution with the function $\frac{1}{\pi z}$
defines an operator from the $\banach{p}$--functions on the circle
$|z'|=r'$ into the $\banach{q}$--functions on the circle
$|z|=r$ with $0\leq r<r'\leq1$, which is bounded by
$$\frac{1}{\pi}\left(\int\limits_{\varphi\in\mathbb{R}/2\pi\mathbb{Z}}
\frac{1}{\left|1-\frac{r}{r'}
\exp\left(2\pi\sqrt{-1}\varphi\right)\right|^{\frac{pq}{p+pq-q}}}
\right)^{1+\frac{1}{q}-\frac{1}{p}},$$
with $1\leq p<q\leq\infty$.
Due to \cite[1.4.10.Proposition]{Ru} this norm is bounded by
$C\left|1-\frac{r^2}{r'^2}\right|^{\frac{1}{q}-\frac{1}{p}}$.
The norm of the restriction of this function to
$r\in[0,r_0]\subset[0,1]$ in the $\banach{q}$--space on $r\in[0,1]$
with respect to the measure $rdr$ is bounded by
$C'\left|r'-r_0\right|^{\frac{2}{q}-\frac{1}{p}}+C''$
with appropriate constants $C'>0$ and $C''>0$.
If $p'$ denotes the dual exponent of $p$ with
$\frac{1}{p}+\frac{1}{p'}=1$, then,
due to H\"older's inequality \cite[Theorem~III.1~(c)]{RS1},
we obtain for all $f\in\banach{p}\left([0,1],\frac{r'dr'}{(1-r'^2)^s}\right)$
\begin{eqnarray*}
\frac{1}{1-r_0}\int\limits_{r_0}^1f(r')
\left|r'-r_0\right|^{\frac{2}{q}-\frac{1}{p}}r'dr'
&\leq&\frac{
\left\|f\right\|_{\banach{p}\left([0,1],\frac{r'dr'}{(1-r'^2)^s}\right)}}
           {1-r_0}
\left\|\left|1-r'^2\right|^s
\left|r'-r_0\right|^{\frac{2}{q}-\frac{1}{p}}
\right\|_{\banach{p'}([r_0,1],r'dr')}\\
&\leq&
\left\|f\right\|_{\banach{p}\left([0,1],\frac{r'dr'}{(1-r'^2)^s}\right)}
\text{\bf{O}}\left(1-r_0\right)^{\frac{s}{p}+\frac{2}{q}-\frac{2}{p}}.
\end{eqnarray*}
With $1\leq q\leq\frac{2p}{2-s}$ this expression
remains bounded in the limit $r_0\rightarrow 1$.
Consequently, for $1\leq q<\frac{2p}{2-s}$
the natural inclusion of the kernel of
$\left(\begin{smallmatrix}
\Bar{\partial} & -\Bar{U}\\
U & \partial
\end{smallmatrix}\right)$ in
$\banach{p}\left(\mathbb{D},(1-|z|^2)^{-s}d^2x,\mathbb{H}\right)$,
into $\banach{q}(\mathbb{D},\mathbb{H})$ is bounded.

Let $\psi$ belong to the kernel of
$\left(\begin{smallmatrix}
\Bar{\partial} & -\Bar{U}\\
U & \partial
\end{smallmatrix}\right)$ without zeros on an open neighbourhood
$\Omega\subset\mathbb{C}$ of the closed unit disk $\Bar{\mathbb{D}}$
with a small $U\in\banach{2}(\mathbb{D})$.
In a second step we claim that for $2<q<\infty$
and $0<s<\frac{q-2}{q-1}$ the subset of all
$\Var U$ in the Banach space
$\banach{q}\left(\mathbb{D},(1-|z|^2)^{-s}d^2x\right)$,
such that
$$\int\limits_{\Omega}\Tilde{\phi}^{t}\begin{pmatrix}
0 & \sqrt{-1}\Var\Bar{U}\\
\sqrt{-1}\Var U & 0
\end{pmatrix}\psi d^2x=0$$
vanishes for all $\Tilde{\phi}$ in the kernel of
$\left(\begin{smallmatrix}
\Bar{\partial} & U+\Var U\\
-\Bar{U}-\Var\Bar{U} & \partial
\end{smallmatrix}\right)$ is a Banach submanifold.
In fact, due to the implicit function theorem \cite[Theorem~S.11]{RS1}
we have to show that for small $\Var U$ these kernels
considered as subspaces of the dual Banach spaces of
$\Var U\in\banach{q}\left(\mathbb{D},(1-|z|^2)^{-s}d^2x\right)$
with respect to the pairing
$$\left(\Tilde{\phi},\Var U\right)\mapsto
\int\limits_{\Omega}\Tilde{\phi}^{t}\begin{pmatrix}
0 & \sqrt{-1}\Var\Bar{U}\\
\sqrt{-1}\Var U & 0
\end{pmatrix}\psi d^2x$$
are isomorphic. Since $\psi^{-1}$ belongs to
$\bigcap\limits_{r<\infty}\banach{r}(\mathbb{D},\mathbb{H})=
\bigcap\limits_{r<\infty}
\banach{r}\left(\mathbb{D},(1-|z|^2)^{-s}d^2x,\mathbb{H}\right)$
these kernels are contained in $\bigcap\limits_{p<\frac{q}{q-1}}
\banach{p}\left(\mathbb{D},(1-|z|^2)^{-s}d^2x,\mathbb{H}\right)$.
The foregoing claim implies that these kernels are contained in
$\bigcap\limits_{p<\frac{2q}{(q-1)(2-s)}}
\banach{p}\left(\mathbb{D},\mathbb{H}\right)$.
The operator $\unity+\Op{I}_{\Omega}(U)\comp\left(\begin{smallmatrix}
0 & -\Var\Bar{U}\\
\Var U & 0
\end{smallmatrix}\right)$ maps these kernels onto the kernel of
$\left(\begin{smallmatrix}
\Bar{\partial} & U\\
-\Bar{U} & \partial
\end{smallmatrix}\right)$. If $q$ and $s$ satisfies
$\frac{2-s}{2}\left(1-\frac{1}{q}\right)+\frac{1}{q}-\frac{1}{2}+
s\left(1-\frac{1}{q}\right)<1-\frac{1}{q}$, then the operator
$\Op{I}_{\Omega}(U)\comp\left(\begin{smallmatrix}
0 & -\Var\Bar{U}\\
\Var U & 0
\end{smallmatrix}\right)$ is a bounded operator from
$\bigcap\limits_{p<\frac{2q}{(q-1)(2-s)}}
\banach{p}\left(\mathbb{D},\mathbb{H}\right)$ into
$\banach{\frac{q}{q-1}}\left(\mathbb{D},(1-|z|^2)^{-s}d^2x,\mathbb{H}\right)$.
Hence for $0<s<\frac{q-2}{q-1}$ these kernels are isomorphic.

Obviously, the same statement holds,
if $\psi$ is replaced by a spinor $\phi$ without zeroes in the kernel of
$\left(\begin{smallmatrix}
\Bar{\partial} & U\\
-\Bar{U} & \partial
\end{smallmatrix}\right)$ and $\Tilde{\phi}$ by spinors in the kernel of
$\left(\begin{smallmatrix}
\Bar{\partial} & -\Bar{U}-\Var\Bar{U}\\
U+\Var U & \partial
\end{smallmatrix}\right)$. Moreover, due to the considerations of
section~\ref{section local}, the intersection of theses two subspaces
of the dual of the Banach space
$\Var U\in\banach{q}\left(\mathbb{D},(1-|z|^2)^{-s}d^2x\right)$
is equal to the linear hull of $\Tilde{\phi}=\phi$
and $\Tilde{\psi}=\psi$. If $\psi$ is a spinor in the kernel of
$\left(\begin{smallmatrix}
\Bar{\partial} & -\Bar{U}\\
U & \partial
\end{smallmatrix}\right)$ and $\phi$ a
spinor in the kernel of
$\left(\begin{smallmatrix}
\Bar{\partial} & U\\
-\Bar{U} & \partial
\end{smallmatrix}\right)$ without zeroes on $\Omega$,
then  the subspace of all
$\Var U\in\banach{q}\left(\mathbb{D},(1-|z|^2)^{-s}d^2x\right)$,
which admit variations with compact support of $\psi$ and $\phi$,
are Banach submanifolds. Furthermore the tangent space of this
manifold is the orthogonal complement of these kernels with respect to
the corresponding pairings. This implies the statement of the Lemma
on the complement of the zeroes of $\psi$ and $\phi$.
Since a $\banach{2}$--functions, which vanishes on this complement,
vanishes on the whole of $\Omega$, the Lemma is proven.
\end{proof}

For all $m\in\mathbb{N}$ and $1<p<2$ the operator
$\Op{I}_{\Omega}(0)$ is a bounded operator from
$\sobolev{m-1,p}(\Omega,\mathbb{H})$ onto
$\sobolev{m,p}(\Omega,\mathbb{H})$ (compare \cite[Chapter~V]{St}).
Therefore the kernels of $\left(\begin{smallmatrix}
\Bar{\partial} & -\Bar{U}\\
U & \partial
\end{smallmatrix}\right)$ belongs to
$\bigcap\limits_{q<\infty}
\sobolev{m,p}_{\text{\scriptsize\rm loc}}(\Omega,\mathbb{H})$,
if the potential $U$ belongs to
$\bigcap\limits_{q<\infty}
\sobolev{m-1,p}_{\text{\scriptsize\rm loc}}(\Omega,\mathbb{H})$.
Therefore Lemma~\ref{local minimizer} implies that local minimizers
of the classes~(ii)--(iii) belong to
$\bigcap\limits_{n\in\mathbb{N},q<\infty}
\sobolev{m,p}_{\text{\scriptsize\rm loc}}(\Omega)$.
With the help of the reality condition for immersion
into the pure imaginary quaternions $\simeq\mathbb{R}^3$
these arguments carry over to case~(i).
\end{proof}

We do not claim that the minimizers are realized by immersions.
They may have branch points. In general they may be compositions
of a finite--sheeted branched covering together with an immersion.
Finally, we remark that the existence of minimizers was proven by
Simon \cite{Si1,Si2} on the class of all smooth immersion
from a compact orientable surface of genus one into the Euclidean spaces
$\mathbb{R}^n$ ($n\geq 3$). Furthermore, Bauer and Kuwert \cite{BK}
extended these arguments to the classes of all smooth immersions
from compact orientable surfaces into the Euclidean spaces
$\mathbb{R}^n$ ($n\geq 3$). It might be possible to deduce
these results for $n=3$ and $n=4$ from our results. In fact, since
the stereographic projections of the minimal surfaces in $S^3$
constructed by Lawson \cite{Law} have Willmore functionals
less than $8\pi$ (compare \cite{Si1,Si2}),
it would suffices to prove that at the boundary of the moduli spaces
$\moduli_g$, which contains stable curves with ordinary double points,
the Willmore functional is at least equal to $8\pi$.
This would follow from \cite{LY},
if the corresponding conformal mappings preserve these double points.
Moreover, with the help of \cite{Kon3} our results might be generalized
to conformal mappings into higher--dimensional Euclidean spaces.

\end{document}